\documentclass[american]{amsart}
\pdfoutput=1

\usepackage[american,english]{babel}
\usepackage{amsmath,amssymb,amsthm}
\usepackage{mathtools}
\usepackage[latin1]{inputenc}
\usepackage[T1]{fontenc}
\usepackage{tikz}
\usepackage{graphicx}
\usepackage{url}
\usepackage{pgfplots}
\usepackage{calc}
\usepackage{intcalc}
\usetikzlibrary{fadings}
\usepackage{mathdots}
\usepackage{ifthen}
\usepackage[all]{xy}
\SelectTips{cm}{10}
\usepackage{color}

\newenvironment{overblik}{\noindent\textbf{Outline of proof of Theorem 1:}}{}

\DeclareMathOperator{\ev}{ev}

\DeclareMathOperator{\dist}{dist}
\DeclareMathOperator{\Th}{Th}

\DeclareMathOperator{\id}{Id}

\DeclareMathOperator{\im}{Im}

\DeclareMathOperator{\re}{Re}
\DeclareMathOperator{\Real}{Re}
\DeclareMathOperator{\Imag}{Im}

\DeclareMathOperator{\inte}{int}

\DeclareMathOperator{\Gl}{Gl}

\DeclareMathOperator{\CW}{CW}
\DeclareMathOperator*{\hocolim}{Hocolim}
\DeclareMathOperator*{\colim}{Colim}

\newcommand{\Or}{\mathrm{O}}

\newcommand{\U}{\mathrm{U}}
\newcommand{\BO}{B\mathrm{O}}
\renewcommand{\epsilon}{\varepsilon}
\newcommand{\arq}{{\vec{q}}}
\newcommand{\arp}{{\vec{p}}}
\newcommand{\arz}{{\vec{z}}}

\newcommand{\arzm}{\vec{z_-}}
\newcommand{\arwm}{\vec{w_-}}

\newcommand{\arzet}{\vec{z_1}}
\newcommand{\arzto}{\vec{z_2}}
\newcommand{\arw}{{\vec{w}}}
\newcommand{\tp}{{\widetilde{p}}}
\renewcommand{\th}{{\widetilde{h}}}
\newcommand{\tH}{{\widetilde{H}}}
\newcommand{\tC}{{\widetilde{C}}}
\newcommand{\tf}{{\widetilde{f}}}
\newcommand{\te}{{\widetilde{\epsilon}}}
\newcommand{\unp}{\underline{P}}
\newcommand{\dd}{d}
\newcommand{\pdd}[1]{\tfrac{\partial}{\partial #1}}
\newcommand{\R}{\mathbb{R}}
\newcommand{\N}{\mathbb{N}}
\newcommand{\ovl}[1]{\overline{#1}}
\newcommand{\C}{\mathbb{C}}
\newcommand{\Z}{\mathbb{Z}}
\newcommand{\La}{\mathcal{L}}
\newcommand{\SLa}{\Gamma}
\newcommand{\Sla}{\SLa}
\newcommand{\LUO}[1]{\U#1/\Or#1}

\newcommand{\param}{{u}}

\newcommand{\Lre}{\Lambda_{r}}

\newcommand{\Lrb}{\Lambda_r^{e<\beta}}

\newcommand{\efed}{\tfrac{1}{5}}
\newcommand{\efd}{\tfrac{1}{4}}
\newcommand{\eh}{\tfrac{1}{2}}
\newcommand{\tfd}{\tfrac{3}{4}}
\newcommand{\ffd}{\tfrac{4}{5}}

\newcommand{\sfd}{\tfrac{6}{5}}

\DeclarePairedDelimiter\absv{\lvert}{\rvert}
\DeclarePairedDelimiter\norm{\lVert}{\rVert}
\DeclarePairedDelimiter\inner{\langle}{\rangle}
\newtheorem{Theorem}{Theorem}
\newtheorem{Proposition}{Proposition}[section]
\newtheorem{Lemma}[Proposition]{Lemma}
\newtheorem{Corollary}[Proposition]{Corollary}

\theoremstyle{remark}
\newtheorem{Remark}[Proposition]{Remark}
\newtheorem{ex}[Proposition]{Example}

\theoremstyle{definition}
\newtheorem{Definition}[Proposition]{Definition}

\DeclarePairedDelimiter\pare{\lparen}{\rparen}
\DeclarePairedDelimiter\brac{\lbrack}{\rbrack}

\newcommand{\no}[1]{\norm{#1}_{C^2}}
\newcommand{\nor}[1]{\norm{#1}_{C^1}}

\usepackage{hyperref}

\def\co{\colon\thinspace}
\begin{document}



\title{The Viterbo Transfer as a Map of Spectra}

\author[T. Kragh]{Thomas Kragh}
\thanks{The author was partially funded by CTQM, OP-ALG-TOP-GEO, Topology in Norway, Aarhus University, Carlsberg, MIT, Oslo University, and Uppsala University.}

\begin{abstract}
Let $L$ and $N$ be two smooth manifolds of the same dimension. Let $j\colon L\to T^*N$ be an exact Lagrange embedding. We denote the free loop space of $X$ by $\Lambda X$. In \cite{MR1617648}, C. Viterbo constructed a transfer map $(\Lambda j)^! \colon H^*(\Lambda L) \to H^*(\Lambda N)$. This transfer was constructed using finite dimensional approximation of Floer homology. In this paper we define a family of finite dimensional approximations and realize this transfer as a map of Thom spectra: $(\Lambda j)_! \colon (\Lambda N)^{-TN} \to (\Lambda L)^{-TL+\eta}$, where $\eta$ is a virtual vector bundle classified by the tangential information of $j$.
\end{abstract}


\maketitle
\setcounter{tocdepth}{1}

\section{Introduction and Statement of Results}

Let $N$ be a closed $d$-dimensional smooth manifold, and let $\pi \colon T^*N \to N$ be the projection from the cotangent bundle of $N$ to $N$. The Liouville form (or canonical 1-form) $\lambda$ is defined by
\begin{align*}
  \lambda_{q,p} (v) = p(\pi_*(v)), \quad q \in N, p\in T^*_qN, v \in T_{q,p}(T^*N).
\end{align*}
The 2-form $\omega=-d\lambda$ is non-degenerate and thus defines a canonical symplectic structure on $T^*N$. Let $L$ be another closed $d$-dimensional smooth manifold. An embedding $j\colon L \to T^*N$ is called Lagrangian if $j^*\omega=-dj^*\lambda=0$ and called exact Lagrangian if $j^*\lambda$ is exact. We assume from now on that $j$ is an exact Lagrangian embedding. The trivial examples of such embeddings are those which are Hamiltonian isotopic to the zero section. To this day no non-trivial examples have been found, and the nearby Lagrangian conjecture states that there are no others. This is trivially true for $N=S^1$.

Recently there has been much progress in this area. Specifically in Nadler \cite{Nadler} and Fukaya, Seidel and Smith \cite{FSS} it is proven, independently, that under certain conditions $j$ is a homology equivalence. This has been extended by Abouzaid in \cite{Abou1} to prove that when the Maslov index vanishes then $j$ is a homotopy equivalence. Finally in \cite{MySympfib} we use results from this paper and some new methods to prove homology equivalence without any assumptions, and with Abouzaid we prove in general that $j$ is a homotopy equivalence. Furthermore, restrictions on the smooth structures and immersions classes has been found in certain cases (mostly spheres) in \cite{MR2874640}, \cite{smooth2} and \cite{immersions}. The latter uses results from this paper. Finally combining the homotopy equivalence result with the dimension dependent argument by Hind in \cite{MR2904911} proves that for $N=S^2$ or $N=\R P^2$ the exact Lagrangian $L$ is in fact Hamiltonian isotopic to the zero-section, thereby confirming the conjecture for these $N$.

We denote the free loop space of a space $X$ by $\Lambda X$. In \cite{MR1617648}, Viterbo constructs a transfer map $(\Lambda j)^!$ on cohomology, such that
\begin{align*}
  \xymatrix{
    H^*(\Lambda L) \ar[r]^{(\Lambda j)^!} \ar[d]^{i^*} & H^*(\Lambda N) \ar[d]^{i^*} \\
    H^*(L) \ar[u]^{\ev_0^*} \ar[r]^{(\pi\circ j)^!} & H^*(M) \ar[u]^{\ev_0^*}
  }
\end{align*}
commutes. Here $(\pi\circ j)^!$ is the standard transfer map on cohomology, $\ev_0$ is the evaluation at base point, and $i$ is the inclusion of constant loops. In this paper we call this map the Viterbo transfer. Viterbo used this transfer as obstruction to the existence of exact Lagrangian embeddings.

Because $j\colon L\to T^*N$ is Lagrangian we get a Maslov class in $H^1(L)$. This defines a map $\Lambda L \to \Z$ called the Maslov index, and it turns out that the Viterbo transfer is shifted in grading on each component by this Maslov index. In this paper we prove the following theorem, which explains this grading shift (when using Thom isomorphism). 

\begin{Theorem} \label{thm:2}
  The Viterbo transfer can be realized as a map of Thom-spectra such that the diagram
  \begin{align*}
    \xymatrix{
      (\Lambda N)^{-TN} \ar[r]^{\Lambda j_!} & (\Lambda L)^{-TL+\eta}   \\
      N^{-TN} \ar[u] \ar[r]^{j_!} & L^{-TL} \ar[u]
    }
  \end{align*}
  commutes. Here $j_!\colon N^{-TN} \to L^{-TL}$ is the usual transfer for manifolds (defined on Thom-spectra) and $\eta$ is a virtual vector bundle classified by the tangential information of the embedding $j\colon L \to T^*N$, and the local dimension of $\eta$ is the Maslov index.

  Furthermore, the definition of the spectra and the map is a contractible choice, and the identification of the homotopy type of all except the top right is also a contractible choice.
\end{Theorem}
The reason for the non-canonicality of the last spectrum is due to a choice of homotopy of Lagrangians along $L$, and if the strong version of the nearby Lagrangian is true (the space of exact Lagrangians are contractible) then this choice will in fact also be canonical.

\begin{Remark}
  \label{rem:Intro:1}
  When defining the standard transfer for manifolds one can alternatively describe this as a map $j_!'\co \Sigma^\infty N_+ \to L^{TN-TL}$, which in the case of non-orientable manifolds is different on homology. Similarly one can make some alternative choices in the construction of this transfer map, and the content of Corollary~\ref{cor:Viterbo:1} and Corollary~\ref{cor:Maslov:3} is that the alternative diagram becomes:
  \begin{align*}
    \xymatrix{
      \Sigma^\infty(\Lambda N)_+ \ar[r]^{\Lambda j'_!} & (\Lambda L)^{TN-TL+\eta}   \\
      \Sigma^\infty N_+ \ar[u] \ar[r]^{j'_!} & L^{TN-TL} \ar[u]
    }
  \end{align*}
\end{Remark}

We have included a short discussion about spectra and CW spectra in Appendix~\ref{cha:app}, and since every spectrum appearing in this paper is homotopy equivalent to a CW spectrum as defined in \cite{MR0402720} we refer to this book for a more thorough introduction to the concept of spectra. However, note that the categories of spectra that are usually used today are much more structured and thus handy for a lot of things. In particular they are symmetric monoidal categories with respect to smash product \emph{before} passing to the homotopy category. This is convenient when considering things like products (ring spectra), but we will not do that here (even though there are natural ring-structures on all the spectra above, and one could conjecture $\Lambda j_!$ to be a ring-spectrum map).

In the original construction by Viterbo, the Thom isomorphism is used on what turns out to be the virtual vector bundle $-TL+\eta$ (a topological $K$-theory class). However, $\eta$ is not necessarily oriented, but if we assume $(\pi \circ j) \colon L \to N$ to be relatively oriented and relative spin, it will be. This has recently yielded a new insight into coherent orientations. Se \cite{MySympfib}, \cite{AbouzaidSymp}, and \cite{Abou2} for more details on this.

The rest of the introduction is an overview of the construction of $(\Lambda j)_!$ using finite dimensional approximations of Floer homology.

\begin{overblik}
  The actual construction does not use Floer homology. However it is very illuminating to sketch the relation. This relation also justifies considering the spectra constructed as representations of the stable homotopy type of Floer homology in cotangent bundles - at least in the oriented and spin case. \nocite{MR1035373}

  For any Hamiltonian $H\colon T^*N \to \R$ we define the action integral
  \begin{align*}
    A_H\colon \Lambda T^*N \to \R
  \end{align*}
  by the formula
  \begin{align*}
    A_H(\gamma) = \int_\gamma \lambda -Hdt.
  \end{align*}
  We give $N$ a Riemannian structure (and wait til the very end of the paper to argue that the Viterbo transfer do not depend on this choice), and we will always assume that $H(q,p)=\mu\norm{p}+c$ for large $\norm{p}$, where $\mu\in\R$ is not the length of any closed geodesic. We say that $H$ is linear at infinity. Floer homology $FH_*(T^*N,H)$ is essentially Morse homology of $A_H$ perturbed on the infinite dimensional manifold $\Lambda T^*N$ (see e.g. \cite{MR2190223}). When $N$ is orientable and spin the linear at infinity case can be calculated to satisfy (see \cite{MR1726235} and \cite{MR1617648} - where the need for the spin assumption was overlooked)
  \begin{align} \label{Floerhom}
    FH_*(T^*N,H) \approx H_*(\Lambda ^\mu N).
  \end{align}
  Here $\Lambda^\mu N$ denotes the space of loops with length less than $\mu$. Define $T^* \Lambda_r N$ as the cotangent space of the manifold of $r$-pieced geodesics each of length less than some fixed $\delta_0>0$. The finite dimensional approximations we define in Section~\ref{Florlike} can be described in the following way: for large $r$  we define embeddings:
  \begin{align*}
    i_r \colon T^*\Lambda_r N \to \Lambda T^*N,
  \end{align*}
  where $T^*\Lambda_r N$ is a finite dimensional manifold, and these satisfy
  \begin{itemize}
  \item{The image of $i_r$ contains all critical points of $A_H$}
  \item{The composition $S_r=A_H \circ i_r$ has no other critical points than those from $A_H$}
    \item{There is a ``consistent'' way of applying Morse theory for $S_r$ on $T^*\Lambda_r N$ and creating a space $Z'$ such that the homotopy type of this space does not change under small compact perturbations of $H$.}
    \item In fact, the Morse theory of $A_H \circ i_r$ captures all of the Morse homology of $A_H$.
  \end{itemize}
  We will not prove the last point. In the actual construction, we use the theory of Conley indices described in section~\ref{homoindex}, but for the purpose of this overview, one may think of $S_r$ as a Morse function, and thus think of $Z'$ as a cell complex with; one cell per critical point, and an extra base-point cell because flow lines can go to $-\infty$.

  In section~\ref{Florlike} we explicitly define a function $S_r$ as above (but skipping the $i_r$ and simply writing down an formula for $S_r$) and in Section~\ref{sec:gener-finite-dimens} prove that $Z'$ (in the case described above) is homotopy equivalent to the Thom space
  \begin{align*}
    \Th(T\Lambda^\mu_r N) \simeq (\Lambda_r^\mu N)^{T\Lambda_r^\mu N} = D(T\Lambda^\mu_r N) / U(T\Lambda^\mu_r N),
  \end{align*}
  where $\Lambda^{\mu}_r N$ is the manifold of piecewise geodesics, with $r$ pieces each having length less than $\mu/r$, $D(\cdot)$ denotes the unit disc bundle, and $U(\cdot)$ denotes the unit sphere bundle. This was already proven by Viterbo in \cite{MR1617648}, but because we have an explicitly defined $S_r$, we can prove it more directly; and we will need this more direct approach to identify the homotopy types later. Using the Thom isomorphism, this is consistent with Equation~\eqref{Floerhom}, but be warned: As interpreted by the new insight into coherent orientations the homology of this is \emph{not} always the Floer homology unless $N$ is oriented and spin. However all of these differences are mostly due to coherent orientations issues and not really important for the heuristical idea.

  This formula suggests that when increasing $r$ the space changes by a (relative) Thom construction using the tangent bundle $TN$, and this is precisely what we prove in Section~\ref{cha:suspmaps}. Since spectra are defined by sequences of spaces up to standard reduced suspensions this does not precisely define a spectrum. So, in Section~\ref{sec:gener-funct-spectr} we describe how to untwist these copies of $TN$ by adding copies of the normal bundle. Subsequently, defining a spectrum out of the collection of Conley indices for all large $r$. We denote this spectrum by $Z^\mu$, but the reader not to comfortable with spectra can continue to consider this as a CW complex with 1 cell per critical point (and a base-point corresponding to $-\infty$).  We will be taking the limit $\mu \to \infty$ and we will denote the limit of these spectra $Z=\lim_{\mu \to \infty} Z^\mu$.

  This construction of a spectrum out of the Morse theory of $S_r$ (or heuristically $A_H$) is completely canonical and natural. Indeed, we argue that all relevant choices leads to spectra with contractible choice of homotopy equivalences between them, and we will see that the natural quotients and inclusions on Conley indices induce natural maps of spectra. This and the usual construction of the Viterbo transfer map (Viterbo functoriality) gives the map $(\Lambda j)_!$ of spectra in Theorem~\ref{thm:2}. However, for the reader unfamiliar with this construction we quickly outline the idea. The full construction is done in Section~\ref{vitconst}. 

  Because $j\colon L \to T^*N$ is a Lagrangian embedding, we can use the Darboux-Weinstein theorem to extend $j$ to a symplectic embedding of a small neighborhood of the zero section in $T^*L$ (which we can assume is $DT^*L$ by choosing the Riemannian structure on $L$ appropriately). Using this neighborhood and the fact that $j$ is exact, we can adjust $H$ such that all the critical points of $S_r$ with critical values above $c$, for some $c$, are loops inside the neighborhood of $L$.
  
  In fact, close to $L$ the Hamiltonian $H$ is defined to be a smooth approximation of the norm function $\norm{p_{L}}$ (using some Riemannian structure on $L$) times some constant ${\mu_{L}}$. We now define $W$ to be the spectrum (for each $\mu$ and $\mu_L$ and then taking the limit as both goes to $\infty$) defined by the usual quotients on Conley indices, which in CW language means that we collapse the sub-complex $Y$ defined by those cells associated to critical points with critical value less than $c$. Because $H$ are close to $\mu_{L}\norm{p_{L}}$ in a neighborhood of $L$ we see that the space $W$ is highly related to the linear at infinity case on $L$.

  In Section~\ref{indexcalc} we identify the homotopy type of the source spectrum (canonically) as $(\Lambda N)^{-TN}$.

  Section~\ref{sec:gener-finite-dimens} and Section~\ref{loclinsec} proves a localization result that makes it possible to also identify the homotopy type of the target (this is also used in Section~\ref{indexcalc}, but a much less general statement is needed for that part). This requires some work. Indeed, when constructing the finite dimensional approximations in $T^*N$ we use the cotangent bundle structure, and it stands to reason that even though we are using a Hamiltonian, which close to $L$ describes a well-known (and similar to the case of $N$) Hamiltonian system, the resulting homotopy type of $W$ could depend on the structure on $T^* N$. This is, indeed, the case.

  Let $W' \simeq (\Lambda L)^{-TL}$ denote the spectrum we get from using the usual structures on $T^*L$ to define a spectrum out of this Hamiltonian system. This has this homotopy because it is the case we computed in Section~\ref{indexcalc}, but with $L$ replacing $N$.

  The most important structure (implicitly) used in the definition of $W$ is the fact that at each point in $T^*N$ we have a Lagrangian subspace defined by vertical vectors. When restricting this to the neighborhood of $L$ this may differ from the Lagrangian subspaces given by vertical vectors in $T^*L$. At this point Viterbo uses a classification result on generating functions to describe the difference of the two spaces (spectra in our case) $W$ and $W'$ as a relative Thom space construction. This is very subtle, and in this paper we instead use the very explicit constructions to actually calculate the stable homotopy type of $W$. This very explicit construction is what led to the new insight into coherent orientations of Floer homology in cotangent bundles mentioned above.

  The calculation involves defining a family of finite dimensional approximations $S_r^{\SLa}$ depending on $\SLa$ which is a section in the bundle.
  \begin{align*}
    \La(T(DT^*L)) \to DT^*L.
  \end{align*}
  Here $T(DT^*L) \to DT^*L$ is viewed as a symplectic vector bundle and $\La(T(T^*L))\to T^*L$ is the associated fibration of Lagrangian Grassmannians, i.e. the fiber of the above fibration is $\La(n) \simeq \LUO{(n)}$, which is the Grassmannian of linear Lagrangians subspaces in $\C^n$.

  In the neighborhood $DT^*L\subset T^*N$ of the zero section of $L$, we have the two canonical sections of this bundle: $\SLa^L$ given by the section in $\La(T(T^*L))$ which to a point associates the vertical directions w.r. to $L$, and similarly $\SLa^N$ (restricted to the neighborhood $DT^*L$). The construction of $S_r^{\SLa}$ implies by homotopy invariance that the spectrum $W$ does not change when perturbing $\SLa$. This implies that if $\SLa^L$ and $\SLa^N$ were homotopic, we would in fact get that $W$ and $W'$ are homotopy equivalent. However as mentioned before this is not the case in general\footnote{unless the nearby Lagrangian conjecture is true, but we can of course not assume that.}. If they are not homotopic we may stabilize by adding trivial factors and get ``stabilized'' finite dimensional approximations
  \begin{align*}
    S_r^{\SLa^N \oplus \R^k} \colon T^*\Lambda_r(N\times \R^k) \cong T^*\Lambda_r N \times (\R^{2kr},\omega_0) \to \R.
  \end{align*}
  We have set up the grading such that the spectrum we get if we apply Morse theory (restricting to values above the constant $c$ from above) in this stabilized case is again $W$. Indeed, we will argue that the Conley indices fore each $r$ is simply a reduced suspension of the old, and by definition the grading of each is shifted in the natural way to compensate. The fact that the Conley indices are reduced suspensions is due to the fact that this function is just a constant quadratic term in the second variable $\R^{2kr}$. It is now a homotopy theoretical fact that: although $\SLa^L$ and $\SLa^N$ are not homotopic we can find a homotopy from $\SLa^N \oplus \R^k$ to $\SLa^L \oplus \ovl{F}$, where $\ovl{F} \colon DT^*L \to \La(k)$ is some smooth map. It is still true that $S_r^{\SLa^L \oplus \ovl{F}}$ is quadratic in the second variable, but this quadratic form is no longer the same for different values of the first variable. Indeed, $\ovl{F}$ is a now a Lagrangian in the second factor depending on the point in the first factor. In Section~\ref{sec:quadr-forms-assoc} we calculate that for $r$ odd the negative eigenbundle of the quadratic form is given by a specific classifying map from $\Lambda L \to \Z\times BO$, which in Section~\ref{maslov} leads to the definition of the Maslov bundle $\eta$, which is classified by the map
  \begin{align*}
    \Lambda L \xrightarrow{\ovl{F}} \Lambda \La(k) \to \Lambda \La \to \Omega \La \simeq \Z \times B\Or.
  \end{align*}
  This map is described in more detail in Section~\ref{sec:quadr-forms-assoc} and Section~\ref{maslov}.
  
  In the end this shows that the spectrum $W'$ is a relative Thom space constructions on the pairs defining $W$, which implies Theorem~\ref{thm:2}.
\end{overblik}

\textbf{Acknowledgments.} I would like to thank Marcel B\"okstedt, John Rognes and Mohammed Abouzaid for many enlightening conversations about material related to this paper.


\section{Conley Indices and Canonicality} 
\label{homoindex}

In this section we introduce the notion of Conley indices (from \cite{MR511133}). To be able to discus how canonical the spectra we define are we recall some proofs. To make things easier we define the notion of a good index pair, for which we prove that the Conley indices are preserved under perturbation.

\subsection{Definitions}

Let $M$ be a smooth open manifold, and let $f\colon M \to \R$ be a smooth function. A pseudo-gradient $X$ for $f$ is a smooth vector field on $M$ such that the directional derivative $X(f)$ is positive at non-critical points and $X=0$ at critical points. The choice of a pseudo-gradient is a contractible choice since a convex combination of pseudo-gradients is a pseudo-gradient, and they exist since a gradient is a pseudo-gradient.

We will denote the flow of $-X$ by $\psi_t$. Let $a<b$ be regular values of $f$ which are isolated from the critical values of $f$. We wish to define the \emph{Conley index} $I_a^b(f,X)$ (when possible).
  
An index pair $(A,B)$ for $(f,X)$ with respect $a<b$ is a pair of subspaces of $M$ satisfying the following properties
\begin{itemize}
\item[I1:]{$B\subset A\subset f^{-1}([a,b])$.}
\item[I2:]{$A$ and $B$ are compact.}
\item[I3:]{$\inte(A-B)$ contains all critical points of $f$ with critical
    values in $]a,b[$.}
\item[I4:]{For each $x\in A$ the pair of spaces
    \begin{align*}
     (\{t\geq 0 \mid \psi_t(x) \in A\},\{t\geq 0 \mid \psi_t(x) \in B\})
    \end{align*}
    is either $(\{t\geq 0\},\emptyset)$ or a pair of closed intervals with the same maximum.
}
\end{itemize}
Notice that I4 means that either; a flow line stays in $A-B$ converging to a critical point or it exits $A$ \emph{through} $B$. The set $B$ is called the \emph{exit set}. When such index pairs exist we define the Conley index
\begin{align} \label{eq:6}
  I_a^b(f,X)=A\big/ B.
\end{align}
If $X=\nabla f$ we write $I_a^b(f)$. If all critical values of $f$ is contained in an interval $]a,b[$ we simply write $I(f,X)$ instead of $I_{a}^{b}(f,X)$ and call this the \emph{total} index.

Since the choice of a pseudo-gradient which admits such index pairs is \emph{not} a contractible choice we fix a given pseudo-gradient $X$ in the rest of this section. The following lemma was due to Conley, however, we recall the proof from \cite{MR1045282}.

\begin{Lemma} \label{lem:indwell}
  When index pairs exist the based homotopy type of $I_a^b(f,X)$ is well-defined, and the homotopy equivalences is induced by the flow $\psi_t$ and hence a contractible choice.
\end{Lemma}

\begin{proof}
  Assume we are given two index pairs $(A_i,B_i)$, $i=1,2$. By definition $\inte(A_2-B_2)$ contains the closed (hence compact) image set of flow lines of $\psi_t$ which converges to critical points at both ends with value in $]a,b[$. This implies that for any $x \in \ovl{(A_1-B_1)}-\inte(A_2-B_2)$ there must be a neighborhood $U$ of $x$ and a $t\geq 0$ such that; for any $y\in U$ we have that one of the two points $\psi_{\pm t}(y)$ are not in $A_1-B_1$ (or not defined). Mutatis mutantis for $x \in \ovl{(A_2-B_2)}-\inte(A_1-B_1)$.

  By compactness of these sets there is a minimal $t_0=t_0(\ovl{A_1-B_1},\ovl{A_2-B_2})\geq 0$ such that when $t > t_0$ we get that; if $\psi_{[-t,t]}(x) \subset A_1-B_1$ then $x \in A_2-B_2$, and if $\psi_{[-t,t]}(x) \subset A_2-B_2$ then $x \in A_1-B_1$.

  For $t > 3t_0$ we now define
  \begin{align*}
    h^{12}_t \co A_1/B_1 \to A_2/B_2
  \end{align*}
  by the map induced by $\psi_t$ when
  \begin{align}\label{eq:13}
    \psi_{[0,2t/3]}(x) \subset A_1-B_1 \quad \textrm{and} \quad \psi_{[t/3,t]}(x) \subset A_2-B_2
  \end{align}
  and sending everything else to the base point $[B_2]$.

  We claim that $h^{12}_t([x])$ is continuous in $([x],t)\in A_1/B_1 \times (3t_0,\infty)$. To see this we first prove that Equation~\eqref{eq:13} is an open condition on the set $(x,t) \in A_1 \times (t_0,\infty)$. Indeed, Since the flow must exit $A_1$ through $B_1$ the first half is by compactness of $B_1$ an open condition, and given this condition we get from the assumptions on $t$ that
  \begin{align*}
    \psi_{[0,2t/3]}(x) \subset A_1-B_1 \quad \Rightarrow \quad \psi_{t/3}(x) \in A_2-B_2.
  \end{align*}
  It again follows that since the flow must exit $A_2$ through $B_2$ that the latter condition is again an open condition inside the set of $x\in A_1$ satisfying the first condition. So, to prove continuity we need only consider an arbitrary sequence $(x_n,t_n) \in (A_1-B_1) \times (3t_0,\infty)$ which all satisfies Equation~\eqref{eq:13}, but the limit $(x,t)=\lim_{n\to \infty}(x_n,t_n) \in A_1\times (3t_0,\infty)$ does not. It follows by compactness of the sets and continuity of the flow that either $\psi_t(x) \in B_2$ or $\psi_{2t/3}(x) \in B_1$. We are finished if $\psi_t(x) \in B_2$ so assume for contradiction that $\psi_t(x) \notin B_2$. This, by the above assumptions on $t>3t_0$, means that
  \begin{align*}
    \psi_{[t/3,t]}(x) \subset A_2-B_2 \quad \Rightarrow \quad \psi_{2t/3}(x) \in A_1-B_1 \qquad \Rightarrow \quad \psi_{2t/3}(x) \notin B_1
  \end{align*}
  which provides the contradiction.

  Since the $t_0$'s satisfy
  \begin{align*}
    t_0(\ovl{A_1-B_1},\ovl{A_2-B_2}) +  t_0(\ovl{A_2-B_2},\ovl{A_3-B_3}) \geq     t_0(\ovl{A_1-B_1},\ovl{A_3-B_3}),
  \end{align*}
  we see that these maps behave well under composition, and we also notice that for $(A_1,B_1)=(A_2,B_2)$ the maps are defined for all $t\geq 0$ and $h^{12}_0$ is the identity.
\end{proof}

\subsection{Quotients and Inclusions}

Some very important aspects of Conley indices are the natural inclusion and quotient maps
\begin{align} \label{eq:70} 
  i \colon I_a^b(f,X) \to I_a^c(f,X) \notag \\
  q \colon I_a^c(f,X) \to I_b^c(f,X),
\end{align}
where $a<b<c$ are all regular for $f$. These maps are constructed as follows. Let $(A,B)$ be an index pair for $I_a^c(f,X)$.
\begin{itemize}
\item The pair $(A\cap f^{-1}([a,b]) , B \cap f^{-1}([a,b]))$ is an index pair for $I_a^b(f,X)$, and the map $i$ is induced by the inclusion of this pair into $(A,B)$.
\item The pair $(A\cap f^{-1}([b,c]),\brac{A\cap f^{-1}(\{b\})}\cup \brac{B\cap f^{-1}([b,c])})$ is an index pair for $I_b^c(f,X)$, and the map $q$ is the map from $A/B$ collapsing the subset $(A\cap f^{-1}([a,b]))/B$.
\end{itemize}
These maps commute on the nose with the homotopy equivalences in Lemma~\ref{lem:indwell}.

\subsection{Good Index Paris and Homotopy Invariance}

It is very convenient to introduce the concept of a good index pair $(A,B)$.

\begin{Definition}
  An index pair $(A,B)$ (for $(f,X)$ with respect to $a<b$) is called \emph{good} if $B \subset f^{-1}(a)$ and for any vector field $X'$ on $M$ sufficiently close to $X$ on $A$ we have; any point $x\in A-B$ will under the flow of $-X'$ for a short positive time stay in $A$.
\end{Definition}

Note that this implies that the flow still exists through $B$ and is thus similar to I4. However, we say nothing about the flow not reentering $A$. Indeed, this would be unreasonable since we only ask that $X'$ is close to $X$ on $A$. In the language of isolated invariant sets (cf \cite{MR1045282}) this definition assures that the isolated invariant set of $(f,X)$ associated to the index pair $(A,B)$ stays within $A$ under small perturbations of $(f,X)$.

There is another (more global) reason why good index pairs are convenient. Indeed, the following lemma is not true if the word ``good'' is removed.

\begin{Lemma}
  \label{lem:Homoindex:2}
  Let $M' \subset M$ be an open submanifold. If all critical points of $f$ lie in $M'$ and $(f',X')=(f_{\mid M'},X_{\mid M'})$ has a good index pair then this is also a good index pair for $(f,X)$. It thus follows that we can canonically identify
  \begin{align*}
    I_a^b(f,X) = I_a^b(f',X').
  \end{align*}
\end{Lemma}

\begin{proof}
  I1 through I3 above is trivial. I4 follows since $B\subset f^{-1}(a)$ and hence any flow line exiting $A$ cannot return since the value of $f$ has gotten to low.
\end{proof}

A small detail, which could be avoided in a different way is also good about good index pairs.

\begin{Lemma}
  \label{lem:Homoindex:3}
  A good index pair is a cofibrant pair.
\end{Lemma}

\begin{proof}
  By using the flow of the negative pseudo-gradient, but stopping it when $f=a$ we get a deformation retraction of $A \cap f^{-1}([a,a+\epsilon[)$ (which is a neighborhood of $B$) onto $B$.
\end{proof}

This has the immediate consequence.

\begin{Corollary}
  \label{cor:Homoindex:1}
  The Conley index $I_a^b(f,X)$ is well-based, and has the homotopy type of a CW complex.
\end{Corollary}

Now let $M \xrightarrow{\pi} I$ be a submersion with each fiber $M^s=\pi^{-1}(s)$ a smooth manifold without boundary (e.g. $M=N\times I$). Let $f\co M \to \R$ be a smooth map, and denote the restriction to $M^s$ by $f^s$. Let $X$ be a vertical vector field on $M$, i.e. it restricts to vector fields $X^s$ on $M^s$ for each $s\in I$. Assume that for each $s\in I$ the vector field $X^s$ is a pseudo-gradient for $f^s$, and that $a<b$ are regular for all $f^s$. Also assume that
\begin{itemize}
\item the union over $s\in I$ of the critical points of $f^s$ (in $M$) is compact and
\item for each $s\in I$ there exist a good index pair $(A^{s},B^{s})$ defining $I_a^b(f^{s},X^{s})$.
\end{itemize}

\begin{Lemma} \label{hominv}
  Under the above assumptions we have
  \begin{align*}
    I_a^b(f^0,X^0) \simeq I_a^b(f^1,X^1).
  \end{align*}
  Furthermore, this homotopy equivalence is a contractible choice and naturally commutes with the homotopy equivalences from Lemma~\ref{lem:indwell} and the quotients and inclusions above.
\end{Lemma}

\begin{proof}
  Given $s_0$, we would like to prove that a good index pair $(A,B)=(A^{s_0},B^{s_0})$ is an index pair defining $I_a^b(f^s,X^s)$ when $s$ is sufficiently close to $s_0$. However this is not exactly possible because we cannot be certain that I1 is satisfied, but because $a$ and $b$ are regular values and we know that the critical points form a compact set, we can replace $a$ and $b$ by $a-\delta$ and $b+\delta$ for some small $\delta$ without changing the Conley indices. Now I1 is not a problem for $s$ close to $s_0$.

  I2 is obvious and I3 follows by compactness of the union of the critical sets. Since the good pair assumption makes sure that we only exit $A$ through $B$ (and transversely so), we  are only left with proving that for any point $x$ in $B$ we have $\{t\geq 0\mid \psi_t^{s}(x)\in A\}=\{0\}$ and I4 will follow. This is equivalent to proving that the flow does not return to $A$ after exiting through $B$.

  Let $(\psi_t)^s$ be the flow of $-X^{s}$ for time $t$. Since $-X^{s_0}(f^{s_0})$ restricted to $B$ is negative, the same is true for $-X^s(f^{s_0})$ for $s$ close to $s_0$. So fix a $\tau>0$ such that $-X^s_x(f^{s_0})$ is negative for $s\in [s_0-\tau,s_0+\tau]$ and $x\in B$. Then (by compactness of $B$) we may choose a $\delta>0$ such that $f^{s_0}((\psi_t)^s(x))$ is strictly decreasing for $t\in[0,\delta], x\in B$, and $s\in [s_0-\tau,s_0+\tau]$. Because $f^{s_0}(B)=\{a\}$ and $f^{s_0}(A)\subset [a,b]$ we see that: for $x\in B$, $t\in (0,\delta]$ and $s\in [s_0-\tau,s_0+\tau]$ we have that $(\psi_t)^s(x)$ is not in $A$. Furthermore, by compactness of $B$ there is an $\epsilon>0$ such that $f^{s_0}((\psi_{\delta})^{s}(B))<a-\epsilon$ when $s\in [s_0-\tau,s_0+\tau]$. By continuity of the family $f^s$ we may now pick a $\tau'>0$ smaller than $\tau$ such that for $s\in [s_0-\tau',s_0+\tau']$ we have that $f^s((\psi_{\delta})^{s}(B)) <a-\epsilon/2$. We can similarly assume (by compactness of $A$) that for small $\tau'$ we have that $f^s(A)>a-\epsilon/3$. So, we conclude that the flow of any $x\in B$ using $-X^s$ for $s\in [s_0-\tau',s_0+\tau']$ immediately exits $A$ and stays out for time $t\in (0,\delta]$, and at time $\delta$ the value of $f^s$ is less than $f^s$ is on all of $A$. So, since $X^s$ is a pseudo-gradient for $f^s$ the flow of $-X^s$ will never reenter $A$.

  To see how this defines a contractible choice of homotopy equivalences from $I_a^b(f^0,X^0)$ to $I_a^b(f^1,X^1)$ we consider finite coverings of $I$ by open intervals $J_\alpha$ over which we have chosen a common index pair $(A_\alpha,B_\alpha)$ for $(f^s,X^s)$ with $s \in J_\alpha$. Such exists because of the above and compactness of $I$. Now choose a subdivision
  \begin{align*}
    0 = s_0 < s_1 < \cdots < s_k = s_{k+1} = 1
  \end{align*}
  such that each closed interval $[s_i,s_{i+1}],i=0,\dots,k$ is contained in a single $J_{\alpha_i}$. Choose such an $\alpha_i$ for each $i=0,\dots,k$. Now choose $t_i \geq 0$ for $i=0,\dots,k$ large enough to use the flow of $-X^{s_i}$ as in Lemma~\ref{lem:indwell} to get homotopy equivalences
  \begin{align*}
    A_{\alpha_{i-1}}/B_{\alpha_{i-1}} \to A_{\alpha_i}/B_{\alpha_i}.
  \end{align*}
  Now the composition of these defines a homotopy equivalence as wanted. This is a contractible choice since we may always introduce new division points and increase flow times. Doing this we can cut up any interval $J_{\alpha_i}$ and replace the index pair we have on that with any refinement of index pairs on a cover of $J_{\alpha_i}$.

  This by construction (the reason for having $s_k=s_{k+1}$) has the flows from Lemma~\ref{lem:indwell} build in at the very beginning and end of the interval, hence naturally commutes with these, by appropriately changing the flow times at each end.
\end{proof}

\begin{Remark}
  \label{rem:Homoindex:1} 
  More generally it is proven in \cite{MySympfib} that for a different base manifold $B$ of a projection $M\to B$ the Conley indices defined as a parameterized based space over $B$ behaves very much like a (based) Serre fibration, and hence the homotopy equivalences from one fiber to another can be thought of as a homotopy lifting property (parallel transport of the fiber). This means that in the above argument it is important that $I$ is contractible - indeed, in the general case the homotopy equivalence would depend on choices of paths in the base $B$.
\end{Remark}

\subsection{Completely Bounded Pseudo-Gradients}

The homotopy type of a Conley index is particular nice to work with in the following case.

\begin{Definition}
  \label{def:Homoindex:1}
  A Pseudo-gradient $X$ for a function $f\co M \to \R$ is said to be completely bounded (CB) if;
  \begin{itemize}
  \item the flow of $-X$ is defined for all times (positive and negative) and
  \item there exists a compact set $K\subset M$ and $k>0$ such that $X(f)>k$ on the complement of $K$.
  \end{itemize}
\end{Definition}

Should is analogous to the Palais-Smale condition.

\begin{Lemma}
  \label{lem:Homoindex:1}
  If $(f,X)$ is CB then there exists a good index pair and
  \begin{align*}
    I_a^b(f,X) \simeq f^{-1}([a,b])/f^{-1}(a)
  \end{align*}
  this homotopy equivalence is a contractible choice (in both directions) compatible with all the above homotopy equivalences.
\end{Lemma}

\begin{proof}
  Using the fact that the flow $\psi_t$ is defined everywhere we get canonical deformation retractions of $f^{-1}([a,b])/f^{-1}(a)$ onto the quotients of the sets
  \begin{align*}
    A_t & = f^{-1}([a,\infty[) \cap (f\circ \psi_{-t})(]-\infty,b]) \\
    B_t & = A_t \cap f^{-1}(a).
  \end{align*}
  Indeed we produce the map by taking the flow on $f^{-1}([a,b])$ and then collapsing everything with $f\leq a$.

  Claim $A_t$ is compact for $t>>0$. To see this let $K\subset M$ and $k>0$ be as in Definition~\ref{def:Homoindex:1}. Define
  \begin{align*}
    K' = \psi_{[0,T]}(K),
  \end{align*}
  with $T>(b-a)k^{-1}$. Then we claim that $A_T\subset K'$ and hence $A_t$ is compact for $t\geq T$. Indeed, let $x\in A_T - K'$ be given. Then by definition of $K'$ the points $\psi_t(x)$ for $t\in[-T,0]$ are not in $K$ hence $X(f)>k$. It follows that the value of $f$ when flowing on $\psi_{-T}(x)$ for time $t\in [0,T]$ decreases faster than $k$ (that is $\pdd{t}$ of this is less than $-k$), but since we are flowing for more time than $(b-a)k^{-1}$ it must decrease totally more than $(b-a)$ - hence $f(\psi_{-t}(x))>b$ which is a contradiction.

  Claim: $(A_t,B_t)$ is a good index pair for $(f,X)$ when $t\geq T$. Indeed, $A_t$ is cut out by two equations which are transversal to the flow, so the fact that the flow of $-X$ points out at the set $f^{-1}(a)$ is preserved under small perturbations of $X$, and the fact that it points in when $f \circ \psi_t(x) = b$ (but not $f(x)=a$) is similarly preserved.

  These homotopy equivalences are compatible (up to contractible choices) with all the above since they are also given by the flow. The homotopy equivalence in the other direction is induced by the inclusion $(A_t,B_t) \subset (f^{-1}([a,b]),f^{-1}(b))$. 
\end{proof}

This lemma also suggest that the index when $X$ is CB does not in fact depend on $X$. We formalize this in the next lemma.

\begin{Lemma}\label{lem:234}
  If two pseudo-gradients $X$ and $X'$ for $f$ are CB, then $I_a^b(f,X)\simeq I_a^b(f,X')$. Again, this choice is contractible compatible with all of the above.
\end{Lemma}

\begin{proof}
  We simply notice that if $(A,B)$ is an index pair for $(f,X)$ and $(A',B')$ for $(f,X')$ then the two inclusions
  \begin{align*}
    (A,B) \subset (f^{-1}[a,b],f^{-1}(b)) \supset  (A',B')
  \end{align*}
  induce the homotopy equivalences, and to get a contractible choice of these we simply compose with the flow as above.
\end{proof}

\begin{Remark}
  \label{rem:Homoindex:2}
  One can view this result in a different way. Indeed, the set of CB pseudo-gradients is in fact contractible so there is a homotopy between $X$ and $X'$ within CB pseudo-gradients. Then Lemma~\ref{hominv} proves the independence.
\end{Remark}

\subsection{Alternative Construction of Good Index Pairs}

Inspired by the above construction we now describe a slightly more general way of producing good index pairs (which is handy when $X$ is not CB) by using what we will call \emph{cut-off} functions.

\begin{Lemma} \label{cutoff}
  Assume that $g_1,g_2,\dots,g_n\colon M\to \R$ are continuous functions such that
  \begin{align*}
    A&=f^{-1}([a,b])\cap \{x\in M\mid g_j(x) \leq 0\} \\
    B& = f^{-1}(a) \cap A
  \end{align*}
  are compact and the interior of $A$ contains all the critical points of $f$ with critical value in $]a,b[$. If; for each $x\in \partial A$ with $g_j(x)=0$ we have that $g_j$ is smooth in a neighborhood of $x$ and
  \begin{align} \label{bound}
     - X_x(g_j) < 0;
  \end{align}
  then $(A,B)$ is a good index pair.
\end{Lemma}

The meaning of \eqref{bound} is the following: Because of the definition of $A$ we see that $-X$ must flow into $A$ at all parts of the boundary except of course when $f$ crosses the value $a$.

\begin{proof}
  At any point $x\in \partial A$ we must have $f(x)-b\leq 0$ and 
  \begin{align*}
    g_j(x) \leq 0
  \end{align*}
  satisfied. Since $b$ is a regular value we have that $-X_x(f-b)<0$ if $f(x)-b=0$ and, by assumption, $-X_x(g_j)<0$ if $g_j(x)=0$. So we see that for any vector $v$ close to $X_x$ we have that: if $f(x)=b$ then $-v(f)<0$ and if $g_j(x)=0$ then $-v(g_j)<0$. So $-v$ points into $A$, except if the equality $f(x)=a$ holds, in which case $-v$ must point out of the set for the same reason. The boundary is compact so there is an $\epsilon>0$ such that this holds for all $v$ and all $x$ if $\norm{v-X_x}<\epsilon$
\end{proof}


\section{The Action Integral in Cotangent Bundles}\label{action}

We will once and for all fix a Riemannian metric on $N$ (and at the end of the paper argue that the entire construction does not depend on this choice). This section recalls some notions and introduces some notation concerning the action integral
\begin{align*}
  A_H(\gamma) = \int_\gamma (\lambda - Hdt),
\end{align*}
where $\gamma$ is a closed curve (of sufficient regularity) in $T^*N$. All parts of this section are well-known, but the methods are vital to the construction.

We denote points in the cotangent bundle $T^*N$ by $(q,p)$, where $q$ is in $N$ and $p$ is a cotangent vector at $q$. Let $\pi\colon T^*N \to N$ be the projection onto the base and define the canonical 1-form $\lambda\in\Omega^1(T^*N)$ and 2-form $\omega\in\Omega^2(T^*N)$ by 
\begin{align*}
  \lambda_{q,p}(v) & = p(\pi_*(v))  \\
  \omega & =\dd (-\lambda).
\end{align*}
The form $\omega$ is non-degenerate and thus defines a canonical symplectic structure on $T^*N$. 

Given any smooth Hamiltonian $H\colon T^*N \to \R$, we may define the associated Hamiltonian vector field $X_H$ by the formula $dH=\omega(X_H,-)$. This is well-defined because $\omega$ is non-degenerate. The flow of $X_H$ will be denoted $\varphi^H_t$ and is called the Hamiltonian flow.

Using the Riemannian structure on $N$ we may induce a Riemannian structure on $T^*N$ in the following way: at each point $(q,p)$ we split the tangent space $T_{(q,p)}(T^*N)$ in two components, the vertical, which is canonically defined without the metric as the kernel of $\pi_*$, and the horizontal defined by the connection given by the metric on $N$. This identifies $T_{(q,p)}T^*N$ with $T_qN\oplus T_q^*N$, on which we use the structure from $N$ to define the inner product on each factor - making this splitting orthogonal. We may also define an almost complex structure $J$ in this splitting by using the isometry $\phi_q\colon T_qN\to T_q^*N$ induced by the metric on $N$
\begin{align*}
  J(\delta q,\delta p)=(-\phi^{-1}(\delta p),\phi(\delta q)).
\end{align*}
This is compatible with the symplectic structure and the induced Riemannian structure. The formula for $X_H$ can be rewritten using these as 
\begin{align} \label{xh}
  X_H=-J\nabla H.
\end{align}
For any smooth manifold $M$ let $\Lambda M$ be the space of piecewise smooth and continuous maps from $S^1=I/\{0,1\}$ to $M$. The action $A_H\colon \Lambda T^*N \to \R$ is defined by
\begin{align*}
  A_H(\gamma)=\int_\gamma \lambda - \int_{S^1}H(\gamma(t))dt =
  \int_\gamma (\lambda - Hdt).
\end{align*}
It is known that the critical points of this integral are precisely the $1$-periodic orbits of the Hamiltonian flow (the calculation in Equation~\eqref{piecedif} in  section \ref{flowline} proves this).

We will often need the special case in which $H$ only depends on the length of the cotangent vector - that is
\begin{align*}
  H(q,p)=h(\norm{p}),
\end{align*}
where $h\colon \R \to \R$. For $H$ to be smooth the germ at $0$ of $h$ needs to be even.

In this case we calculate the gradient of $H$ in the orthogonal splitting:
\begin{align*}
  \nabla H = (0,h'(\norm{p})\frac{p}{\norm{p}})
\end{align*}
We get $0$ in the first factor because parallel transport does not change the norm of $p$. Using equation \eqref{xh} we see that
\begin{align*}
  X_H=-J(0,h'(\norm{p})\frac{p}{\norm{p}})=(h'(\norm{p})
  \phi^{-1}(\frac{p}{\norm{p}})
  ,0)=(h'(\norm{p})\frac{p}{\norm{p}},0).
\end{align*}
As the last equation indicates we will from now on suppress $\phi$ from the notation.

\begin{Remark} \label{geocalc}
  Because this vector field is $0$ on the vertical factor, it will parallel transport $p$ and hence this becomes a reparameterization of the
  \begin{figure}[ht]
    \centering
    \includegraphics{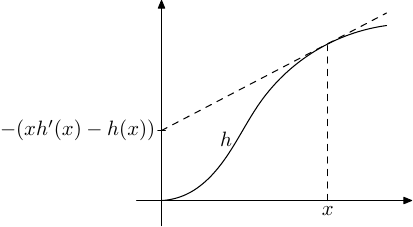}
    \caption{Geometric calculation of critical values.}
    \label{yaction}
  \end{figure}
  geodesic flow on $N$. This describes the $1$-periodic orbits as closed geodesics on $N$ with lengths corresponding to $h'(\norm{p})$. The action of these orbits is easily calculated to be $\norm{p}h'(\norm{p})-h(\norm{p})$. This corresponds to taking minus the intersection of the $y$-axis with the tangent of $h$ at the point $(x,h(x))$ as in figure \ref{yaction}. This geometric formula for calculating the action is very useful for this type of Hamiltonian, and will be used repeatedly.
\end{Remark}


\section{The Gradient of a Flow-line segment in Cotangent
  Bundles} \label{flowline}

We will in this section define what we call a segment function which we will use as building blocks in the finite dimensional approximations in section \ref{Florlike}. We will only define these for Hamiltonians $H \colon T^*N \to \R$ with small $C^2$-norm. The construction may seem technical, but it has the advantage of being very explicit. This section and Section~\ref{Florlike} are inspired by work in \cite{MR1617648} and \cite{MR765426}.

We assume that the injective radius of $N$ is $2\delta_1$. Define 
\begin{align*}
  D_RT^*N = \{(q,p)\in T^*N \mid \norm{p}\leq R\}
\end{align*}
and $DT^*N=D_1T^*N$. Similarly, define
\begin{align*}
  U_RT^*N = \{(q,p)\in T^*N \mid \norm{p} = R\}
\end{align*}
and $UT^*N=U_1T^*N$. We will define the segment functions on the space
\begin{align*}
  W = \{(q',p',q)\in T^*N \times N\mid \dist(q',q) \leq \delta_1\}
\end{align*}
where $\dist(-,-)$ is the distance in $N$ using the Riemannian structure. We also defined the compact sub-space
\begin{align*}
  DW = \{(q',p',q)\in DT^*N \times N\mid \dist(q',q) \leq \delta_1\}
\end{align*}

We will in the entire paper only consider smooth Hamiltonians $H\colon T^*N\to \R$ with the property: there exist $\mu>0$ and $c\in\R$ such that
\begin{align*}
  H(q,p)=\mu\norm{p}+c \qquad \textrm{ for } (q,p) \notin DT^*N.
\end{align*}
For such Hamiltonians we define
\begin{align} \label{eq:12}
  \no{H} = \smashoperator{\sup_{z\in DT^*N}}
  (\norm{\nabla H},\norm{\nabla \nabla H})
\end{align}
where $\norm{-}$, $\nabla$, and $\nabla \nabla$ are defined using the Riemannian structure on $T^*N$ (induced by the Riemannian structure on $N$ - defined in Section~\ref{action}). Notice that we did not include any value of $H$, so a constant function has norm 0 - making this a seminorm. Indeed, all we care about are bounds on the gradient and the second order behavior of $H$. 

\begin{figure}[ht]
  \centering
  \includegraphics{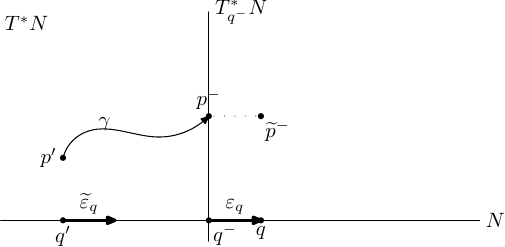}
  \caption{Flow-line segment and related quantities}
  \label{segmentpic}
\end{figure}

In this section we impose the condition $\no{H} < \delta_1/10$ on $H$. This implies in particular that $\mu < \delta_1/10$ and thus the slope at infinity is less than the length of any non-constant geodesic starting and ending at the same point. For such $H$ we define the \textbf{segment function} $S^H \co W\to \R$ by
\begin{align}\label{eq:8}
  S^H(q',p',q) = \pare*{\int_{\gamma}\lambda-Hdt} +
  p^-\epsilon_q,
\end{align}
where $\gamma\colon[0,1]\to T^*N$ is the Hamiltonian flow
curve $\gamma(t)=\varphi^H_t(q',p')$,
\begin{align*}
  (q^-,p^-)=\gamma(1) \qquad\text{and}\qquad
  \epsilon_q=\exp^{-1}_{q^-}(q) \in T_{q^-} N,
\end{align*}
with $\exp\colon TN \to N$ the exponential map. These are illustrated
in Figure~\ref{segmentpic}.

The term $p^-\epsilon_q$ is the pairing of cotangent vectors with
tangent vectors and thus the symplectic area of the rectangle formed
by $p^-,\tp^-,q,q^-$ ($\tp^-$ is defined properly below). The function
$S^H$ is well-defined because $\norm{\nabla H} <\delta_1/10$ implies
that the distance between $q^-$ and $q$ is less than
$\dist(q,q')+\dist(q',q^-) < 11\delta_1/10$ which is less than the
injective radius $2\delta_1$. Notice that if $\norm{p}>1$, then $p^-$
is the parallel transport of $p'$ by a geodesic in the direction of
$p'$ (see Section~\ref{action}), and thus $\norm{p^-}=\norm{p'}$. 

Consider the commutative diagram of isometries 
\begin{align*}\xymatrix{
  T_{q_1}N \ar[r]^{P_{q_1,q_2}} \ar[d]^{\phi} & T_{q_2}N \ar[d]^{\phi} \\
  T_{q_1}^*N \ar[r]^{P_{q_1,q_2}^*} & T^*_{q_2}N,
}\end{align*}
where $P^{(*)}_{q_1,q_2}$ is given by parallel transport along the
unique geodesic when $q_1,q_2\in N$ satisfy $\dist(q_1,q_2)<
2\delta_1$. The isomorphism $\phi$ is the one induced by the metric,
which we suppressed from the notation in the previous section. We will
do so again and thus $P^*_{q_1,q_2}=P_{q_1,q_2}$. We use this to define  
\begin{align*}
  \tp^- = P_{q^-,q}(p^-) \in T^*_qN,
\end{align*}
and a parallel transported version of $\epsilon_{q_j}$ by
\begin{align*}
  \te_q=P_{q^-,q'}(\epsilon_q) \in T_{q'}^*N.
\end{align*}
These are also illustrated in figure \ref{segmentpic}. To control our
finite dimensional approximations defined later we will need the
following facts about the gradient of these segment functions.

We will be using the notation
\begin{align*}
  \nabla S^H = \nabla_{(q',p')}S^H \oplus \nabla_q S^H = \nabla_{q'} S^H \oplus \nabla_{p'} S^H \oplus \nabla_{q} S^H.
\end{align*}
Here the splitting of $\nabla_{(q',p')} S^H$ into two factors is
horizontal and vertical directions as described in
section~\ref{action}.

\begin{Lemma} \label{gradient1half}
  There exists constants $C,\delta>0$ such that for any Hamiltonian
  $H$ with $\no{H} < \delta$ we have
  \begin{align}
    \norm{\nabla_{q'} S^H + p'} \leq & C \norm{\epsilon_q} \label{onehalf} \\
    \norm{\nabla_{p'} S^H - \te_q } \leq &
    \tfrac{1}{4}\norm{\epsilon_q}  \qquad \qquad (=
    \tfrac{1}{4}\norm{\te_q}) \label{twohalf} \\
    \norm{\nabla_{q} S^H -\tp^-} \leq & C \norm{\epsilon_q}, \label{threehalf}
  \end{align}
  on the compact set $DW$.
\end{Lemma}

Notice in particular the very important fact that Equation~\eqref{twohalf} implies that a critical point has $\epsilon_{q_j}=0$, which then have serious implications for the two others at critical points. Indeed, critical points are thus small flow lines starting and ending on the zero-section. Pasting these cyclically together in the next section we get periodic orbits as critical points.

\begin{proof}
  First we consider the integration term of $S^H$. We start by seeing how it depends tangentially on the curve $\gamma=(\gamma_q,\gamma_p)$ -  thinking of $\gamma$ as an independent variable. So in the following $\partial \gamma=(\partial \gamma_q , \partial \gamma_p)$ is a smooth tangent field along $\gamma$.
  \begin{align}
    D_\gamma (\int_{-}\lambda&-Hdt)  (\partial \gamma)  \nonumber \\
    = & \int_0^1
    \gamma_{p}(\nabla_t\partial\gamma_{q}) + (\partial\gamma_{p})\gamma_{q}' -
    (\nabla_{\gamma}H)(\partial\gamma)\dd t \nonumber \\
    = & [\gamma_{p}(t)\partial\gamma_{q}(t)]_0^1 + \int_0^1
    - \gamma_{p}'(\partial\gamma_{q}) + (\partial\gamma_{p})\gamma_{q}' -
    (\nabla_{\gamma}H)(\partial\gamma)\dd t  \nonumber \\
    = & - p' \partial q' + p^-\partial q^- -
    \int_0^1 (J\gamma'+\nabla_{\gamma}H)(\partial
    \gamma)\dd 
    t \label{piecedif} \\
    = & - p' \partial q' + p^-\partial q^- \nonumber
  \end{align}
  The integral vanishes because $\gamma$ is a Hamiltonian flow curve,
  and thus $\gamma'(t)= - J \nabla_{\gamma(t)} H$. This is a standard
  calculation, and it is also a proof that the 1-periodic orbits of
  the flow $X_H$ are the critical points of the action integral.

  Motivated by this fact that the linearization of the action only
  depends on the start and end point of $\gamma$ it is classical to
  extend $S^H$ to a larger manifold where these are considered
  independent variables. Formally we do this by picking for all $(q',p')\in DT^*N$ and
  $(q^-,p^-)\in DT^*N$ with $\dist(q',q^-)< \delta_1/2$ a smooth path
  starting at $(q',p')$ and ending at $(q^-,p^-)$ such that when
  $(q^-,p^-)=\varphi_1^H(q',p')$ this curve is the Hamiltonian flow
  curve of $H$ (it is convenient that the flow of $X_H$ preserves $DT^*N$ since otherwise this would clutter the notation a bit). We may assume that this choice is smooth in all
  variables meaning that the adjoint map is smooth. By abuse of
  notation we denote this choice of curves $\gamma$. Now define an
  extension $G$ of $S^H$ by using the ``same'' formula
  \begin{align*}
    G(q',p',q^-,p^-,q) = \pare*{\int_\gamma \lambda -H dt} +
    p^-\epsilon_q, 
  \end{align*}
  where $\epsilon_q$ is extended simply by $\epsilon_{q}=\exp_{q^-}^{-1}(q')$. So instead of imposing $(q^-,p^-)=\varphi^H_1(q',p')$, as we did in the definitions of $S^H$, we use the chosen smooth family of $\gamma$'s to define $G$ depending on \emph{independent} variables $(q',p',q^-,p^-,q)$. With this definition we see that $S^H(q',p',q)=G(q',p',\varphi_1^H(q',p'),q)$, meaning that if we define an embedding of manifolds
  \begin{align*}
    i\co DW \to DT^*N\times DT^*N\times N
  \end{align*}
  by the formula
  \begin{align*}
    i(q',p',q) = (q',p',\varphi_1^H(q',p'),q)
  \end{align*}
  then $G$ is defined in a neighborhood of the image and $S^H=G \circ i$, which justifies calling this an extension.

  We may use the gradient of $G$ on the image if $i$ to
  calculate the gradient of $S^H$ using the chain rule $DS^H = DG
  \circ Di$ and the fact that the gradient $\nabla f$ is the image
  of 1 under the adjoint $(Df)^\dagger$ of $Df$ for any $f\colon W \to
  \R$. More concisely, we have 
  \begin{align*}
    \nabla S^H = (Di)^\dagger \nabla G,
  \end{align*}
  which more concretely in our case turns into
  \begin{align} \label{gradcal1}
    \nabla_{(q',p')} S^H = \nabla_{(q',p')} G +
    (D_{(q',p')}\varphi^H_1)^\dagger(\nabla_{(q^-,p^-)} G)
  \end{align}
  and
  \begin{align} \label{gradcal2}
    \nabla_q S^H = \nabla_q G.
  \end{align}
  
  So we will calculate $\nabla G$ on the image of $i$. Comparing the calculation in Equation \eqref{piecedif} with the way we defined the Riemannian structure we already calculated this gradient of the extended integration term:
  \begin{align*}
    \nabla_{(q',p',q^-,p^-,q)}\int_\gamma \lambda - Hdt = (-p',0,p^-,0,0).
  \end{align*}
  The last zero is because nothing in this term depends on $q$. Notice that this \emph{only} works on the image of $i$, and the fact that it does not depend on the choice of extension $\gamma$ is because Equation~\eqref{piecedif} shows that it only depends on the end points on the image of $i$.
  
  The gradient of the second term $g(q^-,p^-,q):=p^-\epsilon_q$ is a 
  little more tricky, but at least it does not depend on
  $(q',p')$ nor $H$. The gradient with respect to $p^-$ is easy:
  \begin{align*}
    \nabla_{p^-} g = \epsilon_q.    
  \end{align*}
  For the remaining two factors we first assume that $\epsilon_q=0$,
  i.e. $q^-=q$ then by looking in a normal coordinate chart we see
  that
  \begin{align*}
    \nabla_{(q^-,q)} g = (-p^-,p^-) = (-p^-,\tp^-) \qquad
    \textrm{when} \qquad \epsilon_q = 0.
  \end{align*}
  We rewrite the latter $p^-$ as $\tp^-$ because when looking at
  $\epsilon_q\neq 0$ this vector lives in the correct tangent space
  $T_qN$. We now claim that this implies the bounds
  \begin{align*}
    \norm{\nabla_{q^-} g +p^-} \leq C'\norm{\epsilon_q} \qquad \textrm{and} \qquad \norm{\nabla_{q'} g - \tp^-} \leq C'\norm{\epsilon_q}
  \end{align*}
  for some $C'>0$ for $(q^-,p^-,q) \in DT^*N\times N$. Indeed, we are in the following abstract
  situation: we have a smooth section $\nabla_{q^-}g + p^-$ (similar for $\nabla_{q'} g - \tp^-$) in a vector bundle with a smooth metric on a compact manifold, and this section is 0 when another smooth section $\epsilon_q$ is 0 (in the second case this other smooth section is $\exp^{-1}_{q}(q^-)=P_{q^-,q}(-\epsilon_q)$ which has the same norm), and this other smooth section is transverse to the zero section. In this case we can always locally find $C'$ such that the bounds are true - and we are on a compact set so there is a global $C'$ as well. In particular notice that this $C'$ does not depend on $H$. Indeed, $g$ does not depend on $H$.

  By adding the gradients of the two terms we obtain the gradient of $G$ (on the image of $i$) as
  \begin{align*}
    \nabla_{(q',p',q^-,p^-,q)} G =
    (-p',0,b_1,\epsilon_q,\tp^-+b_2),
  \end{align*}
  where $b_1$ and $b_2$ are smooth sections (vector fields) whose norms are bounded by
  the function $C'\norm{\epsilon_q}$.

  Now we use Lemma \ref{Hamflow} below with
  \begin{align*}
    F_{(q',p'),(q^-,p^-)}=P_{q',q^-}\oplus P_{q',q^-}  
  \end{align*}
  defined using the parallel transport on the usual splitting of the
  tangent spaces
  \begin{align*}
    P_{q',q^-} \oplus P_{q',q^-} \colon
    T_{q'}N \oplus T_{q'}N \to T_{q^-}N \oplus T_{q^-}N \approx
    T_{(q^-,p^-)}(T^*N),
  \end{align*}
  and with $M=DT^*N$. From that lemma we get that if $\no{H} <
  \delta=\delta(\epsilon)$ then
  \begin{align*}
    \norm{D \varphi_1^H-F} \leq \epsilon,
  \end{align*}
  which implies
  \begin{align} \label{opnorm}
    \norm{(D \varphi_1^H)^\dagger - F^{-1}} \leq \epsilon
  \end{align}
  because $F^\dagger = F^{-1}$ since $F$ is an isometry. This means
  that $(D\varphi_1^H)^\dagger$ is $\epsilon$ close in operator norm to
  the isometry $F^{-1}$ which sends $(b_1,\epsilon_q)$ to
  $(P_{q^-,q'}(b_1),\te_q)$. If we had equality $F^{-1}=D\varphi^H_1$
  then we would get from equations \eqref{gradcal1} and
  \eqref{gradcal2} that the gradient of $S^H$ were given by
  \begin{align*}
    \nabla_{q',p',q} S^H = (-p'+P_{q^-,q'}(b_1),\te_q,\tp^-+b_2),
  \end{align*}
  which would easily imply the bounds
  \begin{align*}
    \norm{\nabla_{q'} S^H + p'}  &\leq C' \norm{\epsilon_q} \\
    \norm{\nabla_p S^H - \te_q} &\leq 0 \\
    \norm{\nabla_q S^H-\tp^-}&\leq C' \norm{\epsilon_q}.
  \end{align*}
  However the difference of using $F^{-1}$ and $(D\varphi_1^H)^\dagger$
  in \eqref{gradcal1} is bounded by the operator norm in
  \eqref{opnorm} times the norm of the vector on which we use them:
  \begin{align*}
    \norm{(D\varphi_1^H)^\dagger-F^{-1}}\norm{b_1,\epsilon_q} \leq
    \epsilon (C'+1)\norm{\epsilon_q} \leq
    \tfrac{1}{4}\norm{\epsilon_q}.
  \end{align*}
  The latter if $\epsilon$ was such that $4\epsilon (C'+1)\leq 1$,
  which is true for appropriate $\delta=\delta(\epsilon)$. So we get
  the wanted bounds if we pick $C=(C'+1/4)$ and such a $\delta$.
\end{proof}

The following lemma gives approximations of the same gradients but on the complement of the compact set $DW$.

\begin{Lemma} \label{gradient1}
  For the same $C,\delta>0$ as in Lemma~\ref{gradient1half} we have for any Hamiltonian $H$ with $\no{H} < \delta$ that
  \begin{align}
    \norm{\nabla_{q'} S^H + p'} \leq & C
    \norm{p'} \norm{\epsilon_q} \label{one} \\
    \norm{\nabla_{p'} S^H - \te_q } \leq &
    \tfrac{1}{4}\norm{\epsilon_q}  \qquad \qquad (=
    \tfrac{1}{4}\norm{\te_q}) \label{two} \\
    \norm{\nabla_{q} S^H -\tp^-} \leq & C
    \norm{p'} \norm{\epsilon_q}, \label{three}
  \end{align}
  on $W-DW$.
\end{Lemma}

\begin{proof}
  We will prove this on the entire set $U=\{\norm{p'}\geq 1\}\cap W \supset (W-DW)$.
  Lemma~\ref{gradient1} proves this on the subset $\{\norm{p'}=1\}\cap U$.

  On the set $U$ we have $H(q',p')=\mu\norm{p'}+c$. The description of the flow curves and their actions in section \ref{action} implies that the integration part of $S^H$ is constantly equal to $-c$ on $U$ and that $p^-$ is the parallel transport of $p'$ along the geodesic in direction $p'$ with length $\mu$. This geodesic is also the projection of the flow curve $\gamma$ to $N$. Let $t\geq 1$ be given. We wish to analyze how the term
  $p^-\epsilon_q$ behaves if we multiply the $p'$ coordinate with this
  $t$. Since the projected geodesic is the same for $p'$ and $tp'$
  (both have length $\mu$ and points in the direction given by $p'$)
  and since the parallel transport is linear we see that the term
  $p^-\epsilon_q$ simply gets multiplied with $t$ because $p^-$ does
  so. We have argued that 
  \begin{align*}
    t(S^H(q',p',q)+c) = S^H(q',tp',q)+c \qquad \textrm{for} \qquad
    \norm{p'}\geq 1, t\geq 1.
  \end{align*}
  It is now easy to verify that the gradient of $S^H$ with respect to $q'$ and $q$ scales with $t$ and that the gradient with respect to $p'$ is independent of $t$.
\end{proof}

In the above we used the following lemma, and we will need it
in the following generality later. So let $M$ be any compact almost
k\"ahler manifold (possibly with boundary, corners, etc.). For any
Hamiltonian $H\colon M\to \R$ we define 
\begin{align} \label{eq:19}
  \no{H} = \smashoperator{\sup_{z\in M}}
  (\norm{\nabla H},\norm{\nabla \nabla H})
\end{align}
We also assume that the Hamiltonian flow preserves $M$ (although this is not really necessary if we put $M\subset M'$ where $M'$ is open).

Let $F_{(z_1,z_2)} \colon T_{z_1} M \to T_{z_2}M$ be any smooth
identification of close-by tangent spaces, i.e. $F_{(z_1,z_2)}$ is a
linear isomorphism defined for $\dist(z_1,z_2) \leq \epsilon_1$ and
smooth in $z_1\in M$ and $z_2\in M$. Furthermore, we assume that
$F_{(z,z)}$ is the identity on $T_zM$. 

\begin{Lemma} \label{Hamflow}
  For any $\epsilon>0$ we may find $\delta>0$ such that if 
  $\no{H}<\delta$ then the Hamiltonian flow $\varphi_t^H$
  satisfies
  \begin{align}\label{resul}
    \norm{(D_z\varphi_1^H)-F_{z,\varphi_1^H(z)}} \leq \epsilon
  \end{align}
  for all $z\in M$. Here the norm is the operator norm.
\end{Lemma}

\begin{Remark}
  This is equivalent to the well-known lemma that if $H$ is $C^2$-close to a constant map then the time 1-flow is $C^1$-close to the identity.
\end{Remark}

\begin{proof}
  Since $\norm{\nabla H} < \delta$ implies $\dist(z,\varphi_1^H(z)) < \delta$ we
  see that the left hand side of \eqref{resul} is well-defined for
  small $\delta$. By compactness we may find a finite set of
  symplectic charts $h_i\colon U_i \to M$ with $U_i \subset \R^{2n}$
  such that; for small $\delta>0$ each flow
  curve $\varphi_t(z)$, $t\in [0,1]$ is fully contained in one of
  these charts for all $z\in K$ and all first and second order
  derivatives of all the $h_i$'s are bounded.

  Define $H_i=H \circ h_i$ for any $H$, then the bounds on 
  $h_i$ implies that we can assume that there is a constant $K>0$ such
  that $\no{H_i} \leq K\no{H}$. So by making $\delta$ small we can
  make all these norms small.

  Since we may also assume that the charts have diameter less than
  $\epsilon_1$ we get that $F$ pulled back to any of the charts (in
  the obvious sense), call this $F^i$, defines linear functions
  $F^i_{z_1,z_2}\colon \R^{2n} \to \R^{2n}$ smoothly dependent on
  $z_1,z_2 \in U_i$ such that $F_{z,z}=\id$. This implies that if
  $\dist(z,\varphi_1^H(z)) < \norm{\nabla H_i} \leq \no{H_i}$ is small
  enough we get 
  \begin{align*}
    \norm{F^i_{z,\varphi_1^H(z)}-\id} \leq \epsilon/2,
  \end{align*}
  for all $i$ and all $z\in U_i$ simultaneously by compactness (maybe
  we shrink all $U_i$'s a little).
  
  We may now work entirely in one of these charts, and by abuse of
  notation use $H=H_i$, $F=F^i$, $U=U_i$, $\no{H}=\no{H_i}$.

  Since we are now in the case of the standard flat metric in
  $\R^{2n}$ we see that 
  \begin{align*}
    \norm{\nabla X_H}  = \norm{\nabla^2 H} < \no{H}
  \end{align*}
  implies that $\norm{D_z\varphi_1^H-\id} < \no{H}$ by a standard
  integration argument. So for $\no{H} < \epsilon/2$ we have
  \begin{align*}
    \norm{D_z\varphi_1^H - F_{z,\varphi_1^H(z)}} \leq 
    \norm{D_z\varphi_1^H - \id} +
    \norm{F_{z,\varphi_1^H(z)} - \id} \leq \epsilon,
  \end{align*}
  which is what we set out to prove.
\end{proof}


\section{Finite Dimensional Approximation of the Action Integral in Cotangent
  Bundles} \label{Florlike}

In this section we define finite dimensional approximations $S_r$ to the action $A_H$ by putting several segment functions together. This means we no longer need the Hamiltonian $H$ to be $C^2$ small, but the number $r$ of segment functions needed then depends on the $C^2$-norm. We will then define a pseudo-gradient $X_r$ for this finite dimensional approximation such that there exists good index pairs for large $r>>0$ and hence we have well-defined Conley indices.

As in section \ref{flowline} we assume that all Hamiltonians $H\colon T^*N\to \R$ are smooth and linear outside $DT^*N$ with some slope $\mu$. Again we define $\no{H}$ as in Equation~\eqref{eq:12}. We additionally assume that the slope $\mu$ is not the length of a closed geodesic on $N$.

Let $\delta_1$ be as in Section~\ref{flowline}. Now define
\begin{align} \label{eq:11}
  \delta_0 = \min(\delta_1,\delta,(8C)^{-1}),
\end{align}
where $\delta$ and $C$ are the constants from Lemma~\ref{gradient1half} and Lemma~\ref{gradient1}. We also assume that any ball in $N$ with radius less than $\delta_0$ is geodesically convex. Use this to define the manifold of $r$-piecewise geodesics in $N$ by
\begin{align} \label{eq:4}
  \Lre N = \{(q_j)_{j\in \Z/r}\in N^r\mid
  \dist(q_j,q_{j+1})<\delta_0 \}.
\end{align}
We use $j\in \Z/{r}$ to emphasize that $q_{r}=q_{0}$ and we have a cyclic structure. We see that the cotangent space of this is easily identified by
\begin{align}\label{eq:5}
  T^* \Lre N \cong \{(q_j,p_j)_{j\in \Z/r} \in (T^*N)^r \mid
  \dist(q_j,q_{j+1})<\delta_0\}.
\end{align}
We will denote a point in this space by $\arz=(\arq,\arp)$, and a single coordinate by $z_j=(q_j,p_j) \in T^*N$. These two spaces are given the restriction of the product Riemannian structures from $N^r$ and $T^*N^r$.

We will use the segment functions defined in section \ref{flowline} to define functions resembling the action $A_H$ on $T^*\Lre N$ having the same critical points (the 1-periodic orbits) with the same critical values. In section \ref{loclinsec} we will define generalizations of these, and in Remark \ref{embedrem} we explain how they can be constructed using embeddings $i_r \colon T^*\Lambda_r N \to \Lambda T^*N$ as explained in the overview in the introduction.

\begin{figure}[ht]
  \centering
  \includegraphics{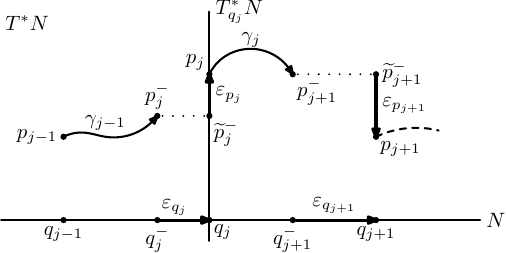}
  \caption{The piece-wise flow and relevant tangent vectors.}
  \label{pflowq}
\end{figure}

Since parameterization of flow curves will be important to handle in the construction later we define such ``approximations'' for certain subdivision of the unit interval $I=[0,1]$. So, let $\alpha=(\alpha_j)_{j\in \Z/r} \in I^r$ with $\sum_j \alpha_j =1$ be given. However, since we will never need to consider subdivisions which does not satisfy
\begin{align} \label{eq:52}
  \alpha_j \leq \tfrac{2}{r}
\end{align}
we will for the sake of simplicity always assume this. Note, that the most important example of these is of course $\alpha_j=1/r$. We will also assume that $r$ is large enough for
\begin{align} \label{eq:10}
  \tfrac{2}{r}\no{H} < \delta_0/10.
\end{align}
Under these two assumptions (and only under these two assumptions) we define
\begin{align} \label{eq:7}
  S_r(\arz) = S_{r,\alpha}^H(\arz) &= \sum_{j\in \Z/r} S^{(\alpha_j H)}(q_j,p_j,q_{j+1}),
\end{align}
where $S^{(\alpha_j H)}$ is the segment function defined in Equation~\eqref{eq:8} for the Hamiltonian $\alpha_jH$. Also define $\gamma_j\colon[0,\alpha_j]\to T^*N$ as the Hamiltonian flow curve given by $\gamma_j(t)=\varphi^{H}_t(q_j,p_j)$, and further define
\begin{align*}
  (q_j^-,p_j^-)&=\gamma_{j-1}(\alpha_{j-1}), \\
  \epsilon_{q_j}&=\exp^{-1}_{q_j^-}(q_j) \in T_{q_j^-} N, \\
  \te_{q_j}&=P_{q_j^-,q_{j-1}}(\epsilon_{q_j}) \in T_{q_{j-1}}N, \\
  \tp_j^- &= P_{q_j^-,q_j}(p^-) \in T_{q_j}^*N,  \qquad\text{and} \\
  \epsilon_{p_j}& =p_j-\tp_j^-
\end{align*}
for all $j \in \Z/r$. Most of these are visualized in figure \ref{pflowq} and $P_{q,q'}$ is the parallel transport used in Section~\ref{flowline}. Finally we define 
\begin{align} \label{Pdef}
  P=\max_j \norm{p_j}.
\end{align}

Since $\gamma_j\colon [0,\alpha_j] \to T^*N$ is a time $\alpha_j$ Hamiltonian flow curve for $H$ it is the obvious reparametrization of a time 1 flow curve for the Hamiltonian $\alpha_j H$. This reparametrization combined with similarly rescaling the Hamiltonian preserves action so it follows that
\begin{align} \label{Sreq}
  S_r(\arz) = \sum_{j\in \Z/r} ( \int_{\gamma_j}\lambda-Hdt) + \sum_{j\in
    \Z/r} p^-_j\epsilon_{q_j}.
\end{align}
\begin{Remark}
  \label{rem:2}
  This approximates the action in the following sense: we integrate the action integral over the small piece-wise flow curves, and then we add the symplectic area of the rectangles with corners $q_j^-,p_j^-,\tp_j^-,q_j$ in Figure~\ref{pflowq} to compensate for the fact that the pieces do not form a closed curve. Indeed, integrating the 1-form over a \emph{closed} curve seems reasonable.
\end{Remark}
We also see that $\gamma_{j-1}$ ends where $\gamma_j$ begins if and only if both $\epsilon_{q_j}$ and $\epsilon_{p_j}$ are 0. We have almost proved the following lemma. Note that we assumed Equation~\eqref{eq:10} in order to define $S_r$.

\begin{Lemma} \label{gradient}
  For any Hamiltonian $H$ and sub-division $\alpha$ (where $S_r$ is defined) we have
  \begin{align*}
    \quad \quad \quad \norm{\nabla_{q_j} S_r + \epsilon_{p_j}} \leq & C
    \max(1,P) (\norm{\epsilon_{q_j}} +
    \norm{\epsilon_{q_{j+1}}}) \\
    \norm{\nabla_{p_j} S_r - \te_{q_{j+1}}} \leq &
    \tfrac{1}{4}\norm{\epsilon_{q_{j+1}}}
    \quad \quad \qquad \quad \quad \qquad \qquad (=\tfrac{1}{4}\norm{\te_{q_{j+1}}})
  \end{align*}
  where $\nabla_{q_j} S_r \oplus \nabla_{p_j} S_r= \nabla_{z_j} S_r$ is the gradient with respect to the $j$th component in $T^*\Lre N$. Here $C$ is the constant from Lemma~\ref{gradient1half}.

  Furthermore, this implies that the critical points of $S_r$ are precisely those where all $\epsilon_{q_j}$ and $\epsilon_{p_j}$ are 0 such that the $\gamma_j$'s fit together to form a 1-periodic orbit of the Hamiltonian flow of $H$, and the critical value is the action of this orbit.
\end{Lemma}

\begin{proof}
  Adding the calculations of the gradients in Lemma~\ref{gradient1half} and Lemma~\ref{gradient1} (which holds on each segment because we did not define $S_r$ otherwise) and using $\max(1,\norm{p_j})\leq \max(1,P)$ proves the first part.
  
  For the second part we use
  \begin{itemize}
  \item the second inequality from the first part,
  \item $\norm{\epsilon_{q_{j+1}}}=\norm{\te_{q_{j+1}}}$ and
  \item $\nabla_{p_j} S_r = 0$
  \end{itemize}
  to conclude that any critical point must have $\epsilon_{q_j}=0$.

  Having this for all $j\in \Z/r$ we then use the first inequality together with $\nabla_{q_j} S_r = 0$ to conclude $\epsilon_{p_j}=0$. The fact that $S_r$ equals the action on these points follows from equation \eqref{Sreq}.
\end{proof}

We will need the following addition at a technical point later.

\begin{Corollary} \label{cor:tech}
  We have
  \begin{align*}
    \norm{\nabla_{q_j} S_r + \epsilon_{p_j}} \leq & \tfrac{1}{2} \max(1,P)
  \end{align*}
\end{Corollary}

\begin{proof}
  This follows easily from the first approximations in the lemma above and because we made $\delta_0$ less than $(8C)^{-1}$. Indeed, by construction we have $\norm{\epsilon_{q_j}} < 2\delta_0 < (4C)^{-1}$ for all $j\in \Z/r$.
\end{proof}

Lemma \ref{gradient} and the description of the 1-periodic orbits (in Section~\ref{action}) now imply that the set of critical points of $S_r$ is a compact set. The function $S_r$ with its gradient does not necessarily have index pairs, but following the idea of Viterbo we define a pseudo-gradient $X_r$ for which it does. On the set where $\max_j \norm{\epsilon_{q_j}} < \delta_0/10$ we use the gradient of $S_r$, and on the set $\max_j \norm{\epsilon_{q_j}} > \delta_0/5$ we keep the non-zero $\arp$-component of the gradient of $S_r$, but use $0$ as the $\arq$-component, i.e. on this set we have
\begin{align*}
  X_r = \bigoplus_j (0,\nabla_{p_j} S_r).
\end{align*}
In between we use some smooth convex combination of them. So by construction we have
\begin{align} \label{Xulig}
  X_r \cdot \nabla S_r \geq \norm{X_r}^2 \geq \sum_j \norm{\nabla_{p_j}S_r}^2 \geq \sum_j \tfrac{9}{16} \norm{\epsilon_{q_j}}^2,
\end{align}
and as we only made $X_r$ different from the gradient on a set where the latter is non-zero (last inequality uses Lemma \ref{gradient}) it is indeed a pseudo-gradient. To prove that Conley indices are well-defined we need the following lemma.

\begin{Lemma} \label{infiniteflow}
  When defined $(S_r,X_r)$ is CB.
\end{Lemma}

See Section~\ref{homoindex} for definition of CB.

\begin{proof}
  We start by proving that the flow of $X_r$ is defined for all times (positive and negative). By construction the flow preserves all the $q_j$ coordinates when $\epsilon_{q_j} > \delta_0/5$, which implies that it preserves the sets $\dist(q_j,q_{j+1})=k$ if $k > 3\delta_0/10$. Indeed, $\norm{\epsilon_{q_j}}>\dist(q_j,q_{j+1})-\delta_0/10$ because the length of any of the flow curves $\gamma_j$ is less than $\delta_0/10$ by Equation~\eqref{eq:10}. So we need only prove that none of the $p_j$ run of to $\infty$ in finite time. However, this follows because Lemma~\ref{gradient} implies
  \begin{align*}
    \norm{\nabla_{p_j} S_r} \leq 5\norm{\epsilon_{q_j}}/4 < 2\delta_0.
  \end{align*}

  Then we prove that $X(S_r)$ has a global lower bound on the complement of a compact set. First, we notice that we can extend the definition of $S_r$ to the set where $\dist(q_j,q_{j+1})$ are all allowed to be equal to $\delta_0$. This means we are done (using compactness of the complement) if we can prove a global lower bound on the set where $P>2$. Using Lemma~\ref{gradient} we see that
  \begin{align*}
    X(S_r) = X \cdot \nabla S_r \geq \norm{\nabla_{p_j} S_r}^2 \geq 9\norm{\epsilon_{q_j}}^2/16.
  \end{align*}
  This means that we can restrict to considering the points where $\max_j\norm{\epsilon_{q_j}}<\delta_0/10$ (any given constant), which means that we only need to consider the case where $X_r=\nabla S_r$.

  We will therefore need a lower bound on the norm squared of the gradient, but this is the same as having a lower bound on the norm. In fact, if we can find a lower bound on $G_q + G_p$ where
  \begin{align*}
    G_p=\sum_j\norm{\nabla_{p_j}S_r} \qquad \textrm{and} \qquad G_q=\sum_j\norm{\nabla_{q_j}S_r}
  \end{align*}
  then we are done. Define
  \begin{align*}
    L_q=\sum_j\norm{\epsilon_{q_j}} \qquad \textrm{and} \qquad
    L_p=\sum_j\norm{\epsilon_{p_j}}.
  \end{align*}
  Because of the approximation of $\nabla_{p_j}S_r$ in Lemma~\ref{gradient} we see that $G_p \geq L_q/2$ and hence
  \begin{align*}
    G_q + G_p \geq G_q + L_q/2,
  \end{align*}
  and we will prove the lemma by finding a lower bound on the latter. We will do this by finding $k_1,k_2>0$ and prove that if 
  \begin{align} \label{eq:16}
    L_q<k_1 \qquad \textrm{then} \qquad G_q>k_2.
  \end{align}

  Define $\unp=\min_j\norm{p_j}$. There are no 1-periodic flow curves on the compact set $1\leq \norm{p} \leq 2$ (this defines a compact set when combined with $\dist(q_j,q_{j+1})\leq \delta_0$), so there must exist $0<c<1$ such that $L_q+L_p>c$ for curves with all $z_j$'s contained in this set. Now we prove that if we define
  \begin{align*}
    k_1=\min(c/2,\frac{1}{4C},\frac{c}{8C}) 
  \end{align*}
  we can find $k_2$ such that the statement in Equation~\eqref{eq:16} is true. Here $C$ is the constant from Lemma \ref{gradient}.

  So, assume that $L_q$ is less than this $k_1$ then we divide the proof that $G_q$ is bounded from below by some $k_2$ into two cases.

  First case: $\unp<P/2$. By assumption we have some $j$ such that $\norm{p_j}=P \geq 2$ and for another $j'$ we have $\norm{p_{j'}}<P/2$. The ``curve'' $\arz$ has to move this distance in $p$-direction and back again. More precisely, the Hamiltonian flow of $H$ when $\norm{p}\geq 1$ is well-understood and we have
  \begin{align*}
    \absv*{\norm{p_i}-\norm{p_{i-1}}} = \absv*{\norm{p_i}-\norm{p_i^-}} < \norm{\epsilon_{p_i}}
  \end{align*}
  when $\norm{p_i}$ and $\norm{p_{i-1}}$ are greater than $1$. If \emph{precisely} one of them is less than $1$ then
  \begin{align*}
    \norm{p_i}-1 < \norm{p_i}-\norm{p_i^-} < \norm{\epsilon_{p_i}} \quad &\textrm{when } \norm{p_{i-1}} < 1 \\
    \norm{p_{i-1}}-1 < \norm{p^-_i}-\norm{p_i} < \norm{\epsilon_{p_i}} \quad &\textrm{when } 
    \norm{p_i} < 1.
  \end{align*}
  We see that for $\norm{p_j}\geq P$ to ``move'' all the way down to $\norm{p_{j'}}<P/2$ we must have
  \begin{align*}
    \sum_{j < i \leq j'} \norm{\epsilon_{p_i}} \geq P/2 \qquad \textrm{and} \qquad
    \sum_{j' < i \leq j} \norm{\epsilon_{p_i}} \geq P/2 
  \end{align*}
  and thus $L_p \geq P$. Note that with the cyclic ordering $j,j'\in \Z/r$ both sums makes sense. The approximation in Lemma \ref{gradient} and the bound $L_q<k_1 \leq 1/(4C)$ now gives
  \begin{align*}
    G_q & = \sum_j \norm{\nabla_{q_j} S_r} \\ &
    > \sum_j \pare*{\norm{\epsilon_{p_j}} -
    CP(\norm{\epsilon_{q_j}}+\norm{\epsilon_{q_{j+1}}})}
    > (P-\frac{P}{2}) \geq 1,
  \end{align*}
  which is a positive constant.

  The second case: $\unp\geq P/2$. In this case we can, because the
  flow is equivariant with respect to the $\R_+$ action on the set
  $\norm{p}\geq 1$, multiply our ``piecewise flow curve'' with
  $2/P$ to obtain a piecewise flow curve on the compact set $1\leq
  \norm{p} \leq 2$. This does not change any of the
  $\epsilon_{q_j}$'s, but it scales the $\epsilon_{p_j}$'s so we can
  conclude that the original curve satisfies
  \begin{align*}
    \frac{2}{P} L_p + L_q > c.
  \end{align*}
  Because $L_q<c/2$ this implies that $L_p > \frac{cP}{4}$, which
  implies by using the bound $L_q< k_1 \leq \frac{c}{8C}$ that
  \begin{align*}
    G_q     > \sum_j \norm{\epsilon_{p_j}} -
    CP(\norm{\epsilon_{q_j}}+\norm{\epsilon_{q_{j+1}}})  
    > \frac{cP}{4} - \frac{cP}{8} > \frac{cP}{4}.
  \end{align*}
  This is again a positive constant.
\end{proof}


\section{The Suspension Maps}\label{cha:suspmaps}

In this section we prove that when increasing $r$ by 1 we get a relative Thom space construction (defined in this section) of the Conley indices, whose existence was guaranteed in the previous section. The first part is producing an explicit (and hence canonical) map realizing this homotopy equivalence, this will respects quotients and inclusions of index pairs. The second part is a concrete construction of index pairs which proves that the maps induce a homotopy equivalence. We will need this map for any $H$ (linear at infinity with slope not a geodesic length) and the Conley index with respect to any fixed interval $[a,b]$ where $a$ and $b$ are regular values for the finite dimensional approximation $S_r$. The result of this section is summarized in Proposition~\ref{prop:Suspension:1}, and an even shorter summary is given by Equation~\eqref{eq:38}.

\subsection{Definitions and Preliminaries}

Define the relative Thom construction of a metric vector bundle $E\to M$ on a pair $(A,B)$ in $M$ by
\begin{align*}
  (A,B)^{E-} =  (DE_{\mid A}, UE_{\mid A} \cup DE_{\mid B}).
\end{align*}
Here $UE$ denotes the unit sphere bundle of the vector bundle $E$. We also define a shorthand for the quotient of the pair by
\begin{align*}
    (A,B)^{E/} =  DE_{\mid A}/\pare*{UE_{\mid A} \cup DE_{\mid B}}.
\end{align*}
We use these notations with $-$ and $/$ because we are dealing with pairs. However, as is standard we will use the notation $A^E=(A,\varnothing)^{E/}$ for the Thom-space (and sometimes even Thom-spectra, but we will make this clear from the context) when we are dealing with a single unbased space $A$. If the space is based at $*\in A$ we use $A^E = (A,\{*\})^{E/}$. This generalizes that $\Sigma A$ usually means two different things when $A$ is based and unbased - i.e. the usual suspension and the reduced suspension.

The short version of what we prove in this section is
\begin{align} \label{eq:38}
  I_a^b(S_{r+1},X_{r+1}) \simeq (A,B)^{TN/}
\end{align}
when $(A,B)$ is a good index pair for $(S_r,X_r)$. To make sense of $TN$ as a vector bundle on $T^*\Lre N$ we define
\begin{align*}
  \ev_j \co T^*\Lre N \to N
\end{align*}
to be the map given by $\ev_j(\arz)=q_j$. Then by abuse we could (and will in later sections) define $TN=\ev_0^*TN$ as a metric vector bundle over $T^*\Lre N$. However, for convenience we will alter the notation a bit in this section.

Indeed, we have thus far indexed the coordinates of a point in $T^*\Lambda_r N$ by $j\in\Z/r$. However since we are comparing this construction for different $r$'s this is inconvenient in this section. The main idea in this section is to insert an extra point somewhere in the ``cycle'' of points. This becomes notationally messy if we insert the point at $j=0$ or $j=r$ due to the reindexing combined with the change of relations in the groups ($r=0$ is changed to $r+1=0$). So, to make the argument more transparent we identify $\Z/(r+1)$ with $\{0,\dots,r\}$ and we identify $\Z/r$ with the sub\emph{set} not containing $j$ for $0 < j < r$. So the $j^\textrm{th}$ point in $T^*\Lambda_{r+1}N$ is the ``new'' point, and we fix this $j$ throughout the section. Since we have cyclic symmetry this covers all cases even $j=0$ and $j=r$. For $j=0$ this is like inserting an extra $z_0$, and is why we would define $TN \to T^* \Lre N$ as above, but in this section (to ease notation) we thus use $TN=\ev_{j+1}^*TN$.

Fix the Hamiltonian $H$ (linear at infinity with slope not a geodesic length). Then notice that: for any sub-divisions such that both $S_r$ and $S_{r+1}$ are defined their critical values coincide. Indeed, Lemma \ref{gradient} identifies the critical points as dissections of the 1-periodic Hamiltonian orbits, and the critical value is the action. However; if we use $\alpha_i=1/r$ and $\alpha'_i=1/(r+1)$ as sub-divisions when defining $S_r$ and $S_{r+1}$ respectively, then in the two cases this orbit is dissected in very different ways. This makes them difficult to compare. It is therefore convenient to \emph{not} use these standard choices. The more convenient choice is any choice of $\alpha=(\alpha_0,\dots,\alpha_r)$ where $\alpha_j=0$ then this works simultaneously for defining $S_r$ and $S_{r+1}$ (provided $\alpha_i\leq \tfrac{3}{(r+1)}$ as assumed in Equation~\eqref{eq:52}) using the indexing described above. This makes the Conley indices comparable. Indeed, the piece $\gamma_j$ is not used in the definition of $S_r$ and this fits well with the fact that $\gamma_j$ is constant and integrates the $0$ Hamiltonian in the definition of $S_{r+1}$.

Since we want the homotopy equivalence in Equation~\eqref{eq:38} to be a contractible choice (compatible with other structure) we will define it rather explicitly. Indeed, define the proper embedding
\begin{align} 
  h_0 \co TN \to T^*\Lambda_{r+1} N
\end{align}
by the simple formula
\begin{align} \label{eq:14}
  h_0(\arz,v) = (z_0,\dots,z_{j-1},(q_{j+1},v),z_{j+1},\dots,z_r)
\end{align}
for $\arz \in T^*\Lre N$ and $v\in TN_{\arz}$. That is, $h_0$ inserts the new point $z_j=(q_{j+1},v)$, which makes sense since by definition of $TN$ we have $v\in T_{q_{j+1}}N=T^*_{q_{j+1}}N$. In the old $\Z/r$ and $\Z/(r+1)$ notation and with $j=0$ (which is what we use outside of this section) this map is defined by
\begin{align}\label{eq:41}
  h_0(z_0,\dots,z_{r-1},v)=((q_0,v),z_0,\dots,z_{r-1}),
\end{align}
where $v\in T_{q_0} N$. The reason why this is notational messy is that we move the index on all the points up, which we avoid having to do in the notation in this section.

\begin{Lemma}\label{lem:commu}
  The diagram
  \begin{align} \label{comar}
    \xymatrix{
      TN \ar[r]^{h_0} \ar[d] & T^* \Lambda_{r+1}N
      \ar[d]^{S_{r+1}} \\
      T^* \Lre N \ar[r]^{S_r} & \R
    }
  \end{align}
  commutes.
\end{Lemma}

\begin{proof}
  In the definition of $S_{r+1}$ we have $\gamma_j$ is constant (in fact - it is parametrized by a point $\{0\}$). This implies that $(q_{j+1}^-,p_{j+1}^-)=(q_j,p_j)$. On the image of $h_0$ we also have $q_j = q_{j+1}$ implying that $\epsilon_{q_{j+1}}=0$ and hence both of the two new terms
  \begin{align*}
    \int_{\gamma_j} \lambda - Hdt \qquad \textrm{and} \qquad p_{j+1}^-\epsilon_{q_{j+1}}
  \end{align*}
  in $S_{r+1}$ (Equation~\eqref{Sreq}) are 0 - independent of $p_{j+1}^-=p_j=v$.
\end{proof}

If $(A_r,B_r)$ is a good index pair for $S_r$ then this lemma implies that 
\begin{align}
  S_{r+1}(h_0(D_RTN_{\mid B_r})) \subset \{a\} \qquad &\textrm{using} \qquad S_r(B_r)\subset \{a\} \textrm{ and} \\
  S_{r+1}(h_0(U_RTN_{\mid A_r})) \subset [a,b] \qquad &\textrm{using} \qquad S_r(A_r)\subset [a,b].
\end{align}
The goal is to get an induced map from the relative Thom construction pair. However, for this we will need to the sphere bundle to be mapped to points where $S_{r+1}$ takes values less than or equal to $a$, which by the above formula it is not, and for this we need the $R$ factor to choose the discs big enough. We therefore pre-compose with the canonical homeomorphism of the pair
\begin{align} \label{eq:39}
  (DTN,UTN) \cong (D_RTN,U_RTN)
\end{align}
given by scaling with $R>0$. This does not change the fact that the unit sphere is not mapped to values less than $a$, but to get this and induce maps on index pairs we also modify $h_0$ by using the negative pseudo-gradient flow. So, define
\begin{align*}
  h_t \co TN \to T^*\Lambda_{r+1} N
\end{align*}
by $h_0$ composed with the flow of $-X_{r+1}$ for time $t$.

\begin{Lemma}\label{lem:flowbelow}
  There exist $t_0>0$ and $R_0>0$ such that for $t\geq t_0$ and $R\geq R_0$ we have
  \begin{align*}
    h_t(U_RTN_{\mid A_r}) \subset S_{r+1}^{-1}(]-\infty,a]).
  \end{align*}
\end{Lemma}

\begin{proof}
  Let $k>0$ be such that Lemma~\ref{infiniteflow} gives that
  \begin{align*}
    X_{r+1}(S_{r+1}) \geq k^{-1}
  \end{align*}
  on the complement of a compact set. Now as in the proof of Lemma~\ref{lem:Homoindex:1} we see that for $t_0 > (b-a)k$ only a compact subset of $(S_{r+1})^{-1}([a,b])$ is not flowed to having the value of $S_{r+1}$ less than $a$. Since $h_0$ is proper this implies that we may find $R_0>0$ large enough such that $h_0^{-1}$ of this compact set is in the interior of $D_{R_0}TN$.
\end{proof}

These lemmas now imply that we get an induced map of pairs
\begin{align*}
  h_t \co (D_RTN_{\mid A_r},U_RTN_{\mid B_r} \cup D_RTN_{\mid A_r}) \to (S_{r+1}^{-1}(]-\infty,b]),S_{r+1}^{-1}(]-\infty,a]))
\end{align*}
for $t>t_0$ and $R>R_0$. Furthermore, if $(A_{r+1},B_{r+1})$ is any good index pair for $(S_{r+1},X_{r+1})$ then by flowing further (see Section~\ref{homoindex} in particular Lemma~\ref{lem:Homoindex:1}) we can get that
\begin{align*}
  h_t(D_RTN) \subset A_{r+1} \cup S_{r+1}^{-1}(]-\infty,a]).  
\end{align*}
This is enough to get maps induced on the quotients to $A_{r+1}/B_{r+1}$. Indeed, since the pair is good we have
\begin{align} \label{eq:40}
  \pare*{A_{r+1} \cup S_{r+1}^{-1}(]-\infty,a])}/S_{r+1}^{-1}(]-\infty,a]) \cong A_{r+1}/B_{r+1}.
\end{align}
We thus define the induced composition
\begin{align*}
  \th_t \co (A,B)^{TN/} \to A_{r+1}/B_{r+1}
\end{align*}
for $t$ and $R$ large. Here we are pre-composing $h_t$ with the homeomorphism in Equation~\eqref{eq:39} and post-composing with the identification in Equation~\eqref{eq:40}.

The rest of this section is devoted to proving the following Proposition.
\begin{Proposition}   \label{prop:Suspension:1}
  Let $(A_r,B_r)$ be a good index pair for $(S_r,X_r)$, and $(A_{r+1},B_{r+1})$ be a good index pair for $(S_{r+1},X_{r+1})$ (both defined using the same $H$ and a compatible subdivision $\alpha$ as above). Then for large $t$ and $R$ we have that the induced map
  \begin{align*}
    \th_t \co (A,B)^{TN/} \to  A_{r+1}/B_{r+1}
  \end{align*}
  is a homotopy equivalence.
\end{Proposition}

We will not explicitly use the following, but when $N$ is oriented we have a Thom isomorphism
\begin{align*}
  H_*(A,B) \cong H_{*+d}((A,B)^{TN-})
\end{align*}
However, this means that the result in this section implies that: the Morse homology (shifted in degree by $rd$) of $S_r$ does not depend on the choice of $r$ (see \cite{MR1045282} and Appendix~\ref{cha:app} for more on the relation between Morse homology and the homology of index pairs).

\subsection{The homotopy type of $A_{r+1}/B_{r+1}$}

In this subsection we construct another CB pseudo-gradient $Z_{r+1}$ for $S_{r+1}$ and an index pairs $(A',B')$ for $(S_{r+1},Z_{r+1})$. These will help in proving Proposition~\ref{prop:Suspension:1}.

Consider the subspace
\begin{align*}
  O \subset T^*\Lambda_{r+1} N
\end{align*}
given by the equation $\dist(q_{j-1},q_{j+1}) < \delta_0$. There is a canonical projection
\begin{align*}
  \pi \co O \to T^*\Lre N
\end{align*}
given by forgetting $z_j$. This makes $O$ a fiber-bundle with contractible fibers. Indeed, the fiber $O_{\arz} = \pi^{-1}(\arz)$ is symplectomorphic to $T^*U$ where
\begin{align*}
  U = \{q_j \in N \mid \dist(q_j,q_{j-1}) < \delta_0 \textrm{ and } \dist(q_j,q_{j+1}) < \delta_0 \},
\end{align*}
which is convex and non-empty for all $\arz=(z_0,\dots,z_{j-1},z_{j+1},\dots,z_r) \in T^*\Lre N$.

Now consider the function $S_{r+1}'=S_{r+1\mid O}$ and its restricted pseudo-gradient $X_{r+1}'=X_{r+1\mid O}$.
\begin{Lemma}
  \label{lem:Suspension:1}
  The restricted function and pseudo-gradient $(S_{r+1}',X_{r+1}')$ is CB
\end{Lemma}

Note that this lemma together with Lemma~\ref{hominv} and Lemma~\ref{lem:Homoindex:2} shows that we may as well replace $I_a^b(S_{r+1},X_{r+1})$ with $I_a^b(S_{r+1}',X_{r+1}')$. Indeed, we have a canonical homotopy equivalence from the latter to the former.

\begin{proof}
  We already know that $(S_{r+1},X_{r+1})$ is CB. So all we need to prove is that the flow is defined for all time on the restriction. This is the case because the boundary is given by
  \begin{align*}
    \dist(q_{j-1},q_{j+1}) = \delta_0
  \end{align*}
  and the flow of $X_{r+1}$ preserves this equation. Indeed, by definition $X_{r+1}$ preserves all $q_i$ when $\max_i \epsilon_{q_i} \geq \delta_0/5$ and the above equation implies this since
  \begin{align*}
    \dist(q_{j-1},q_{j+1}) \leq \dist(q_{j-1},q_j^-) + \dist(q_j^-,q_j) + \dist(q_j,q_{j+1}^-) + \dist(q_{j+1}^-,q_{j+1}),
  \end{align*}
  and by construction we have
  \begin{itemize}
  \item $\dist(q_{j-1},q_j^-) \leq \delta_0/10$ by Equation~\eqref{eq:10},
  \item $\dist(q_j^-,q_j) = \norm{\epsilon_{q_j}}$,
  \item $q_{j+1}^- = q_j$ (the curve $\gamma_j$ is constant) and
  \item $\dist(q_{j+1}^-,q_{j+1}) = \norm{\epsilon_{q_{j+1}}}$.
  \end{itemize}
  These imply that $\norm{\epsilon_{q_{j+1}}} + \norm{\epsilon_{q_j}} \geq 9\delta_0/10$.
\end{proof}

We now consider how $S_{r+1}'$ looks fiber-wise over a point $\arz \in \Lre N$. Indeed, we define
\begin{align*}
  S^{\arz} = S_{r+1\mid O_{\arz}} \co O_{\arz} \to \R
\end{align*}
as a function depending on coordinates $z_j=(q_j,p_j)$ such that $\dist(q_{j-1},q_j)<\delta_0$ and $\dist(q_j,q_{j+1})<\delta_0$.

\begin{Lemma} \label{lem:Suspension:3}
  There exists a smooth section $s\co T^*\Lre N \to O$ such that $s(\arz)$ is $\underline{\textrm{the}}$ unique critical point for $S^\arz$. Furthermore, this unique critical point is non-degenerate and has $q_j=q_{j+1}$ (i.e. lies in the image of $h_0$), and the Hessian of $S^\arz$ at this point is given by
  \begin{align*}
    \left[\begin{array}{cc} 
      Q & -I \\
      -I & 0 \end{array} 
    \right]
  \end{align*}
  in the $(q_j,p_j)$ splitting. Hence the vector bundle of positive (resp. negative) eigen\-spaces of this Hessian over $T^* \Lre N$ are both canonically isomorphic to $TN=\ev_{j+1}^*TN$.
\end{Lemma}

Note that, here $Q$ is an arbitrary symmetric bilinear form on $T_{q_j}N$, and the definition of $-I$ formally uses the canonical pairing of tangent with cotangent vectors.

\begin{proof}
  Firstly we note that by the inequality in Lemma~\ref{gradient} and because
  \begin{align*}
    \nabla_{p_j} S^{\arz} = \nabla_{p_j} S_{r+1}
  \end{align*}
  we have this component of the gradient is equal to zero if and only if $q_j=q_{j+1}$.

  Since the only term in Equation~\eqref{Sreq} for $S_{r+1}$ which involves $p_j$ is
  \begin{align} \label{eq:36}
    p_{j+1}^-\epsilon_{q_{j+1}} = p_j\exp_{q_j}^{-1}(q_{j+1})
  \end{align}
  it follows that
  \begin{align} \label{eq:35}
    \nabla_{q_j} S^{\arz} = p - p_j \qquad \textrm{on the set where } q_j=q_{j+1}
  \end{align}
  for some fixed $p$ (independent on $p_j$) equal to the gradient of the remaining terms in $S_{r+1}$ with respect to $q_j$. It follows that $(q_j,p_j)$ is a critical point for $S^{\arz}$ if and only if $(q_j,p_j)=(q_{j+1},p)$, which defines the smooth section $s$.

  The fact that the Hessian has a zero matrix in its bottom right follows from $S^{\arz}$ not depending on $p_j$ (when $q_j=q_{j+1}$). The two copies of $-I$ in the Hessian follows from Equation~\eqref{eq:35}.

  The last statement follows from the fact that for $v\in T_{q_{j+1}}N$ this Hessian is positive on vectors of the type $(v,-kv)$ for large $k>>0$ and negative on $(v,kv)$ for large $k>>0$.
\end{proof}

We wish to combine a construction of index pairs for $S^{\arz}$ for each $\arz$ and an index pair for $(S_r,X_r)$ to get an index pair for $(S_{r+1}',X_{r+1}')$. That is, we will combine index pair for the base of $\pi$ with fiber-wise index pairs to get a global index pair. So, we now explicitly construct CB pseudo-gradients for $S^{\arz}$. We do this in a way smoothly dependent on $\arz$.

For $\arz \in \Lre N$ define $\delta(\arz) = \delta_0 - \dist(q_{j-1},q_{j+1}) > 0$. Define the pseudo-gradient $X^{\arz}$ by
\begin{align*}
  X^{\arz} = (g(\arz,q_j,p_j)\nabla_{q_j} S^{\arz} ,\nabla_{p_j} S^{\arz})
\end{align*}
where $g\co O \to [0,1]$ is smooth such that $g=1$ when 
\begin{align} \label{eq:34}
  \dist(q_j,q_{j+1}) < \frac{\delta(\arz)}{2}
\end{align}
and $g=0$ when
\begin{align*}
  \dist(q_j,q_{j+1}) > \tfrac{2}{3}\delta(\arz) < \delta_0.
\end{align*}

\begin{Lemma}
  \label{lem:Suspension:2}
  For any $\arz \in \Lre N$ we have that the function and pseudo-gradient $(S^{\arz},X^{\arz})$ is CB.
\end{Lemma}

\begin{proof}
  Firstly $X^{\arz}$ is a pseudo-gradient since we only changed the $q_j$-component when $q_j\neq q_{j+1}$ which implies by Lemma~\ref{gradient} that $\nabla_{p_j}S^{\arz} \neq 0$.

  Furthermore, we made sure to define $X^{\arz}$ such that its flow preserves the set where $\dist (q_j,q_{j+1})=\delta_0$, but also the set where $\dist(q_{j-1},q_j)=\delta_0$. Indeed, if $\dist(q_{j-1},q_j)=\delta_0$ then we see that
  \begin{align*}
    \dist(q_j,q_{j+1}) \geq \dist(q_{j-1},q_j) - \dist(q_{j-1},q_{j+1}) = \delta_0- \dist(q_{j-1},q_{j+1})  \geq \delta(\arz)
  \end{align*}
  and hence $g=0$. Since Lemma~\ref{gradient} shows that the norm of $\nabla_{p_j} S^{\arz}$ is bounded we get that the flow of $X^{\arz}$ is defined for all times (positive and negative).

  Lemma~\ref{gradient} also shows that when $\dist(q_j,q_{j+1}) \geq \tfrac{1}{2}\delta(\arz)$ then we have
  \begin{align*}
    X^{\arz}(S^{\arz}) \geq \norm{\nabla_{p_j} S^{\arz}}^2 \geq \tfrac{9}{16} \norm{\epsilon_{q_{j+1}}} = \tfrac{9}{16}\dist(q_j,q_{j+1}) \geq \tfrac{9}{32}\delta(\arz).
  \end{align*}
  Hence we can assume, when proving a lower bound on the compliment of a compact set, that $g=1$. This implies that $X^\arz=\nabla S^\arz$, and that we may assume that $\norm{p_j} > R$ for some $R$, which we may pick much larger than all the fixed numbers $3\norm{p_i},i\neq j$. Now, Corollary~\ref{cor:tech} provides the lower bound:
  \begin{align*}
    \norm{\nabla_{q_j}S^{\arz}} \geq \norm{\epsilon_{p_j}} - \tfrac12 \norm{p_j} \geq \tfrac 12 \norm{p_j} - \norm{p_{j+1}} \geq \tfrac{1}{2}R -\tfrac{1}{3}R  \geq \tfrac{1}{6}R.
  \end{align*}
  Indeed, $\norm{\epsilon_{p_j}} \geq \norm{\tp_{j+1}}-\norm{p_{j+1}} = \norm{p_{j+1}^-}-\norm{p_{j+1}}= \norm{p_j}-\norm{p_{j+1}}$, the latter since $\gamma_j$ is constant.
\end{proof}

We now define some explicit index pairs for $(S^\arz,X^\arz)$. Indeed, for each  $\epsilon>0$ and $c>0$ and with $s$ the section in Lemma~\ref{lem:Suspension:3} we define
\begin{align*}
  A^{\arz}_{\epsilon} = \exp_{s(\arz)} (D_{c\epsilon} F_- \times D_\epsilon F_+) \subset O_{\arz}\\
  B^{\arz}_{\epsilon} = \exp_{s(\arz)} (U_{c\epsilon} F_- \times D_\epsilon F_+) \subset A^{\arz}_{\epsilon}.
\end{align*}
Here $F_- \oplus F_+$ is the canonical orthogonal decomposition of the Hessian of $S^{\arz}$ at $s(\arz)$ into negative and positive eigenspaces (using the Riemannian structure). Note that this is well-defined due to Lemma~\ref{lem:Suspension:3}. The $c$ will not be useful nor relevant until the next subsection - so we fix this $c$ as a positive constant in the rest of this section.

\begin{Lemma}
  \label{lem:Suspension:4}
  The pair $(A^{\arz}_{\epsilon},B^{\arz}_{\epsilon})$ is for small $\epsilon$ an index pair for $(S^{\arz},X^{\arz})$.
\end{Lemma}

\begin{proof}
  By Lemma~\ref{lem:Suspension:3} we see that $\dim(F_-)=\dim(F_+)=d$. So pick a linear isometry $\psi \co \R^{2d} \to T_{s(\arz)} O_{\arz}$ such that the usual inclusion $\R^d \subset \R^{2d}$ is identified with $F_-$ - hence the orthogonal complement of $\R^d\subset \R^{2d}$ is identified with $F_+$. We may even assume that the standard coordinate axes are mapped to eigenvectors. For some small $e>0$ we use this to define a normal coordinate chart
  \begin{align*}
    \exp_{s(\arz)} \circ \psi \co D_e^{2d} \to O_{\arz}.
  \end{align*}
  For small $\epsilon>0$ the pair $(A^{\arz}_{\epsilon},B^{\arz}_{\epsilon})$ is identified in this chart with
  \begin{align} \label{eq:37}
    (A_\epsilon,B_\epsilon)=(D_{c\epsilon}^d \times D_\epsilon^d,U_{c\epsilon}^d \times D_\epsilon^d)
  \end{align}
  Since it is a normal chart the pull back of the Hessian of $S^{\arz}$ at $s(\arz)$ equals the Hessian of the pull back function at $0$. This implies that the pull back of $X^{\arz}$ (which is the gradient of $S^{\arz}$) is equal to the usual gradient (in $\R^{2d}$) of the pull back function to the first order at $0$. We may thus reduced the lemma to proving that the pair in Equation~\eqref{eq:37} is for small $\epsilon>0$ an index pair for any function
  \begin{align*}
    f(x) = f(0) + \sum_{i=1}^{2d} \lambda_i x_i^2 + O(\norm{x}^3)
  \end{align*}
  (Restricted to a small enough neighborhood around $0$) with a pseudo-gradient
  \begin{align*}
    X_x = (2 \lambda_i x_i,\dots, 2\lambda_{2d} x_{2d}) + O(\norm{x}^2),
  \end{align*}
  where $\lambda_i < 0$ for $i=1,\dots,n$ and $\lambda_i >0$ for $i=n+1,\dots,2d$. We may assume that $\lambda_1$ is the negative eigenvalue closest to $0$ and that $\lambda_{2d}$ is the positive eigenvalue closest to $0$.

  To prove this consider the two functions
  \begin{align*}
    f_1(x) = \sum_{i=1}^d x_i^2 \qquad \textrm{and} \qquad f_2(x) =\sum_{i=d+1}^{2d} x_i^2.
  \end{align*}
  The change of these when flowing with $-X$ is given close to $0$ by
  \begin{align*}
    -X(f_1)(x) = -X_x \cdot \nabla f_1 = -\sum_{i=1}^d 4 \lambda_i x_i^2 + O(\norm{x}^3) \geq - 4 \lambda_1 f_1(x) + O(\norm{x}^3)
  \end{align*}
  and
  \begin{align*}
    -X(f_2)(x) = -X_x \cdot \nabla f_2 = -\sum_{i=d+1}^{2d} 4 \lambda_ i x_i^2 + O(\norm{x}^3) \leq -4\lambda_{2d} f_2(x) + O(\norm{x}^3).
  \end{align*}
  The boundary of $A_\epsilon$ has two parts (not disjoint) given by:
  \begin{align*}
    B_\epsilon = \{f_1=c^2\epsilon^2, f_2\leq \epsilon^2\} \quad \textrm{and} \quad
    W_\epsilon = \{f_1\leq c^2\epsilon^2, f_2 = \epsilon^2\}    
  \end{align*}
  For $\epsilon$ small enough $-X(f_1)$ is positive on $B_\epsilon$ and $-X(f_2)$ is negative on $W_\epsilon$.

  It follows that for such $\epsilon$ the flow does, indeed, exit through $B_\epsilon$, and by making $\epsilon$ even smaller the flow will not reenter before the value of $f$ becomes to low for it to ever reenter.
\end{proof}

\begin{Lemma}
  \label{Lem:homtyp}
  Let $(A,B)$ be a good index pair for $(S_r,X_r)$ then there exists an $\epsilon>0$ small enough and a CB pseudo-gradient $Z_{r+1}$ for $S_{r+1}$ such that
  \begin{align*}
    (A',B') = (\bigcup_{\arz \in A} A^\arz_\epsilon, \bigcup_{\arz \in A} B_\epsilon^\arz \cup \bigcup_{\arz\in B} A_\epsilon^\arz)
  \end{align*}
  is an index pair for $(S_{r+1},Z_{r+1})$.
\end{Lemma}

\begin{proof}
  We still only focus on the set $O$. Define for any set $C\subset T^* \Lre N$ and $\epsilon>0$ the sets
  \begin{align*}
    A_C' = \cup_{\arz\in C} A_{\epsilon}^{\arz} \qquad \textrm{and} \qquad B_C' = \cup_{\arz \in C} B_\epsilon^\arz.
  \end{align*}
  Then $(A',B')=(A_A',B_A'\cup A_B')$.

  Pick a neighborhood $W$ with compact closure around $A$ (not containing any additional critical points).  Let $\delta>0$ be such that $S_r$ evaluated on the flow of $-X_r$ on $B$ goes below $a-\delta$ before the flow exits $W$. Since the section in Lemma~\ref{lem:Suspension:3} is in the image of $h_0$ Lemma~\ref{lem:commu} tells us that $S_{r+1}(s(\arz))=S_r(\arz)$ and thus for small enough $\epsilon$ the diagram in Lemma~\ref{lem:Suspension:3} implies that we can assume that:
  \begin{align}
    \label{eq:2dd}
    S_{r+1}(A_\epsilon^\arz)\subset [S_r(\arz)-\tfrac{\delta}{3},S_r(\arz)+\tfrac{\delta}{3}] \textrm{ for all } \arz \in \ovl{W},
  \end{align}
  which we will use at the very end of the proof. Notice that $A'_{\ovl{W}}$ is compact.

  Now, pick a neighborhood $U$ around $\ovl{W}$ also with compact closure (and with no new critical points). Let $g\co T^*\Lre N \to I$ be smooth and such that $\{g=1\}=\ovl{W}$ and $\{g>0\}=U$. Now define the vector field $Y$ by
  \begin{align*}
    Y_x = (1-g(\pi(x)))\cdot (X_{r+1})_x + g(\pi(x)) \cdot (X^{\pi(x)})_x. 
  \end{align*}
  Here $X^{\pi(x)}$ is the fiber-wise pseudo-gradient for $S^{\pi(x)}$ from above - viewed as a vector field on all of $O$ parallel to the fibers of $\pi$. So this is the fiber-wise pseudo-gradient precisely on $\pi^{-1}({\overline{W}})$, and on the ``buffer'' $\pi^{-1}(\overline{W}-\overline{U})$ it convexly interpolates to $X_{r+1}$ so that outside it is $X_{r+1}$. This vector field $Y$ is \emph{not} a pseudo-gradient because by construction $Y=0$ at all points $x\in \im(s)\cap A_{\ovl{W}}$. However, it does satisfy
  \begin{align*}
    Y(S_{r+1}) =& (1-g)X_{r+1}(S_{r+1})+gX^{\pi(-)}(S_{r+1}) = \\
    = &(1-g)X_{r+1}(S_{r+1})+gX^{\pi(-)}(S^{\pi(-)}) \geq 0.
  \end{align*}
  A more serious, but related, problem is that the exit set of $A'$ for $Y$ is not $B'$. Indeed, the exit set is precisely $B'_A$ since we are flowing fiber-wise over $A$, and these are the exit sets for $X^{\arz}$. So, we need to modify $Y$ a bit.

  Let $Y'$ be a vector field on $O$ which is a lift of $X_r$. I.e. $\pi_*(Y')=X_r$ at all points in $O$. By Lemma~\ref{lem:Suspension:3} and Lemma~\ref{lem:commu} we have that $Y'(S_{r+1})=X_r(S_r)$ on the image of the section $s$. Indeed, this uses that in the fibers $s(\arz)$ is a critical point for $S^\arz$. This means that the open set $\{Y'(S_{r+1})>0\}$ contains all the non-critical points where $Y=0$. Since $X_r(S_r)>0$ on the compact set $\partial A$ we have for small $\epsilon>0$ the inclusion 
  \begin{align} \label{eq:42}
    A_{\partial A}' \subset \{Y'(S_{r+1})>0\}.
  \end{align}
  Since we only wish to change $Y$ on a compact set we pick a smooth function $f \co O \to I$ such that
  \begin{align*}
    \{f>0\} = \{Y'(S_{r+1})>0\} \cap \inte(A'_U)
  \end{align*}
  Here $\inte(-)$ means interior. Note that this depends on $\epsilon>0$ and that this set still contains all the non-critical points of $S_{r+1}$ where $Y=0$. We now define
  \begin{align*}
    Z_{r+1} = Y + f \cdot Y'=(1-g)X_{r+1} + gX^{\pi(-)} + f Y',
  \end{align*}
  which is now a pseudo-gradient for $S_{r+1}$.

  Note that, for $x\in \pi^{-1}(\ovl{W})$ we have (since $g=0$) that 
  \begin{align*}
    \pi_*((Z_{r+1})_x)=\pi_*(X^{\pi(x)}_x+f(x)Z_x)=f(x)(X_r)_{\pi(x)},
  \end{align*}
  which implies that the projection under $\pi$ of flow lines (over $\ovl{W}$) of $-Z_{r+1}$ are reparameterizations (may even stand still if $f(x)=0$) of the flow of $-X_r$. We now argue that this implies that $(A',B')$ is an index pair.

  Indeed, any $x\in A'$ where the flow of $-Z_{r+1}$ exits $A'$ must be in the boundary of $A'$ which is given by
  \begin{align*}
    \partial A' = \pare*{\bigcup_{\arz\in \partial A} A_\epsilon^\arz }\cup \pare*{\bigcup_{\arz \in A} \partial A_{\epsilon}^\arz}.
  \end{align*}
  Now if $x$ is in the latter of these two unions then by construction $Z_{r+1}=X^{\pi(-)}$ (indeed, $f=0$ and $g=1$), and we see that $x$ is an exit point if and only if $x\in B_\epsilon^{\pi(x)}$. Indeed, at this point (and outside of $A_\epsilon^\arz$) we have $Y=X^{\pi(x)}$ the fiber-wise CB pseudo-gradient (Lemma~\ref{lem:Suspension:2}) and $(A_\epsilon^{\pi(x)},B_\epsilon^{\pi(x)})$ an index pair in this fiber. This implies that the flow of $x$ stays in this fiber and outside of $A_\epsilon^\arz$ where $Z_{r+1}=X^{\pi(x)}$, so we conclude that the flow never returns.

  If on the other hand $x$ is an exit point in the first union but \emph{not} in the second then we see that $x\in A_{\partial A}' \cap \inte(A_U')$ and hence $f>0$ (this was why we wanted Equation~\eqref{eq:42}). Thus the flow of $-Z_{r+1}$ applied to $x$ projects under $\pi$ to a non-zero reparameterized flow line for $-X_r$. So it immediately exits $A'$ if and only if $x\in \pi^{-1}(B)$. When this is the case, the flow of $x$ never reenters. Indeed, the flow exists $A_{\ovl{W}}'$ either through a point in $B'_{\ovl{W}}$ in which case it never returns for the same reason as the first case - or the projected flow exits $W$, which by Equation~\eqref{eq:2dd} means that the value of $S_{r+1}$ is to low for it to ever reenter $A'$.

  Proof that $Z_{r+1}$ is CB (on $O$). This follows precisely as the proof of the fact that $X_{r+1}$ and $X^{\arz}$ are CB. Indeed, both preserve the sets where $\epsilon_{q_i}=k$ when this is close to the ``boundary'' of $O$ and both have $p_i$ components bounded - so the same is true for the convex combination defining $Y$ (and hence $Z_{r+1}$ outside a compact set). Also, a lower bound on $Z_{r+1}(S_{r+1})$ outside a compact set follows from the fact that we have this for $X^{\arz}$ (over the \emph{compact} set $\ovl{U}$ by Lemma~\ref{lem:Suspension:2}) and $X_{r+1}$ by Lemma~\ref{infiniteflow}.

  Formally we have only defined $Z_{r+1}$ on $O$, but we extend it using the following idea. By the fact that $Z_{r+1}$ is CB we see that the union of the flow of $-Z_{r+1}$ over all non-negative time of $B'$ intersected with $S_{r+1}^{-1}([a,b])$ is compact in $O$. Hence we can change $Z_{r+1}$ outside a compact set without changing that $(A',B')$ is an index pair. Now outside this compact set we interpolate between $Z_{r+1}$ and $X_{r+1}$ and extend by putting $Z_{r+1}=X_{r+1}$ outside of a slightly larger compact set and outside of $O$. This is CB because it equals $X_{r+1}$ outside a compact set.
\end{proof}

\begin{proof}[Proof of Proposition~\ref{prop:Suspension:1}]
  Firstly, we analyze the index pair $(A_\epsilon^\arz,B_\epsilon^\arz)$ for small $c>0$ (recall that it in fact was made to depend on a constant $c$). Let $H^\arz$ be the Hessian from Lemma~\ref{lem:Suspension:3} (at $s(\arz)$). Since $(0,v)H^\arz(0,v)^T=0$ we see that the vertical vectors (fiber directions of $T^*N$) forms a complement to the positive eigenspace of $H^\arz$. For $c>0$ small the pair is a thin ``tube'' around $F_+$ (the positive eigenspace of $H^{\arz}$), and as Figure~\ref{fig:fiberint} 
  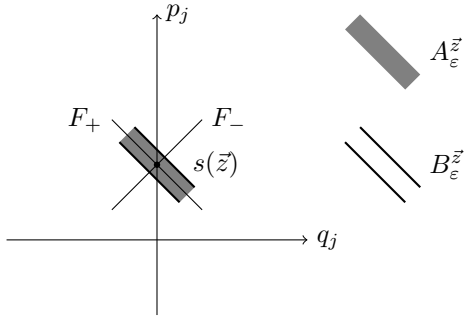
\begin{figure}[ht]
    \centering
    \begin{tikzpicture}
      \fill[gray] (-0.5,1.3) -- (-0.3,1.5) -- (0.5,0.7) -- (0.3,0.5) -- cycle;
      \draw[->] (0,-1) -- (0,3) node [right] {$p_j$};
      \draw[->] (-2,0) -- (2,0) node [right] {$q_j$};
      \fill (0,1) circle (1.2pt) node [right] {$\quad s(\arz)$};
      \draw (-0.6,0.4) -- (0.6,1.6) node [right] {$F_-$};
      \draw (0.6,0.4) -- (-0.6,1.6) node [left] {$F_+$};
      \draw[thick] (-0.5,1.3) -- (0.3,0.5);
      \draw[thick] (-0.3,1.5) -- (0.5,0.7);
      \fill[gray] (2.5,2.8)  -- (2.7,3.0) -- (3.5,2.2) -- (3.3,2.0) -- cycle;
      \draw[thick] (2.5,1.3) -- (3.3,0.5);
      \draw[thick] (2.7,1.5) -- (3.5,0.7);
      \draw (3.5,2.5) node [right] {$A_\epsilon^\arz$};
      \draw (3.5,1.0) node [right] {$B_\epsilon^\arz$};
    \end{tikzpicture}
    \caption{The index pair $(A_\epsilon^\arz,B_\epsilon^\arz)$ and the intersection with the fiber.}
    \label{fig:fiberint}
  \end{figure}
  illustrates we see that this means that the image of $h_0$ intersected with the pair is in fact a small disc and its sphere. Call these $(DTN^\arz,UTN^\arz)$. Now the pair $(A_\epsilon^\arz,B_\epsilon^\arz)$ deformation retracts onto this pair. Indeed, in the normal neighborhood in the proof of Lemma~\ref{lem:Suspension:4} all these sets are convex in $\R^{2d}$ and the convex deformation retraction parallel to the positive eigenspace $F_+$ is such a deformation retraction.
  
  Let $(A,B)$ be a good index pair for $(S_r,X_r)$, and pick $c>0$ so small that this works for all $\arz \in A$, and then let $\epsilon>0$, $(A',B')$, and $Z_{r+1}$ be as in Lemma~\ref{Lem:homtyp} above. The deformation retraction above is smoothly dependent on $\arz$ and thus defines a deformation retraction of $(A',B')$ onto a pair homeomorphic to $(A,B)^{TN-}$. We thus see that the map
  \begin{align*}
    (A,B)^{TN/} \to A'/B'
  \end{align*}
  given by using Equation~\eqref{eq:39} for large $R$ and by collapsing everything outside the interior of the small discs $DTN^\arz\subset D_RTN_{\arz}$ to the base-point $[B']$ is a homotopy equivalence.
  
  Since $Z_{r+1}$ is CB there is a $T>>0$ such that when we flow on $B'$ using $-Z_{r+1}$ for time $T$ we get the value of $S_{r+1}$ below $a$. It follows that we get an induced homotopy equivalence:
  \begin{align*}
    (A,B)^{TN/} \to S_{r+1}^{-1}((-\infty,b])/    S_{r+1}^{-1}((-\infty,a]).
  \end{align*}
  
  We can interpolate between the flow of $-Z_{r+1}$ and $-X_{r+1}$ in the following way. For all $s\in I$ the pseudo-gradient $X^s=sX_{r+1}+(1-s)X_{r+1}$ is CB (CB pseudo-gradients form a convex set) and there is a compact set $K\subset T^*\Lambda_{r+1} N$ and a $k>0$ such that $X^s(S_{r+1})>k$ for all $s$. Hence by properness of $h_0$ there is a large $R>0$ and a $T>>0$ such that the flow of $-X^s$ for time $T$ induces a homotopy of maps
  \begin{align*}
    g_S \co (A,B)^{TN/} \to S_{r+1}^{-1}((-\infty,b])/    S_{r+1}^{-1}((-\infty,a]).
  \end{align*}
  For $s=0$ this is the map argued to be a homotopy equivalence above, and for $s=1$ this is the map $h_T$ in Proposition~\ref{prop:Suspension:1}.
\end{proof}


\section{The Generating Function Spectrum} \label{sec:gener-funct-spectr}

In this section we construct a spectrum ``representing'' Floer homology. Concisely, we let $H\co T^*N \to \R$ be a Hamiltonian linear outside of $DT^*N$ (with slope not a geodesic length), and we let $a<b$ be regular values for the action associated to $H$. With this data we construct a spectrum $Z_a^b(H)$ out of the sequence ($r\to \infty$) of associated Conley indices of the finite dimensional approximations from the previous sections. The homology of this spectrum is not always isomorphic to Floer homology (it is in the oriented \emph{and} spin case), however, the definition of this spectrum is a contractible choice. This is due to the canonical structures we have on the cotangent bundle $T^*N$. This construction will be natural with respect to the inclusion and quotient maps defined when changing the action intervals. When $[a,b]$ contains all critical points we simply write $Z(H)$ for this spectrum. Note that this spectrum depends on the Riemannian structure, because the slope condition at infinity depends on this. However, the transfer map constructed in the next section will not depend on this (up to contractible choices).

In this paper we will use the following rather simple definition of spectra. A spectrum $Z=(Z_n,\sigma_n)_{n\in \N}$ is a sequence of based spaces $Z_n$ (well based - i.e. having cofibrant inclusion of base-point) and cofibrant structure maps
\begin{align*}
  \sigma_n \co \Sigma Z_n \to Z_{n+1},
\end{align*}
Here $\Sigma(-)$ is the reduced suspension of a based space. We will use the notation:
\begin{align} \label{eq:29}
  \sigma_n^m \co \Sigma^{m-n}Z_n \xrightarrow{\Sigma^{m-n-1} \sigma_n} \Sigma^{m-n-1}Z_{n+1} \to \cdots \to \Sigma^1Z_{m-1} \xrightarrow{\sigma_{m-1}} Z_m.
\end{align}
for the composition of the (suspended) structure maps. The homology of a spectrum is defined as
\begin{align} \label{eq:33}
  H_*(Z) = \colim_{n\to \infty} \tH_{*+n}(Z_n)
\end{align}
using the maps
\begin{align*}
  \tH_{*+n}(Z_n) \cong \tH_{*+n+1}(\Sigma Z_n) \xrightarrow{\sigma_{n*}} \tH_{*+n+1}(Z_{n+1}),
\end{align*}
where the first map is the suspension isomorphism (uses well-based). The reader unfamiliar with spectra can consult Appendix~\ref{cha:app} for some properties. This also contains a description of spectra using CW complexes which describes the relation with Morse homology, and hence heuristically explains why the homology of the spectrum $Z_a^b(H)$ is Floer homology (up to orientations).

It will not in the given context be natural to construct a space $Z_a^b(H)_n$ for all $n\in \N$. So, in the following we construct every ${k}^{\textrm{th}}$ space in the spectrum $Z_a^b(H)$ and fill in the gaps afterwards.

Let $r_0\geq 1$ be so large that for any $r\geq r_0$ we have that the finite dimensional approximation $S_r$ and its pseudo-gradient $X_r$ are defined using $H$ and the sub-division $\alpha_j=1/r$. By Lemma~\ref{infiniteflow} (and Lemma~\ref{lem:Homoindex:1}) their exist good index pairs $(A_r,B_r)$ for each $r\geq r_0$. To be able to define spectra we need to compensate for the fact that in Section~\ref{cha:suspmaps} we got a relative Thom construction using the vector bundle of $TN$ and not a standard $d$-fold reduced suspension. Note, that the standard $k$-fold suspension $\Sigma^k (A/B)$ of $A/B$ is canonically homeomorphic to $(A,B)^{\zeta^k/}$, where $\zeta^k$ denotes the trivial metric bundle (we will use this notation over any base).

Let $F, E \to M$ be two metric vector bundles over $M$. Let $(A,B)$ be a pair in $M$. We can iterate the relative Thom space construction as follows.
Let $\pi \co E\to M$ be the projection to the base. Then we may define the pair
\begin{align*}
  ((A,B)^{E-})^{\pi^*(F)-}
\end{align*}
The total space of $\pi^*(F)$ is canonically identified with $E\oplus F$, and we have a canonical homeomorphism
\begin{align*}
  DE \oplus DF \cong D(E\oplus F)
\end{align*}
by scaling each line in $E\oplus F$. This takes $(SE \oplus DF) \cup (DE \oplus SF)$ to $S(E\oplus F)$, and by putting these together we get a canonical identification
\begin{align} \label{eq:1:pairident}
    ((A,B)^{E-})^{\pi^*(F)-} \cong (A,B)^{E\oplus F-}.
\end{align}
This is what we will use to ``untwist'' the tangent bundles $TN$.

Pick an isometric embedding $TN\subset \zeta^k$ of the tangent bundle into the trivial metric bundle, and let $\nu$ denote the normal bundle, and assume it has dimension at least 2. By abuse of notation we define $TN=\ev_0^*TN$ which is a special case in Section~\ref{cha:suspmaps}. We also define
\begin{align*}
  \nu = \ev_0^*\nu
\end{align*}
as a vector bundle over $T^*\Lre N$. We have (by the choices made above) a canonical isomorphism of metric vector bundles:
\begin{align} \label{eq:15}
  \nu \oplus TN \cong \zeta^k \qquad (\textrm{as metric vector bundles over } T^* \Lre N).
\end{align}
We will get back to canonicality of this choice, but for now we consider these choices fixed.

We will almost define the $(r+1)k^{\textrm{th}}$ space in the spectrum $Z_a^b(H)$ for $r$ by
\begin{align} \label{eq:1spect}
  Z_a^b(H)^p_{(r+1)k} = (A_r,B_r)^{\nu^{(r+1)}/}.
\end{align}
Here the ${}^p$ refers to preliminary, and we will change this slightly (not up to homotopy). Here $\nu^{r+1}=\nu^{\oplus (r+1)}$, and $(A_r,B_r)$ is a good index pair for $(S_r,X_r)$.

The map $h_0$ from Equation~\eqref{eq:41} (inducing the homotopy equivalence in Proposition~\ref{prop:Suspension:1} after applying the flow) can be extended to include these normal bundles in the following way. Define a lift $f_0$ fitting into the diagram
\begin{align*}
  \xymatrix{
    \nu^{r+1} \oplus \nu \oplus TN \ar[d] \ar[r]^-{f_0} & \nu^{r+2} \ar[d] \\
    TN \ar[r]^{h_0} & T^* \Lambda_{r+1} N
  }
\end{align*}
by the formula
\begin{align*}
  f_0(\arz,w,v) = (h_0(\arz,v),w).
\end{align*} 
Here $w\in (\nu^{r+2})_\arz$ and $v\in (TN)_{\arz}=T_{q_0}N$  and since $h_0$ commutes with $\ev_0$ - we may naturally consider $w$ as a vector in $(\nu^{r+2})_{h_0(\arz,v)}$. Note that viewing this as a vector bundle map over the bases in the bottom of the diagram this is a linear isometry in each fiber. Hence it is a pull back of metric vector bundles.

It is easy to incorporate the flow of $-X_{r+1}$ into this lift. Indeed, we lift the flow of $-X_{r+1}$ to the bundle $\nu^{r+2}$ by choosing a compatible connection (a contractible choice). So by composing with such a flow we construct a lift $f_t$ of $h_t$ such that the diagram
\begin{align*}
  \xymatrix{
    \nu^{r+1} \oplus \nu \oplus TN \ar[d] \ar[r]^-{f_t} & \nu^{r+2} \ar[d] \\
    TN \ar[r]^{h_t} & T^* \Lambda_{r+1} N
  }
\end{align*}
commutes. Again this is a fiber-wise linear isometry so it is a metric vector bundle pullback diagram. In Proposition~\ref{prop:Suspension:1} we saw that $h_t$ induces a homotopy equivalence
\begin{align*}
  \th_t \co (A_r,B_r)^{TN/} \to A_{r+1}/B_{r+1}.
\end{align*}
Putting the isometry $f_t$ of vector bundles over a map like this induces a new map on the ``untwisted'' indices:
\begin{align*}
  \tf_t \co (A_r,B_r)^{\nu^{r+1}\oplus \nu \oplus TN /} \to (A_{r+1},B_{r+1})^{\nu^{r+2}/}.
\end{align*}

\begin{Lemma}
  \label{lem:Spectrum:1}
  The map $\tf_t$ is a homotopy equivalence (in the oriented case - in the unoriented case we will consider odd $r$ and increasing it 2 at a time).
\end{Lemma}

\begin{proof}
  Start by assuming that $N$ and hence $TN$ and $\nu$ are orientable. Since good index pairs are cofibrant (Lemma~\ref{lem:Homoindex:3}) we get the commuting diagram:
  \begin{align*}
    \xymatrix{
      \tH_*((A_r,B_r)^{TN/}) \ar[r]^{h_{t*}} \ar[d]^{\cong} & \tH_*(A_{r+1}/B_{r+1}) \ar[d]^{\cong} \\
      H_*((A_r,B_r)^{TN-}) \ar[r]^{h_{t*}} \ar[d]^{\cong} & H_*((A_{r+1},B_{r+1})) \ar[d]^{\cong} \\
      H_{*+(r+2)(k-d)}((A_r,B_r)^{\nu^{r+2}\oplus TN-}) \ar[d]^{\cong}\ar[r]^{f_{t*}}& H_{*+(r+2)(k-d)}((A_{r+1},B_{r+1})^{\nu^{r+2}-}) \ar[d]^{\cong}\\
      \tH_{*+(r+2)(k-d)}((A_r,B_r)^{\nu^{r+2}\oplus TN/}) \ar[r]^{f_{t*}}& \tH_{*+(r+2)(k-d)}((A_{r+1},B_{r+1})^{\nu^{r+2}/})
    }
  \end{align*}
  Indeed, this uses (in vertical order); excision, then the Thom-isomorphism for the \emph{orientable} bundle $\nu^{r+2}$ (plus naturality of the Thom-isomorphism), and then excision again.

  Since the top map is an isomorphism by Proposition~\ref{prop:Suspension:1} it follows that the lower horizontal map is an isomorphism, and since Thom-spaces of vector bundles with dimension at least 2 are simply connected the lemma follows. Indeed, the space are simply connected since any representative of an element in $\pi_1$ of the Thom space can be made transversal to the zero section - and hence not intersect the zero section for dimension reason, and then pulled off to the sphere (which is part of the base-point in the quotient) by radial projection homotopy.

  The case where $N$ is not orientable is \emph{not} completely straightforward. Indeed, it is not generally true that a map of pairs $g\co (A,B) \to (A',B')$ which induces a homotopy equivalence on the quotients induces a relative homology equivalence on the pair with any choice of coefficients. This means that one can find an example of this with a bundle $E\to A'$, where the induced map $(A,B)^{g^*E/} \to (A',B')^{E/}$ is not even a homology equivalence.

  So, in the case where $N$ is not orientable it is convenient to jump two $r$'s at a time. Indeed, the proof of Proposition~\ref{prop:Suspension:1} generalizes to proving that the ``composed maps''
  \begin{align*}
    (A_r,B_r)^{TN\oplus TN} \to (A_{r+2},B_{r+2})
  \end{align*}
  is a homotopy equivalence. Note, however, that this requires doing everything in the previous section again, but with the slightly more complicated projection $T^*\Lambda_{r+2} N \to T^*\Lre N$. However, one can still identify a section of this as the fiber-wise unique critical points, and the Hessian of the normal bundle is now two copies of the Hessian in Lemma~\ref{lem:Suspension:3} (both these statements follows by using that lemma twice). This homotopy equivalence now implies using Thom-isomorphism (as above) on the oriented bundle $\nu^{\oplus 2}$ (oriented because of the factor 2) that
  \begin{align*}
    (A_r,B_r)^{TN\oplus TN\oplus\nu \oplus \nu} \to (A_{r+2},B_{r+2})^{\nu\oplus \nu}
  \end{align*}
  is a homology equivalence - hence as above a homotopy equivalence.
\end{proof}

\begin{Remark}
  \label{rem:Spectrum:2}
  We could have chosen to add $\nu^r$ instead of $\nu^{r+1}$. In fact, for any vector bundle $V$ and any $l\in \Z$ we could have added $\nu^{r+l}\oplus V$ and get a lot of different spectra. However, the choice we have is the most natural choice; indeed, it fits with the standard spectrum transfer map $N^{-TN} \to L^{-TL}$, and we even conjecture $\Lambda j_!$ to be a ring-spectrum homomorphism (using twisted Chas-Sullivan products) - and a sketch of a proof of this is contained in \cite{mythesis}. However, the specific alternative of adding $\nu^r$ is what gives rise to the alternative discussed in Remark~\ref{rem:Intro:1}. We will explain this and the relation in a series of remarks and corollaries by considering what happens if we add $\nu^r$ instead of $\nu^{r+1}$. The spectra constructed this way will be decorated with primes. I.e. denoted ${Z'}_a^b(H)$.
\end{Remark}

The source of $\tf_t$ is canonically identified as
\begin{align*}
  (A_r,B_r)^{\nu^{r+1}\oplus \nu \oplus TN /} \cong  (A_r,B_r)^{\nu^{r+1}\oplus \zeta^{k} /} \cong  \Sigma^{k} (A_r,B_r)^{\nu^{r+1}/} = \Sigma^{k} Z_a^b(H)^p_{(r+1)k}.
\end{align*}
The target needs a little adjustment to be similarly recognized. Indeed, recall that $(A_{r+1},B_{r+1})$ was in Proposition~\ref{prop:Suspension:1} an index pair defined for the finite dimensional approximation using the subdivision $\alpha=(0,1/r,\dots,1/r)$ (which satisfies Equation~\eqref{eq:52} for all $r\geq 1$). However, using the obvious convex homotopy from this sub-division to
\begin{align*}
  \alpha=(1/(r+1),\dots,1/(r+1))
\end{align*}
we may use Lemma~\ref{hominv} to (with a contractible choice) identify the target with $Z_a^b(H)_{(r+2)k}$ and $\tf_t$ thus induces a map
\begin{align} \label{eq:17}
  \tau_r \co \Sigma^k Z_a^b(H)^p_{(r+1)k} \to Z_a^b(H)^p_{(r+2)k},
\end{align}
which in the oriented case is a homotopy equivalence by the above lemma (and composing two with the right suspensions added is in the unoriented case a homotopy equivalence). The choices to construct this is contractible (still considering the embedding $N\subset \zeta^k$ of vector bundles fixed and not part of the choices made). We now define the spectrum and fill in the ``gaps''.

\begin{Proposition} \label{prop:spectrum}
  The sequence of maps $\tau_{r}, r\geq r_0$ defines a spectrum $Z_a^b(H)$, and another definition of this (using other choices) is related by a contractible choice of homotopy equivalences. Furthermore, it is naturally compatible with inclusions and quotients of Conley indices.
\end{Proposition}

\begin{proof}
  Let the normal bundle $\nu$ and the isomorphism $\nu \oplus TN \cong \zeta^k$ be fixed as above. The only reason why we won't use the preliminary spaces defined above is that technically it is easier to work with spectra where the structure maps are cofibrations (which we do). So since the maps defined above are not cofibrations, we replace the spaces with the iterated mapping cylinders. That is, we define $Z_a^b(H)_{(r+1)k}$ as the iterated mapping cylinder of the maps
  \begin{align} \label{eq:67}
    \Sigma^{(r-r_0-1)k}\tau_{r_0},\Sigma^{(r-r_0-2)k}\tau_{r_0+1}, \dots, \Sigma^k\tau_{r-2}, \tau_{r-1}
  \end{align}
  See Appendix~\ref{cha:app} for a description of this and a discussion about contractible choices. The gaps in between every $k$th space are filled by making the spectrum ``constant''. Indeed,
  \begin{align*}
    Z_a^b(H)_n = \left\{
      \begin{array}{ll}
        \{*\} & n < (r_0+1)k \\
        \Sigma^{n_1} Z_a^b(H)_{n_2 2k} & n\geq (r_0+1)k
      \end{array}
    \right.
  \end{align*}
  Here $n_1$ is the remainder in $\{0,\dots,k-1\}$ of $n$ when diving by $k$ and $n_2$ is the integral division so that $n= n_1 + n_2k$.
  
  The structure maps $\sigma_n \co \Sigma Z_a^b(H)_n \to Z_a^b(H)_{n+1}$ of the spectrum are defined using the inclusion of the mapping cylinders, which defines the structure maps $\sigma_{(r+1)k}^{(r+2)k}$. The gaps are again filled with constants (identities):
  \begin{align*}
    \sigma_n =  \left\{\begin{array}{ll}
        \id & k \nmid n+1 \\
        \sigma_{n+1-k}^{n+1} & k \mid n+1
      \end{array} \right.  
  \end{align*}
  for $n\geq (r_0+1)k+1$ (since $\Sigma \{*\} \cong \{*\}$ is an initial object in the category of based spaces - the structure maps for $n\leq (r_0+1)k$ are canonically defined).

  We get a canonically homotopy equivalent spectrum if we increase $r_0$ and forget a finite number of Conley indices. Indeed, since we are leaving out the first part of the sequence of the maps in Equation~\eqref{eq:67} on the remaining non-trivial levels, we see that the mapping cylinders gets shorter. However, since the last space is still there the inclusion is a homotopy equivalence. So we have a canonical level-wise cofibrant inclusion of one into the other, which from a certain level (the new $r_0$) is a level-wise homotopy equivalence - hence a homotopy equivalence of spectra.

  So the only part of the construction of $Z_a^b(H)$ that is not at this point a contractible choice is the choice of $k$ and the embedding $N\subset \zeta^k$. However, if we increase $k$ to $k'$ the following two things happen:
  \begin{itemize}
  \item the space of embeddings $TN\subset \zeta^{k'}$ is more connected than the old (connectivity goes to $\infty$),
  \item we add a trivial factor $\zeta^{k'-k}$ to $\nu$ and this corresponds to adding trivial suspensions (a total of $(r+1)(k'-k)$ ) to each of the spaces and consequently mapping cylinders above.
  \end{itemize}
  It follows that modulo the usual reordering of suspension factors in the definition of spectra this is a (weakly) contractible choice. This reordering of suspension factors can be handled by introducing e.g. symmetric spectra (see \cite{SchwedeSymSpec}) or EKMM spectra (see \cite{MR1417719}).
\end{proof}

\begin{Corollary}
  \label{cor:Spectrum:1}
  The spectrum $Z_a^b(H)$ is homotopy equivalent to the shifted suspensions sub-spectrum
  \begin{align*}
    \Sigma^{\infty-(r+1)k}Z_a^b(H)^p_{(r+1)k} \subset Z_a^b(H)
  \end{align*}
  for all $r\geq r_0$.
\end{Corollary}

Note that we could have defined $Z_a^b(H)$ as this shifted suspension spectrum, but then one needs to go trough arguments similar to the content of the above to argue that canonically this does not really depend on $r_0$. Also, our choice makes the spectra easier to handle - since they are more functorially defined; and later, when we will be taking a limit of these we will loose this property anyway.

\begin{proof}
  In the oriented case this follows from Lemma~\ref{lem:Spectrum:1}. Indeed, by definition the inclusion is a homotopy equivalence at level $(r+1)k$ (the mapping cylinder deformation retracts onto this part), and by the lemma all higher structure maps are homotopy equivalences. In the non-oriented case it follows from the fact that the sequence of suspended structure maps
  \begin{align*}
    Z_a^b(H)_{(r+1)k} \xrightarrow{\Sigma^{2k}\tau_r} 
    Z_a^b(H)_{(r+2)k} \xrightarrow{\Sigma^{k}\tau_{r+1}} 
    Z_a^b(H)_{(r+3)k} \xrightarrow{\tau_{r+2}} 
    Z_a^b(H)_{(r+4)k} 
  \end{align*}
  satisfies that the composition of the first two and the last two are both homotopy equivalences - hence they are all homotopy equivalences. So in fact Lemma~\ref{lem:Spectrum:1} is true also for odd $r$.
\end{proof}

\begin{Corollary} \label{cor:hominv}
  Any smooth homotopy of the Hamiltonian $H^s, s\in I$ and smooth homotopy $a_s<b_s$ of regular values for the action associated to $H^s$ induces a (contractible choice) homotopy equivalence
  \begin{align*}
    Z_{a_0}^{b_0}(H^0) \to Z_{a_1}^{b_1}(H^1).
  \end{align*}

  Furthermore, for $c_s$ another regular value such that $a_s<c_s<b_s$ this is compatible with the natural quotients and inclusion from Conley indices induced on the spectra.
\end{Corollary}

\begin{proof}
  This is almost Lemma~\ref{hominv}. However, in that lemma $a,b$ and $c$ were fixed values. However, the lemma is easily generalized to values depending on $s$ by applying a diffeomorphism $\phi_s \co \R \to \R$ depending on $s$ such that $a=\phi_s(a_s),b=\phi_s(b_s)$ and $c=\phi_s(c_s)$ are constants.

  The first part of the corollary now follows from this generalized Lemma~\ref{hominv}. Indeed, we have the diagram
  \begin{align*}
    \xymatrix{
      \Sigma^k Z_{a_0}^{b_0}(H^0)_{(r+1)k} \ar[r]^{\tau_r} \ar[d]^{\simeq} & Z_{a_0}^{b_0}(H^0)_{(r+2)k} \ar[d]^{\simeq} \\
      \Sigma^k Z_{a_1}^{b_1}(H^1)_{(r+1)k} \ar[r]^{\tau_r} & Z_{a_0}^{b_0}(H^1)_{(r+2)k}
    }
  \end{align*}
  which is a homotopy commutative (with a contractible choice of homotopy). So this induces a contractible choice of homotopy equivalence from one to the other.

  The second part follows since the homotopy equivalences from Lemma~\ref{hominv} and the maps $\tau_r$ are compatible with inclusions and quotients.
\end{proof}


\section{The Construction of the Transfer as a Map of Spectra}\label{vitconst}

In this section we construct the map $(\Lambda j)_!$ of spectra in Theorem~\ref{thm:2}. However, we will not yet identify the stable homotopy types of the source and target as the Thom-spectra:
\begin{align*}
  \Lambda N^{-TN}  \qquad  \textrm{and} \qquad \Lambda L^{-TL+\eta}.
\end{align*}
In fact the second identification will not be canonical, and we postpone defining the virtual vector bundle $\eta$ (given a virtual vector bundle as $-TL+\eta$ or $-TN$ we describe in Appendix~\ref{cha:app} how to define these Thom-spectra). The method of construction is similar to that of Viterbos, and we will use a limit of certain Hamiltonians to define a map of spectra
\begin{align*}
  \Lambda j_! \co Z \to W,
\end{align*}
where $Z$ and $W$ later will be proven to be homotopy equivalent to the above. However, it will follow rather directly from the construction that there is a commutative diagram
\begin{align*}
  \xymatrix{
    Z \ar[r]^{\Lambda j_!} & W \\
    N^{-TN} \ar[u] \ar[r]^{j_!} & L^{-TL} \ar[u]
  }
\end{align*}
where $j_!$ is the usual transfer map of manifolds.

By the Darboux-Weinstein Theorem (see e.g. \cite{MR1698616}) we can by choice of Riemannian structures on $L$ and scaling in $T^*N$ assume that
\begin{align*}
  j\colon D_{1/2}T^*L \subset D_{1/2}T^*N
\end{align*}
is a symplectic embedding. To distinguish between coordinates in $T^*N$ and $T^*L$, we denote them by $(q_N,p_N)$ and $(q_L,p_L)$ respectively. So when we write $\norm{p_L}$ we mean the $L$-norm and similar for $\norm{p_N}$, which thus defines two different functions on $D_{1/2}T^*L$. It is very important for the construction that exactness of the embedding implies that $p_N\dd q_N-p_L\dd q_L = \lambda_N - \lambda_L$ defined on $D_{1/2}T^* L$ is exact (we will often omit the $j$ from the notation). This implies that the two action integrals $\int_{\gamma} \lambda_N -Hdt$ and $\int_{\gamma} \lambda_L-Hdt$ are equal on closed curves in $D_{1/2}T^*L$. This means that if we have a Hamiltonian on $T^*N$, which restricted to $D_{1/2}T^*L$ depends only on $\norm{p_L}$, then we can use the method in Remark~\ref{geocalc} to calculate the action integral on closed 1-periodic orbits. In the following this is important to keep in mind.

\begin{figure}[ht]
  \centering
  \begin{tikzpicture}[scale=3]
    \draw[->] (-0.2,0) -- (1.2,0);
    \draw[->] (0,-0.2) -- (0,1.2);
    \draw[dashed] (0,0.275) -- (1,1.075);
    \draw[dotted] (0,0.25) -- (0.66,1.1);
    \draw[dotted] (0,0.8) -- (1.2,1.05);
    \draw (-0.02,1) node [left] {$1$} -- (0.02,1);
    \draw (-0.02,0.8) node [left] {$\ffd$} -- (0.02,0.8);
    \draw (-0.02,0.25) node [left] {$\efd$} -- (0.02,0.25);
    \draw[thick] (-0.2,0) --  (0,0) to[out=0,in=180] (0.48,0.75) --  (0.50,0.75) to[out=0,in=215] (0.75,0.875) to[out=35,in=180] (1,1) --  (1.2,1);
    \draw (1,-0.02) node [below] {$1$} -- (1,0.02);
    \draw (0.5,-0.02) node [below] {$\eh$} -- (0.5,0.02);
    \draw (0.25,-0.02) node [below] {$\efd$} -- (0.25,0.02);
    \draw (0.75,-0.02) node [below] {$\tfd$} -- (0.75,0.02);
    \draw [->] (1.4,1.2) node [right] {slope $\mu_N$}  -- (1.2,1.05);
    \draw [->] (0.7,1.3) node [right] {slope $\mu_L$}  -- (0.66,1.1);
  \end{tikzpicture}
  \caption{The function $f$}
  \label{fig:functionf}
\end{figure}
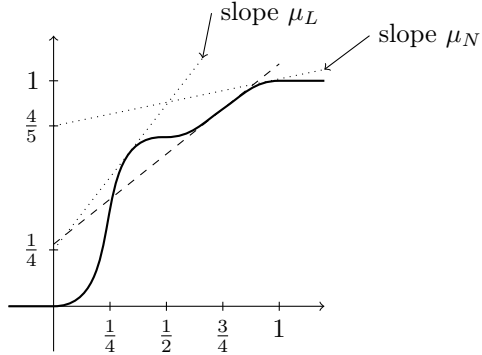

First we choose a function as in figure~\ref{fig:functionf}. That is, construct $f\co \R \to \R$ smooth such that
\begin{itemize}
\item $f(x)=0$ when $x\leq 0$,
\item $f(x)=1$ when $x\geq 1-\epsilon$,
\item $f(x)=\tfd$ in a neighborhood of $x=\eh$,
\item $f$ is convex on the intervals $]0,\efd[$ and $]\eh,\tfd[$,
\item $f$ is concave on the intervals $]\efd,\eh[$ and $]\tfd,1[$, and
\item the inflection point at $x=\tfd$ has tangent intersecting the 2. axis above $\efd$.
\end{itemize}
By construction there are unique tangents to $f$ which intersects the 2. axis at $\efd$ and $\ffd$. These are the dotes lines in the figure, and we denote the slopes of these by $\mu_L$ and $\mu_N$ respectively.

\begin{Remark}
  \label{rem:Viterbo:1}
  We may and will assume that the two tangents to $f$ with slope $\mu_N$ and $\mu_L$ are in fact tangents to all orders. We do this because it will make a technical point (in later sections) much easier.
\end{Remark}

Now we use $f$ to define a smooth family of smooth Hamiltonians $H^s \co T^*N \to \R$ depending smoothly on $s>0$. Indeed, define
\begin{align*}
  H^s(z) = \left\{
    \begin{array}{ll}
      sf(\norm{p_L}) & z \in \im(j) \\
      s\tfd & z \in D_{1/2}T^*N - j(D_{1/2}T^*L) \\
      sf(\norm{p_N}) + h_\infty(\norm{p_N}) & z \notin D_{1/2}T^*N
    \end{array}
  \right.
\end{align*}
Here $h_{\infty} \co \R \to \R$ is such that
\begin{itemize}
\item $h_{\infty}(x)=0$ when $x<1-\epsilon$,
\item $h_{\infty}$ is convex and
\item $h_{\infty}(x) = \mu_\infty t - c_\infty$ for $x\geq 1$.
\end{itemize}
Here the constants $\mu_\infty,c_\infty>0$ are so small that for
\begin{align} \label{eq:18}
  s_1 = 6 c_\infty  
\end{align}
the finite dimensional approximations $S_1$ ($r=1$) in Section~\ref{Florlike} is well-defined for $H^{s_1}$. This is, indeed, possible since we can make $\norm{H^{6c_\infty}}_{C^2}$ small by making both $c_\infty$ and $\mu_\infty$ small. We may also assume that there are no 1-periodic orbits of the flow - hence the critical points of $S_1$ are the same as the critical points of $H^{s_1}$.

Note that we are only adding this small $h_\infty$ so that we get a slightly positive slope at infinity, which is not a geodesic length, for any $s>0$. Heuristically one may ignore this detail, but to be absolutely precise we have added it, and note that the proof of Lemma~\ref{infiniteflow} was made significantly easier by adding this (although a similar yet slightly more general statement is true without adding this slope at infinity).

By the assumptions on $f$ and since the tangents of $t\mapsto s+h_\infty(t)$ intersects the 2. axis in the interval $[s,s-c_\infty]$ we get by the calculation of the action of orbits for $H^s$ using tangents (described in section~\ref{action}) that as long as $s\geq s_1=6c_\infty$ we have
\begin{itemize}
\item All critical values of the action (and hence the finite dimensional approximations) from orbits outside $D_{1-\epsilon}T^*N$ lies in $[-s,-s + c_\infty] \subset [-s,-\tfrac56 s]$.
\item All critical values of the action with action in $]-\infty,-s\ffd]$ comes from orbits outside of $D_{3/4}T^*N$.
\item $-s\ffd$ is a critical value if and only if $s \mu_N$ is not the length of a closed geodesic on $N$.
\item All critical values of the action with action in $[-s\efed,\infty[$ comes from orbits inside of $D_{1/2} T^*L$.
\item $-s\efd$ is a critical value if and only if $s \mu_L$ is not the length of a closed geodesic on $L$.
\end{itemize}

Since the set of lengths of closed geodesics (both for $N$ and $L$) is closed and has measure 0 we see that for almost all $s \geq s_1$ the values $-s\ffd$ and $-s\efed$ are regular for the action. We thus pick a strictly increasing sequence $s_l$ tending to $\infty$ (with $s_1$ as above) such that
\begin{itemize}
\item $-s_l\ffd < -s_l\efed$ are regular values for the action associated to $H^{s_r}$.
\end{itemize}
Now let $a_s^L = -\efed s$, $a_s^N = -\ffd s$ and let $b_s$ denote an upper bound on the critical values of $A_{H^s}$ smoothly depending on $s$. Now define the spectra (depending on $l\in \N$):
\begin{align} \label{Zdef}
  Z(l) &= Z_{a_{s_l}^N}^{b_{s_l}}(H^{s_l}) \qquad \textrm{and} \\
  W(l) &= Z_{a_{s_l}^L}^{b_{s_l}}(H^{s_l}). \label{Wdef}
\end{align}
By construction there is the canonical map of spectra $Z(l) \to W(l)$ given by quotients of Conley indices. We now define a spectrum version of symplectic homology of $T^*N$ by constructing a homotopy colimit of spectra (this again means a mapping cylinder construction - see appendix~\ref{cha:app} for a concrete description)
\begin{align} \label{eq:62}
  Z=\hocolim_{l\to \infty} Z(l).
\end{align}
This will have an essentially canonical map using the quotients above to a spectrum version of symplectic homology of $T^*L$, which we define as a similar limit
\begin{align}\label{eq:64}
  W=\hocolim_{l\to \infty} W(l).
\end{align}
We therefore need to define maps of spectra
\begin{align*}
  \kappa_l \co Z(l) \to Z(l+1)
\end{align*}
compatible with the quotient maps induced by the natural quotient maps on Conley indices.

For this we consider the homotopy $H^s$ for $s\in [s_l,s_{l+1}]$ (we will consider this interval instead of $I$ as to not clutter notation). The concavity of $f$ on the intervals $[\efd,\eh]$ implies that for such an $s$ there is a unique tangent of $s f$ in the interval $[\tfd s,s]$ with slope $s_l\mu_N$. Minus the intersection of this tangent with the 2. axis thus defines a regular value, say $d_s^N$ for the action $A_{H^s}$. Similarly, there is a unique tangent in $[\efd s,\eh s]$ with slope $s_l\mu_L$ whose negative intersection with the 2. axis defines a regular value $d_s^L$. Note that by definition we have
\begin{itemize}
\item $d_{s_l}^N = - s_l \ffd =a_{s_l}^N$ \qquad and \qquad $d^L_{s_l} = - s_l \efd = a_{s_l}^L$.
\end{itemize}
However, since we are moving the tangents (defining these values) up (see Figure~\ref{fig:functionf}) as we increase the multiplication factor $s$ we have, by the concavity on the intervals that
\begin{itemize}
\item $d^N_{s} < -s\ffd = a_{s}^N$ and $d^L_{s} < -s\efd = a_{s}^L$ both for all $s\in ]s_r,s_{r+1}]$.
\end{itemize}
In particular we have
\begin{itemize}
\item $d^N_{s_{l+1}} < a^N_{s_{l+1}}$ and $d^L_{s_{l+1}} < a_{s_{l+1}}^L$.
\end{itemize}
We thus define the map of spectra
\begin{align*}
  \kappa_l \co Z(l) \to Z(l+1)
\end{align*}
by using Corollary~\ref{cor:hominv} on this homotopy $s\in [s_l,s_{l+1}]$ with regular values $d^N_s<b_s$, and then compose with the natural quotient
\begin{align*}
  Z_{d^N_{s_{l+1}}}^{b_{s_{l+1}}}(H^{s_{l+1}}) \to Z_{a^N_{s_{l+1}}}^{b_{s_{l+1}}}(H^{s_{l+1}}) = Z(l+1).
\end{align*}
Since the map from Corollary~\ref{cor:hominv} is compatible with quotients we see that we can construct these (and it is a contractible choice) such that we get commutative diagrams
\begin{align} \label{eq:69}
  \xymatrix{
    Z(l) \ar[r]^{\kappa_l} \ar[d] & Z(l+1) \ar[d] \\
    W(l) \ar[r]^{\kappa'_l} & W(l+1).
  }
\end{align}
Here $\kappa_l'$ is induced by restricting $\kappa_l$ on each Conley index, which means that the diagram commutes on the nose. Making it easy to verify that we get a map on the homotopy colimits
\begin{align} \label{eq:32}
  \Lambda j_! \co Z \to W.
\end{align}

\begin{Proposition} \label{prop:Viterbo:1}
  The map of spectra defined above fits into a commutative diagram
  \begin{align*}
    \xymatrix{
      Z \ar[r]^{\Lambda j_!} & W \\
      N^{-TN} \ar[u] \ar[r]^{(\pi\circ j)_!} & L^{-TL} \ar[u].
    }
  \end{align*}
  Here $(\pi \circ j)_!$ is the usual transfer map for a map of manifolds $(\pi \circ j) \co L \to N$.
\end{Proposition}

\begin{proof}
  By construction above we picked $s_1$ such that $Z(1)$ (the first spectrum in the homotopy colimit in Equation~\eqref{Zdef}) is defined by Conley indices already at the level $r=1$ ($r$ as in Section~\ref{Florlike}). Also, the only other critical values (periodic orbits) are; the constants in $D_{1/2}T^*N-D_{1/2}T^*L$ which has critical value $-\tfd s_1$, and the constants on $L \subset DT^*N$ which has critical value $0$. One way to think of this is that the action approximation
  \begin{align*}
    S_1 \co T^* \Lambda_1 N = T^*N \to \R
  \end{align*}
  is approximately minus the Hamiltonian $-H^{s_1}$. In fact we claim the following: The pair $(DT^*N,UT^*N) \subset T^*N$ is an index pair for $(S_1,X_1)$ containing all critical points, and $(D_{1/4} T^*L,U_{1/4} T^*L) \subset D_{1/2}T^*L \subset T^*N$ is an index pair containing the critical point set $L$. To see this claim we prove that the negative pseudo-gradient $-X_1$ of $S_1$ points out of these sets. Indeed, in the case $r=1$ we have no length conditions so $X_1=\nabla S_1$, and we have
  \begin{align} \label{eq:68}
    S_1(z_0) = \int_{\gamma_0} \lambda_0 + p_0^-\epsilon_{q_0} - H(z_0)
  \end{align}
  Here since $r=1$ the notation implies $(q_1^-,p_1^-)=(q_0^-,p_0^-)$. Now for $Z_0=(q_0,p_0)$ with $\norm{p_0}=1$ the Hamiltonian flow is geodesic flow with speed $\mu_\infty$. This implies $\epsilon_{q_0}=\exp_{q_0}^{-1}(q_0) = - \mu_\infty p_0$ (parallel transported to $q_0^-$) and therefore $\te_{q_0}=-\mu p_0$. The approximation of the gradients in Lemma~\ref{gradient} then shows that
  \begin{align*}
    \norm{\nabla_{p_j} S_1 + \mu_\infty p_0} \leq \tfrac{\mu_\infty}{4} \Rightarrow 
    \inner{\nabla_{p_j} S_1, p_0} > 0.
  \end{align*}
  This shows that the gradient of $S_1$ is inward pointing at the boundary of $DT^*N$ hence the first pair is an index pair (the negative gradient points out). Inside $D_{1/2}T^*L$ we have a slightly different setup. Indeed, the gradient $\nabla H^{s_1}$ is orthogonal to the codimension 1 manifold $U_{1/4}T^*L$ (pointing out of $D_{1/4}T^*L$), and we will prove that for small $s_1$ (which we may assume with out loss of generality) the dominating term in the gradient of $S_1$  Equation~\eqref{eq:68} is $-\nabla H$, which means that for small enough $s_1$ the gradient of $S_1$ will point into the index set - hence the negative gradient points out. To see this we realize that the term $p_0^-\epsilon_{q_0}$ is the integration of $\lambda_0$ over the horizontal geodesic going from $(q_0^-,p_0^-)$ to the fiber over $q_0$. Hence we can write the two first terms in Equation~\eqref{eq:68} as the sum of integrating $\lambda_0$ over 2 curves. In fact we can write this sum as the integration of the closed piece-wise smooth curve given by:
  \begin{itemize}
  \item First part is simply $\gamma_j$ which is a curve from $z_0$ to $z_0^-$ (contributing $\int_{\gamma_0} \lambda_0$),
  \item the second part is the horizontal geodesic from $z_0^-$ to the fiber over $q_0$ (contributing $p_0^-\epsilon_{q_0}$),
  \item and the last part (which contributes 0) is the line (geodesic) in the fiber $T_{q_0}^*N$ from the point the second part arrived at (which is $\tp_0^-$) back to $(q_0,p_0)$.
  \end{itemize}
  This closed curve is a geodesic triangle with side lengths bounded by $\norm{ \nabla H^s}_\infty \leq \norm{f} _\infty s$. So, the enclosed symplectic area is of order this squared (or smaller). Also, moving the point $(q_0,p_0)$ does not violently change these curves (the endpoints are smooth functions in $z_0$ and $s$ - even for $0$ and negative $s$) and hence we conclude that the gradient of the two first term in Equation~\eqref{eq:68} is bounded by some constant times $s_1^2$. This is dominated by $-\nabla H^{s_1}$ which is non-zero on the boundary of the proposed index pair, and scales with $s_1$.

 It follows that the map at level $r=1$ on Conley indices (without the added normal bundles) is given by:
  \begin{align*}
    N^{TN} \to L^{TL}.
  \end{align*}
  This realizes the Pontryagin-Thom collapse map, which realizes the transfer map (see e.g. \cite{MR1942249}). Warning: it is not standard that the bundles showing up here are $TN$ and $TL$. However, since we are adding two copies of the normal bundle and desuspending (shifting by $2k$) we get a map of spectra of the type:
  \begin{align*}
    Z(1) \simeq \Sigma^{-2k}N^{TN+2\nu} = N^{-TN} \to \Sigma^{-2k}L^{TL+2 j^*\nu} \simeq W(1)
  \end{align*}
  Now since $L$ is Lagrangian in $T^*N$ we see why the above looked slightly confusing compared to the standard transfer map $N^{-TN} \to L^{-TL}$. Indeed, as virtual bundles we have:
  \begin{align*}
    TL+2j^*\nu \cong & -TL+2TL+2j^*\nu \cong -TL+TL\otimes \C + 2j^*\nu \cong \\
    \cong & -TL+2j^*TN+j^*\nu \cong -TL+\zeta^{2k}.
  \end{align*}
  Hence $W(1) \simeq L^{-TL}$.
\end{proof}

Note that the bundle isomorphisms used in the later part are canonical so we can in fact identify this part canonically, but for the spectrum $W$ we will run into trouble. In Remark~\ref{rem:Spectrum:2} we discussed an alternate possible definition of the spectra. Indeed, let $Z' \to W'$ be the map of spectra constructed as above, but adding only $\nu^r$ copies of the normal bundle (as opposed to $\nu^{r+1}$).

\begin{Corollary}
  \label{cor:Viterbo:1}
  The alternate transfer map fits into a diagram:
  \begin{align*}
    \xymatrix{
      Z' \ar[r]^{\Lambda j_!} & W' \\
      \Sigma^\infty N_+ \ar[u] \ar[r]^{(\pi\circ j)_!'} & L^{TN-TL} \ar[u].
    }
  \end{align*}
  Here $(\pi \circ j)_!'$ is the usual transfer map for a map of manifolds $(\pi \circ j) \co L \to N$.
\end{Corollary}

\begin{proof}
  Same proof, but the virtual bundle classes turns out different (and precisely like this) because we are adding one less copy of $\nu$.
\end{proof}


\section{Generalized Finite Dimensional Approximations} \label{sec:gener-finite-dimens}

In this section we introduce a ``generalization'' of the finite dimensional approximations considered in Section~\ref{Florlike}, and we prove an ``energy bound'' on these; which will make us able to bound gradient trajectories and prove the localization results we need to be able to identify the homotopy types of $Z$ and $W$ from Section~\ref{vitconst}. The general situation we will consider is; given a compact exact symplectic manifold $M$ (satisfying some topological condition, which we address in Remark~\ref{rem:Genfin:1}), with a compatible Riemannian structure $g$ and a Hamiltonian $H$, we will define functions $S_r$ on finite dimensional approximations of the loop space of $M$. The finite dimensional approximations of the loop space will be denote $\Lrb M$ and is essentially the space of $r$-pieced geodesics with energy less than $\beta$. The reason we put the word ``generalization'' in quotes above is that the functions defined in Section~\ref{Florlike} where defined on bigger approximations of the loop space. More concisely, we have $\Lrb DT^*N \subset T^*\Lre N$. However, other than this restriction the approximations here are more general. Indeed, they will depend on a time-dependent choice of Lagrangian $\SLa$ at every point on $M$, and the case of this being the fiber directions in $T^*N$ we recover the old $S_r$ from Section~\ref{Florlike}. We will introduce the important energy type function $E$, which measures how far a piecewise path is from being a periodic orbit of the Hamiltonian flow, and prove the following important proposition (the constants are various bounds on the structure discussed so far, and will be defined below).
\begin{Proposition} \label{energy}
  There exists a $K>1$ (only dependent on $g$) large enough such that for
  \begin{align} \label{eq:53}
    r>K\pare*{\no{H} + C_\SLa^2(\beta+\nor{H}^2)}
  \end{align}
  we have
  \begin{align*}
    \norm{\nabla E}^2 \leq 20E \leq 40\norm{\nabla S_r}^2 \leq 80E
  \end{align*}
  on $\Lrb M$. Equality holds if and only if $E=0$.
\end{Proposition}
This proposition has a very important implication: the critical points of $S_r$ are the 1-periodic orbits regardless of the Riemannian structure and $\SLa$ (for large $r$). However, we will later see that in the case of a non-degenerate critical point, the Morse index will depend on $\SLa$ (its Maslov index is important).

During this section we will slowly put more and more lower bounds on $K$, but we will make sure that these ``adjustments'' does not depend on $\beta$,$\Sla$ nor $H$ when we do. However, to make the formulation of the lemmas and corollaries in this section more palpable we will not mention any adjustment needed in the formulation of the lemma/corollary, but simply adjust it in the proof.

\begin{Remark}
  \label{rem:Genfin:2}
  The constants $K$, $C_\SLa$, $\nor{H}$ and $\beta$ will all be assumed to be greater than 1. Indeed, this is not going to influence our ability to use the result, and the equation that $r$ should satisfy without this assumption is much more complicated than Equation~\eqref{eq:53}.

  Furthermore, in Section~\ref{Florlike} we defined finite dimensional approximations and index pairs only for $r>C\no{H}$ (for $C$ from Section~\ref{flowline}), and since the goal is to compare this to a more general construction we will assume $K$ to be bigger than this $C$. In fact there is a $C$ coming from $T^*N$ with its induced Riemannian structure, but we will also assume this from the one associated with $T^*L$ using the Riemannian structure we picked in Section~\ref{vitconst}).
\end{Remark}

\subsection{Preliminaries}
Let $(M^{2d},\partial M)$ be a smooth compact manifold with boundary. Let $(M',\lambda)$ be an open exact ($\omega=d(-\lambda)$ is non-degenerate) symplectic manifold without boundary containing $M\subset M'$. So $M$ is an exact symplectic manifold inside a slightly larger $M'$ acting as a ``buffer'' around the boundary of $M$. Let $g$ be a Riemannian structure on $M'$ compatible with $\omega$, and let $J$ be the associated almost complex structure. Notice that any compact exact symplectic $M$ has such a buffer. Associated to this structure we have a constant
\begin{align} \label{eq:54}
  \delta_M > 0
\end{align}
which should be smaller than the injective radius of the exponential function on $M$ (mapping into $M'$), but we will need to make it even smaller later. However, when doing so we make sure that it only depends on $M\subset M'$ and their structures.

We assume we are given a Hamiltonian $H \co M' \to \R$ such that the Hamiltonian flow preserves $M$. As before we denote the (semi) $C^2$-norm of $H$ by $\no{H}$ (Equation~\eqref{eq:12}). We will, however, also have to involved the (semi) $C^1$-norm:
\begin{align*}
  \nor{H} = \max(\max_{z\in M} \norm{\nabla H},1).
\end{align*}
Again we made it bigger than 1 to not make Equation~\eqref{eq:53} more complicated.

\begin{ex}
  \label{exa:Genfin:1}
  Two important examples to keep in mind are $M=DT^*N \subset T^*N = M'$ and $M=D_{1/2}T^*L \subset T^*N=M'$ (as in Section~\ref{vitconst}), with any compatible Riemannian structure. The Hamiltonians that this will be used on are not precisely those from that section, but some related Hamiltonians (and we will relate them later). Indeed, for Proposition~\ref{energy} to be useful we will need to narrow the Hamiltonians such that the action we consider is a small interval proportional to $1/r$. However, we will not consider this until the next section, which uses a family of narrowing Hamiltonians to construct good index pairs and later spectra as in Section~\ref{sec:gener-funct-spectr}.
\end{ex}

As mentioned above $S_r$ will depend on a choice of Lagrangian at each point $z\in M$. To formalize this we introduce the following notions.

\begin{Definition}\label{ladefn}
  For any $d\in \N$ we let $\La(d)$ denote the Grassmannian of Lagrangian sub-spaces in $\R^{2d}=\C^d$. Using the standard inner product on $\R^{2d}$ we may induce a canonical Riemannian structure on the manifold $\La(d)$.

  For any symplectic vector bundle $\xi \to M$ denote by $\La(\xi) \to M$ the fiber bundle with fiber $\La(\xi)_q\cong \La(\dim_{\C}(\xi))$ the Grassmannian of Lagrangian sub-spaces in $\xi_q$ for $q\in M$. If $\xi$ has a fiber-wise compatible inner product and the manifold has a Riemannian structure, we may choose a Riemannian structure on $\La(\xi)$ as follows: each fiber is a Grassmannian of Lagrangian subspaces of a vector space with a compatible inner product, which means it has an induced Riemannian structure. We then choose an arbitrary orthogonal complement to the fiber, and use the Riemannian structure on $M$ to define the inner product on this complement. We denote these the horizontal directions in $\La(\xi)$ and we may choose these smoothly. Note that there might be a canonical choice of horizontal directions when $\xi=TM$, but in the following that will not matter, and we will simply fix any such choice.
\end{Definition}

We now consider as part of the given data (needed to define a generalized approximation of the action) a \emph{time-dependent} smooth section
\begin{align*}
  \SLa \co M\times S^1 \to \La(TM).
\end{align*}
The reason that we need this to be time-dependent (the $S^1$ factor) will not be clear until Section~\ref{stabloc}, and may seem weird since we did not consider time dependent Hamiltonians (although we easily could).
For each $t\in I$ we use $S^1=I/\{0,1\}$ and consider
\begin{align*}
  \SLa_t = \SLa(-,t) \co M \to \La(TM).
\end{align*}

\begin{Remark}
  \label{rem:Genfin:1}
  Not all Liouville domains has such a section. Indeed, if any of the odd Chern classes are non-torsion this cannot exist. Indeed, the map $B\Or \to BU$ given by $\otimes \C$ (which is the structure we need to lift to define a single $\SLa_t$) has torsion odd Chern classes. However, in cotangent bundles such a structure always exists - in fact canonically so.
\end{Remark}

Since we, at a technical point later, will be working with a non-compact family of such sections we will need to assume a specific bound. Indeed, assume that
\begin{align*}
  \SLa(z,-) \co S^1 \to \La(TM)
\end{align*}
has energy bounded by some fixed constant $C_{\SLa}>1$. That is,
\begin{align}\label{eq:31}
  e(\SLa(z,-) = \int_0^1 \norm{\pdd{t}\SLa(z,-)}^2dt \leq C_{\SLa}
\end{align}
for all $z\in M$, and assume that for fixed $t\in S^1$ we have the first and second derivatives bounded by
\begin{align} \label{eq:30}
  \norm{D\SLa_t} \leq C_\SLa \qquad \norm{D^2\SLa_t} \leq C_\SLa.
\end{align}
Of course, for a single $\SLa$ the existence of such a constant follows by compactness and smoothness, but to make the results in this section work (for the non-compact family we will consider later) we will in the following use these concrete bounds on $\SLa$, and we will in fact \emph{not} be assuming that this is smooth in $t$ only in $z$. Note that these bounds imply
\begin{align} \label{eq:49}
  \dist(\SLa_t(z),\SLa_{t'}(z')) \leq C_\SLa\dist(z,z')+\sqrt{C_\SLa\absv{t-t'}}
\end{align}
so continuity is automatic from the bounds.

\begin{ex} \label{SLdef}
  With $M=DT^*N$ as in the above examples we define the time independent $\SLa^N$ as the canonical section in
  \begin{align*}
    \La(T(DT^*N)) \to DT^*N,
  \end{align*}
  given by the vertical directions (the $p$-directions). Indeed, this is a canonical Lagrangian in the tangent space at each point in $(q,p)\in T^*N$. We may restrict this to $D_{1/2}T^*L \subset DT^*N$, but there we also have the section $\SLa^L$ by using vertical directions in $T^*L$. It is in fact the difference in these two choices we are going to explain explicitly.
\end{ex}

Our approximations will again depend on a subdivision
\begin{align*}
  \alpha=(\alpha_0,\dots,\alpha_{r-1})  \qquad \textrm{with} \qquad \sum_j \alpha_j = 1.
\end{align*}
Precisely as in Section~\ref{Florlike}, and again we assume for simplicity Equation~\eqref{eq:52}.

The finite dimensional approximation of loops we will use for the finite dimensional approximations of the action is
\begin{align*}
  \Lrb M = \{\arz \in (\inte M)^r \mid e(\arz) < \beta \}.
\end{align*}
where $\arz=(z_j)_{j\in \Z/r}$ and $e$ is the energy given by
\begin{align} \label{eq:25}
  e(\arz) = r\sum_j \dist(z_j,z_{j+1})^2.
\end{align}
Here $\inte M$ denotes the interior of $M$. Note that this is in fact the usual energy of loops if one interprets $\arz$ as a piece-wise geodesic (each parameterized by an interval of length $1/r$). Notice that each $z_j$ and $z_{j+1}$ will be closer than $\delta_M$ (from Equation~\eqref{eq:54}) if we assume
\begin{align*}
  \sqrt{\beta/r} \leq \delta_M.
\end{align*}
In this case $\Lrb M$ is a well-defined open and finite dimensional manifold. By Equation~\eqref{eq:53} we can assume this if $K>\delta_M^{-2}$. Note that we are \emph{not} using the sub-division $\alpha$ in the definition of $\Lrb M$. We could have done this, and in some ways this might have been more natural, but the formulas turn out easier this way.

\begin{ex} \label{ex:newnew}
  With $M$ and $M'$ as in Example~\ref{exm:4} we can if we also assume that $K>\delta_0^{-2}$ (from Equation~\eqref{eq:11}) see that
  \begin{align*}
    \Lrb M = \Lrb DT^*N \subset T^*\Lre N
  \end{align*}  
  is an open submanifold. Here the latter was defined in Section~\ref{Florlike} using $\delta_0$. Moreover, if we assume that $K>9\delta_0^{-2}$ we get that this inclusion is inside the set where the pseudo-gradient $X_r$ in that section where defined to be equal to the gradient (defined right before Equation~\eqref{Xulig}).
\end{ex}

\begin{figure}[ht]
  \centering
  \includegraphics{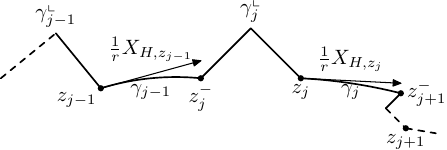}
  \caption{Curves involved in definition of finite dimensional approximation.}
  \label{wedgedcurve}
\end{figure}

To define the finite dimensional approximations depending on $\SLa$ we will also need the following geometric construction. Given a Lagrangian subspace $L\subset T_z M$ we define for any close by point $z^- \in M$ with $\dist(z^-,z)$ small enough the \emph{L-curve}
\begin{align*}
  \gamma^\llcorner(z^-,z,L) \co I \to M
\end{align*}
as the continuous path from $z^-$ to $z$ defined by:
\begin{itemize}
\item parameterized by constant arc length on $[0,\eh]$ we go from $z^-$ along a geodesics to the closest point in $D=\exp_z(D_{2\delta_M} L)$ (the path meets $D$ orthogonally at a point closer than $2\delta_M$ to $z$),
\item parameterized by constant arc length on $[\eh,1]$ we follow a geodesic from that closest point to $z$ (this is inside $\exp_z(D_{2\delta_M}L)$).
\end{itemize}
Notice that we may assume (by possibly making $\delta_M$ smaller) that this is well-defined for all $\dist(z^-,z) \leq \delta_M$ and all $L$. We will use these to close up piecewise flow curves. The actual parameterizations of the L-curves are unimportant since we will use them only for integrating 1-forms. Notice that although $z$ and $z^-$ are both in $M\subset M'$ we allow this L-curve to exit and reenter $M$. Figure~\ref{wedgedcurve} illustrates many aspects of how we will use these L-curves to define $S_r$.

Observe that if $r> 4 \nor{H} \delta_M^{-1}$ and we define for each $\arz \in \Lrb M$
\begin{align} 
  \gamma_j(t)&=\varphi_t(z_j),\quad t\in [0,\alpha_j], \label{flowpice}
\end{align}
then each $\gamma_j$ is shorter than $\nor{H} \alpha_j < 2\nor{H}/r < \delta_M/2$ (follows from H2a) above). So, we adjust $K$ to satisfy
\begin{align*}
  K>\max(4\delta_M^{-2},4\delta_M^{-1})
\end{align*}
so that this is true (uses Equation~\eqref{eq:53} and $\nor{H}>1$), and so that $\sqrt{\beta/r} < \delta_M/2$. This implies that the distance from
\begin{align}
  z_j^-&=\gamma_{j-1}(\alpha_j) \label{zjminus}
\end{align}
to $z_j$ is less than $\delta_M$ and we may thus define
\begin{align}
  \gamma^\llcorner_j &=
  \gamma^\llcorner(z_j^-,z_j,\SLa_{j/r}(z_j)). \label{wedgecurvedef}
\end{align}
With these we may finally define the finite dimensional approximation of the action as
\begin{align}  \label{genSrdef}
  S_r(\arz) = S_{(r,g,\SLa,H)}(\arz) = \sum_{j\in \Z/r} \pare*{
    \int_{\gamma_j}(\lambda-H dt) +
    \int_{\gamma^\llcorner_{j+1}} \lambda }.
\end{align}
Note how the curves all fit together (as pictured in Figure~\ref{wedgedcurve}) to integrate $\lambda$ over a closed curved. This is a very important point; indeed, the gradient of $S_r$ now only depends on $\omega$ and not $\lambda$, which is an important point for the usual action $A_H$.

To analyze the gradient of $S_r$ we also define
\begin{align}
  \epsilon_j &= -\exp_{z_j}^{-1}(z_j^-) \label{ej}
\end{align}
This is basically the vector pointing from $z_j^-$ to $z_j$, but moved to the tangent space at $z_j$.
\begin{ex}
  \label{exrem:1}
  As in the above examples where $M=DT^*N$, $\SLa=\SLa^N$ is the vertical directions and the metric $g$ is induced from a metric on $N$, we may compare this to the previous definition of $S_r$ in Equation~\eqref{eq:7}. Indeed, the curve $\gamma^\llcorner_j$ will because $\exp(\SLa(q,p))=T_q^*N$ be the curve first going in horizontal direction from the fiber over $q_j^-$ to the fiber over $q_j$ (this is the geodesic $\epsilon_{q_j}$ lifted horizontally to start at $(q^-_j,p_j^-)$, which ends at $(q_j,\tp_j^-)$ by definition of $\tp_j^-$) then it goes in the fiber from $\tp_j^-$ to $p_j$. Integrating this over $\lambda$ we precisely get the term $p_j^-\epsilon_{q_j}$ (the movement in the fiber direction does not contribute to the integral). So Equation~\eqref{genSrdef} generalizes the definition from Equation~\eqref{eq:7} - albeit only on the subset $\Lrb DT^* N \subset T^*\Lre N$.

  In this example $\epsilon_j \approx (\epsilon_{q_j},\epsilon_{p_j})$ and these two components are basically the tangents to the two pieces of the L-curve.
\end{ex}

\begin{Remark}\label{embedrem}
  The function $S_r$ can be approximated by $S_r\approx A_H\circ i_r$ where $i_r$ is an embedding
  \begin{align*}
    i_r \colon \Lrb M \to \Lambda M.
  \end{align*}
  We may define $i_r$ as the curve depicted in figure~\ref{wedgedcurve} with parameterization on the flow curves $\gamma_j$ ``almost'' as defined, but leaving a little parameterization room for the L-curves to be parameterized by a \emph{very} short interval. Because this is short we get almost no contribution from the integration of the $H dt$ term over the L-curve part and we approximately get the expression for $S_r$. This can be made more rigorous such that the Conley indices defined by such embeddings is the same as the one defined by $S_r$. We will not need this, but it is a good justification for the name finite dimensional approximation, and describes the relation with Floer homology discussed in the introduction.
\end{Remark}  

We now introduce the important energy type functional appearing in Proposition~\ref{energy}:
\begin{align} \label{eq:2EE}
  E(\arz)=\sum_j \norm{\epsilon_j}^2.
\end{align}
We have defined $E$ such that it is zero if and only if the curves $\gamma_j$ fit together to a 1-periodic orbit for the Hamiltonian flow of $H$. Indeed, $E$ can be thought of as a finite version of the energy relative to the Hamiltonian flow
\begin{align*} 
   \int_{s^1} \norm{\gamma'(t)-X_{H}(\gamma(t))}^2dt.
\end{align*}
However when considering this and comparing with $e$ in Equation~\eqref{eq:25} we note that a factor $r$ has been omitted in the expression for $E$. This is evident in the follow lemma which tells us that $rE$ is in a sense equivalent to $e$, i.e. bounding one bounds the other.

\begin{Lemma} \label{lem:Ebound}
  We have
  \begin{align*}
    rE(\arz) \leq 2e(\arz) + 8\nor{H}^2 \qquad \textrm{and}
    \qquad e(\arz) \leq 2rE(\arz) + 8\nor{H}^2.
  \end{align*}
\end{Lemma}

\begin{proof}
  The length of the Hamiltonian flow is bounded by $\nor{H}$, and each
  small flow curve $\gamma_j$ is bounded in length by $\alpha_j \nor{H} \leq 2\nor{H}/r$. This
  implies
  \begin{align*}
    E(\arz) = &\sum_j \dist(z_j^-,z_j)^2 \leq \sum_j (\dist(z_j,z_{j+1})+2\nor{H}/r)^2 \leq \\
    \leq & \sum_j 2(\dist(z_j,z_{j+1})^2 + 4\nor{H}^2/r^2) \leq 2e(\arz)/r+8\nor{H}^2/r.
  \end{align*}
  Similarly $\dist(z_j,z_{j+1}) \leq \dist(z_j^-,z_j) + 2\nor{H}/r$ proves the other inequality.
\end{proof}

\subsection{Approximations in local coordinates}

We will essentially have to reduce the proof of Proposition~\ref{energy} to the flat case in local coordinates. However, to get all the bounds we need it is convenient to make sure that we can for any $z\in M$ find ``good'' coordinates with certain bounds. The following lemma takes care of this.

\begin{Lemma}
  \label{lem:Loclin:2}
  By making $\delta_M$ smaller we can assume that: for any $z\in M$ and any Lagrangian subspace $L \subset T_zM$ there exists a symplectic chart $h \co D_{\epsilon}^{2d}(0) \to M'$ with $h(0)=z$ and $h^*g$ equal to the standard structure at $0$ and $h^*(L)=i\R^d$ also at $0$. Furthermore, we may assume that $B_{\delta_M}(z) \subset \im(h)$ and that there are bounds independent of $z$ and $L$ (using the Riemannian structure on $M$ and the standard on $D_{\epsilon}^{2d}(0)$) on the first and second derivatives of $h$ and $h^{-1}$.
\end{Lemma}

\begin{proof}
  Cover $M$ by finitely many open Darboux charts $h_i \co D_{\epsilon'}(0) \subset \R^{2d} \to M'$ which extends smoothly to the boundary so that we have bounds on all derivatives. Then pick a smooth isometric and symplectic (Hermitian) trivialization
  \begin{align*}
    \phi_i \colon h_i^* TM' \cong U_i \times C^n
  \end{align*}
  of the tangent bundle $TM'$ pulled back to each of these charts (and their closures - so as to have global bounds on
  derivatives). Using this we can for each $z\in h_i(U_i)$ define a new chart $h_i^z \colon U^z_i \to M'$ by
  \begin{align*}
    h^z_i(w) = h_i((\phi_{i\mid z})^{-1}(w)+h_i^{-1}(z)).
  \end{align*}
  This sends $0$ to $z$ and the pull back of $g$ is the standard Riemannian structure at $0$ and this choice depends smoothly on $z$ for fixed $i$. It follows that for small enough $\delta_M$ small balls around $0$ of these cover $M$ (in the way the lemma specifies - ignoring the $V$) and we get the global bound with property that $h^*g$ is standard at $0$.

  We may make sure that the pull back of $L$ at $0$ is $i\R^d$ by multiplying the entire chart with an element in $U(d)$. This does not change bounds on the derivatives.
\end{proof}

When working in local charts in this subsection we will be assuming the local chart comes from Lemma~\ref{lem:Loclin:2}. That is, we assume that $g$ is a Riemannian structure on $D_{\epsilon}^{2d}(0)$ standard at $0$, and that the first and second derivatives of $g$ is bounded. We will need to understand what happens if we locally vary the Lagrangians $\SLa$ so we will not generally assume compatibility with $\SLa$ (as in the lemma). Because of the bounds we have on the charts and their first and second derivatives the constants in this subsection can be chosen as global constants working on all charts in $M$ from Lemma~\ref{lem:Loclin:2}.

Define for small $\epsilon>0$ the functions (with compact domain)
\begin{align*}
  F,F_g \co D_\epsilon^{2n}\times D_\epsilon^{2n} \times \La(d) \to \R
\end{align*}
by
\begin{align*}
  F_g(z^-,z,L) = \int_{\gamma^{\llcorner}(z^-,z,L)} \lambda_0
\end{align*}
using $g$ to define the L-curve, and define $F$ by the same formula but we use the standard Riemannian structure to define the L-curve. Notice that this ``standard'' L-curve has its two geodesic parts parallel to $\R^n$ and $i\R^n$ respectively.
 
\begin{Lemma}\label{wedgelem}
  There is a constant $c>0$ (depending on the bounds on $g$) such that 
  \begin{align*}
    \norm{F(z^-,z,L)-F_g(z^-,z,L)} \leq c\norm{(z^-,z)}^3.
  \end{align*}
\end{Lemma}

\begin{proof}
  Since $\La(d)$ is compact we fix an $L$ and the lemma is equivalent to showing that
  \begin{align*}
    v^{-2}(F(vz^-,vz,L)-F_g(vz^-,vz,L)) \to 0
  \end{align*}
  for $v\to 0 \in \R$, which is what we will show.
  
  Note that
  \begin{align} \label{eq:21}
    v^{-2}F_g(vz^-,vz,L)= F_{g^v}(z^-,z,L),
  \end{align}
  where $g^v$ is the Riemannian structure given at the point $z$ by taking $g$ at the point $vz$. Indeed, this is because:
  \begin{itemize}
  \item Scaling a geodesic for $g^v$ by $v^{-1}$ gives a geodesic for $g$ so the L-curves for $g^v$ scale with $v^{-1}$ to the L-curves of $g$.
  \item Integrating $v^{-2}\lambda_0$ over a path $\gamma$ in $\R^{2d}$ gives the same as integrating $\lambda_0$ over the path $v^{-1}\gamma$.
  \end{itemize}
  Equation~\eqref{eq:21} implies for $F$ that
  \begin{align*}
    v^{-2}F(vz^-,vz,L) =F(z^-,z,L).
  \end{align*}
  So $F(-,-,L)$ is a quadratic form. Since $g^v\to g_0$ ($C^{\infty}$ on compact sets since $g$ at $0$ is standard) for $v\to 0$ we see that this is in fact the limit of Equation~\eqref{eq:21} for $v\to 0$, and is therefore the Hessian of $F_g(\cdot,\cdot,L)$ at $(0,0)$ (independent of $g$), and hence the lemma follows.
\end{proof}

We will need the following corollary of this.

\begin{Corollary}
  \label{cor:Loclin:1}
  There is a constant $C>0$ (depending on the bounds on $g$) such that
  \begin{align*}
    \norm{\nabla (F-F_g)} \leq C (\norm{z^-,z}^2) = C(\norm{z^-}^2+\norm{z}^2).
  \end{align*}
\end{Corollary}

Notice here that $\nabla$ is taken with respect to $z^-,z$ and $L$.

\begin{proof}
  By the above lemma we have a constant $c>0$ such that $\absv{F-F_g} \leq c\norm{z^-,z}^3$. That is: the function is bounded by a constant times the distance cubed to the compact submanifold $\{(0,0)\}\times \La(d)$. Hence the gradient (with respect to all directions) will be bounded by a constant times the distance squared.
\end{proof}

\begin{Corollary}
  \label{cor:Loclin:4}
  There is a constant $C'>0$ (depending on the bounds on $g$) such that
  \begin{align*}
    \norm{\nabla_L F_g} \leq C'(\norm{z^-}^2+\norm{z}^2).
  \end{align*}
\end{Corollary}

\begin{proof}
  The above corollary shows that this is true for $F_g$ if it is true for $F$. It is true for $F$ since $\norm{F(z^-,z,L)-F(z^-,z,L')} \leq C''\dist(L,L')\norm{z^-,z}^2$ for some $C''$ as illustrated in Figure~\ref{fig:area}. 
  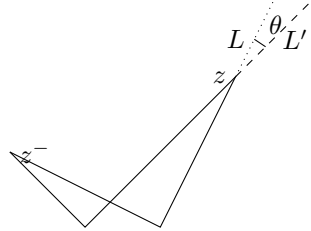
\begin{figure}[ht]
    \centering
    \begin{tikzpicture}
      \draw (0,0) -- (2,-1) -- (3,1);
      \draw [dotted] (3,1) -- node [left] {$L$} (3.5,2);
      \draw (0,0) node [right] {$z^-$} -- (1,-1) -- (3,1) node [left] {$z$};
      \draw [dashed] (3,1) -- node [right] {$L'$} (4,2);
      \draw (3.25,1.5) -- node [above right] {$\theta$} (3.4,1.4);
    \end{tikzpicture}
    \caption{Difference of the symplectic area for $L$ and $L'$ with fixed endpoints. The area of each triangle is bounded by $2\theta\norm{z^--z}^2$ for small $\theta$.}
    \label{fig:area}
  \end{figure}
\end{proof}

We will also need the actual gradient of $F$ for $L=i\R^d$ (the flat case with standard Lagrangian).

\begin{Lemma}
  \label{lem:Loclin:1}
  The gradient of $F^s=F(-,-,i\R^d)$ is given by
  \begin{align*}
    \nabla F^s = (-y^-,x-x^-,y^-,0)
  \end{align*}
  at the point $(z^-,z)=(x^-,y^-,x,y)$.
\end{Lemma}

\begin{proof}
  By the very explicit way L-curves look in the flat structure and the definition of $\lambda_0$ we see directly that
  \begin{align*}
    F(z^-,z) = \inner{y^-,(x-x^-)}.
  \end{align*}
  Note that this is the linear version of the term $p^-_j\epsilon_{q_j}$.
\end{proof}

The following corollary specializing this is in fact all we will need for this gradient.

\begin{Corollary}
  \label{cor:Genfin:1}
  At $z=0$ the gradient of $F^s$ with respect to $z$ is given by
  \begin{align*}
    \nabla_z F^s = (y^-,0).
  \end{align*}
  At $z^-=0$ the gradient of $F^s$ with respect to $z^-$ is given by
  \begin{align*}
    \nabla_{z^-} F^s = (0,x).
  \end{align*}
\end{Corollary}

This corollary is what inspires the next subsection. Indeed, in the flat case in charts around $z_j=0$ this corollary says that the gradient w.r. to $z_j$ of the L-curve integration part of $S_r$ is equal to minus the imaginary part of $\epsilon_j=z_j-z^-_j$. Similarly it says that in charts around $z^-_j=0$ the gradient w.r. to $z^-_j$ (if this could move freely) of the L-curve integration part of $S_r$ is the real part of $\epsilon_j=z_j-z^-_j$.

\subsection{Local approximations of the energy}

Because of the above observation it is convenient to approximate $E$ by some slightly different functions $E'$ (depending on some local choices). So let $\arz\in \Lrb M$ be given. Pick charts as in Lemma~\ref{lem:Loclin:2} $h_j \co D_{\epsilon}^{2d}(0) \to M$ around $z_j$ (pulling back $\SLa_{j/r}(z_j)$ to $i\R^d$) and $h_j^- \co D_{\epsilon}^{2d}(0) \to M$ around $z_j^-$ (pulling back $\SLa_{j/r}(z_j^-)$ to $i\R^d$) for each $j\in \Z/r$. Now define
\begin{align*}
  \epsilon_{x_j}^{h_j^-} & = \Real((h^-_j)^{-1}(z_j)) \in \SLa_{j/r}(z_j^-)^{\perp} \subset T_{z_j^-}M \\
  \epsilon_{y_j}^{h_j} & = -i\Imag((h_j)^{-1}(z_j^-) \in \SLa_{j/r}(z_j) \subset T_{z_j}M.
\end{align*}
Here we consider $\Real \co \C^d \to \R^d \subset \C^d$ and $i\Imag \co \C^d \to i\R^d \subset \C^d$ as orthogonal real projections to real part and imaginary part, and since the chart identified $\SLa_{j/r}$ (at different points) with $i\R^d$ we can interpret these (as indicated) as tangent vectors inside the Lagrangians. These are important because Corollary~\ref{cor:Genfin:1} tells us that in the flat case the real parts and imaginary parts of $\epsilon_j$ are important for the gradient of $S_r$. These are approximate orthogonal projections to $\SLa_{j/r}(z_j)^\perp$ and $\SLa_{j/r}(z_j^-)$ of $\epsilon_j$ and we thus have the following heuristical description of these:
\begin{itemize}
\item The vector $\epsilon^{h_j^-}_{x_j}$ approximates the tangent to the first part of the L-curve from $z_j^-$ to $z_j$ (with length the length of this geodesic).
\item The vector $\epsilon^{h_j}_{y_j}$ approximates the tangent to the second part of the L-curve from $z_j^-$ to $z_j$ (with length the length of this geodesic).
\end{itemize}
In fact we will later see that $\epsilon_j \approx \epsilon_{x_j}^{h_j^-} + \epsilon_{y_j}^{h_j}$ almost as an orthogonal decomposition - so these are to be thought of as linear versions of the L-curve. So, in Examples~\ref{exrem:1} these approximates the components of $\epsilon_j$ given by $\epsilon_{q_j}$ and $\epsilon_{p_j}$. However, we change the notation to $x$ and $y$ because in the general case it may not be compatible with the cotangent bundle structure, and even when it is the usual structure it is not clear that these are exactly equal to $\epsilon_{q_j}$ and $\epsilon_{p_j}$ - only approximately. We will make several of these statements more explicit in the following proof. However, the only result we will explicitly need for these involves comparing $E$ to the function
\begin{align} \label{eq:9}
  E'(\arz) = \sum_{j} \norm{\epsilon_{x_j}^{h_j^-}}^2 + \norm{\epsilon_{y_j}^{h_j}}^2.
\end{align}

\begin{Lemma}
  \label{erew}
  For $r$ as in Equation~\eqref{eq:53} we have
  \begin{align*}
    \absv{E-E'} \leq \tfrac{1}{100} E
  \end{align*}
  independent on the choice of charts $h_j$ and $h_j^-$ (as long as they satisfy the derivative bounds that we assume by Lemma~\ref{lem:Loclin:2}).
\end{Lemma}

\begin{proof}
  Consider any point $z^- \in M$ and a fixed $t\in S^1$. Define $\epsilon_z=-\exp^{-1}_{z}(z^-)$ (similar to $\epsilon_j$) for \emph{any} close by $z\in M$. Assume $h^-\co D_{\epsilon}^{2d}(0) \to M$ is local coordinates around $h^-(0)=z^-$ from Lemma~\ref{lem:Loclin:2} pulling back $\SLa_t$ to $i\R^d$. Now define
  \begin{align*}
    \epsilon_x = \Real((h^-)^{-1}(z))
  \end{align*}
  (similar to above but suppressing the dependency on the charts). As above this approximates the tangent of the first part of the L-curve $\gamma^\llcorner(z^-,z,\SLa_t(z))$.
  
  To make this statement explicit consider everything in the local coordinates $h^-$ (which we now suppress from the notation). Consider, as in the proof of Lemma~\ref{wedgelem}, the ``zoom'' in the sense that we may change these two structures depending on $v\in I$ as follows:
  \begin{align*}
    (g^v)_w = g_{vw} \qquad \textrm{and} \qquad (\SLa^v_t)(w) = \SLa_t(vw).
  \end{align*}
  So $g^1=g$, $\SLa^1_t=\SLa_t$, but $g^0$ and $\SLa^0_t = i\R^d$ are the standard structures. Now as $v$ tends to zero all structures converge uniformly on any compact set. So it follows that we get the following limit behavior:
  \begin{align*}
    \frac{\pi_{\SLa_t(z)^\perp}(\epsilon_z) - \epsilon_x}{\norm{z}} = \pi_{\SLa_t(z)^\perp}(\tfrac{\epsilon_z}{\norm{z}}) - \tfrac{\epsilon_x}{\norm{z}} \to 0
  \end{align*}
  for $z\to 0$. Indeed, for small $z$ this happens in a very small ball, and for the standard structures this formula is equal to 0. Here $\pi_{\SLa_t(z)^\perp}$ is the orthogonal projection using $g$ of the tangent vectors of $T_zM$ onto the subspace $\SLa(z)$. It follows by smoothness of the numerator that we get a bound
  \begin{align*}
    \norm{\pi_{\SLa_t(z)^\perp}(\epsilon_z) - \epsilon_x} \leq C_M C_\SLa\norm{z}^2.
  \end{align*}
  Here $C_M>1$ is a bound only depending on $M\subset M'$ and their structures (not the section $\Gamma$). We can assume this specific bound since the constant here only depends on the second order behavior of the numerator at $z=0$, and this only depends on the second order behavior of $g$ and $\SLa_t$. Hence by Equation~\eqref{eq:30} and the bounds assumed by Lemma~\ref{lem:Loclin:2} we get such a bound. In fact, we can replace $\norm{z}$ with $\dist(z^-,z)$ by again making $C_M$ larger. That is, we have a global bound (for close by $z^-$ and $z$):
  \begin{align*}
    \absv*{\norm{\pi_{\SLa_t(z)^\perp}(\epsilon_z)}-\norm{\epsilon_x}} \leq C_M C_\SLa\dist(z^-,z)^2.
  \end{align*}
  For this arbitrary $t$, and for any chart $h^-$ around $z^-$ satisfying the bounds from Lemma~\ref{lem:Loclin:2}

  Similarly we get using a chart $h$ at any $h(0)=z$ (also from the lemma) by considering varying $z^-$ the bound
  \begin{align*}
    \absv*{\norm{\pi_{\SLa_t(z)}(\epsilon_z)}-\norm{\epsilon_y} } \leq C_M C_\SLa\dist(z^-,z)^2.
  \end{align*}
  with $\epsilon_y=-\im(h^{-1}(z^-))$.

  Now use these bounds for each pair $z_j^-,z_j$, and the charts fixed before the lemma, together with
  \begin{align*}
    \absv*{\norm{a}^2-\norm{b}^2}=\absv*{\inner{a-b,a+b}} \leq \norm*{a-b}(\norm*{a}+\norm*{b}),
  \end{align*}
  $\dist(z_j^-,z_j)=\norm{\epsilon_j}$, and $\norm{\epsilon_j}^2=\norm{\pi_{\SLa_{j/r}(z_j)}(\epsilon_j)}^2 + \norm{\pi_{\SLa_{j/r}(z_j)^\perp}(\epsilon_j)}^2$ to conclude:
  \begin{align*}
    \absv*{E(\arz)-E'(\arz)} &= \sum_j \pare*{\norm{\epsilon_j}^2-\norm{\epsilon_{x_j}}^2-\norm{\epsilon_{y_j}}^2} \leq \\
     &\leq \sum_j (2C_MC_\SLa \norm{\epsilon_j}^2) 3\norm{\epsilon_j} \leq 6C_MC_\SLa\sqrt{\frac{2\beta+8\nor{H}^2}{r}}E.
  \end{align*}
  Here the $3\norm{\epsilon_j}$ comes from the factor ($\norm{a}+\norm{b}$) and:
  \begin{itemize}
  \item the fact that the orthogonal projections have length less than $\norm{\epsilon_j}$ and
  \item $\norm{\epsilon_{x_j}} \leq 2\norm{\epsilon_j}$ and $\norm{\epsilon_{y_j}} \leq 2\norm{\epsilon_j}$, which are easy consequences of the bounds above (for $\norm{\epsilon_j} < \tfrac{1}{C_MC_\SLa}$, which we can get by adjusting $K$).
  \end{itemize}
  We now see that if we pick $K>50(100)^2C_M^2$ we get (when $r$ satisfies Equation~\eqref{eq:53}) that
  \begin{align*}
    (E-E')^2 \leq 6C_M^2C_\SLa^2 \frac{2\beta+8\nor{H}^2}{r} E^2 \leq 50C_M^2 \frac{C_\SLa^2(\beta+\nor{H}^2)}{r} E^2 \leq \frac{1}{100^2} E^2.
  \end{align*}
\end{proof}

At the end of the proof we saw the following, which will be useful again later.

\begin{Corollary}
  \label{cor:Loclin:3}
  For $r$ as in Equation~\eqref{eq:53} we have:
  \begin{align*}
    \norm{\epsilon_{x_j}^{h_j^-}} \leq 2\norm{\epsilon_j} \qquad \textrm{and} \qquad
    \norm{\epsilon_{y_j}^{h_j}} \leq 2\norm{\epsilon_j}
  \end{align*}
\end{Corollary}

\subsection{Gradient approximations using extensions}
We will analyze several gradients in the following, and in more than one case it is convenient to consider the same trick as we employed in the proof of Lemma~\ref{gradient1half}. There we extended our function to a larger manifold where the $z_j^-$ coordinates did not depend on $z_j$. That is, define the manifold
\begin{align} \label{eq:22}
  W \subset \Lrb M \times (\inte M)^r
\end{align}
by $(\arz,\arzm) \in W$ if $\dist(z_j,z_j^-) < \delta_m$. In particular, if $r > 4\nor{H}^{-1}\delta_M$ (as assumed to define $S_r$) we have the embedding
\begin{align} \label{eq:23}
  \chi = \chi_{s,r} \co \Lrb M \to W
\end{align}
given by
\begin{align}
  \label{eq:564}
  \chi(\arz) =(\arz,\varphi_{\alpha_0}(z_0),\dots,\varphi_{\alpha_{r-1}}(z_{r-1})).   
\end{align}
Now, let $f\co \Lrb M \to \R$ be a smooth function. If
\begin{align} \label{eq:24}
  f^e \co W \to \R
\end{align}
is \emph{any} extension of $f$ in the sense that $f^e\circ \chi =f$ then we may calculate the gradient of $f$ by the formula
\begin{align*}
  \nabla f = (D\chi)^\dagger(\nabla f^e),
\end{align*}
which coordinate wise may be written as
\begin{align} \label{gradcalc}
  \nabla_{z_j} f = \nabla_{z_j} f^e +
  (D_{z_j}(\varphi_{\alpha_j}))^\dagger(\nabla_{z_j^-} f^e).
\end{align}
The following lemma is an easy consequence, and it is a proof of the first part of Proposition~\ref{energy} above.

\begin{Lemma} \label{lem:Enablabound}
  For $r$ as in Equation~\eqref{eq:53} we have that
  \begin{align*}
    \norm{\nabla E}^2 \leq 20 E.
  \end{align*}
\end{Lemma}

\begin{proof}
  As discussed above we extend the definition of $E$ to a function
  \begin{align*}
    E^e \co W \to \R
  \end{align*}
  by the simple formula
  \begin{align*}
    E^e(\arz,\arzm) = \sum_j \dist(z_j,z_j^-)^2.
  \end{align*}
  So that $E^e(\chi(\arz))=E(\arz)$. The gradient of $E^e$ is easily calculated (see e.g. \cite{MR0163331}) to be
  \begin{align*}
    \nabla_{(z_j,z_j^-)} E^e= (-2\exp_{z_j}^{-1}(z_j^-),-2\exp_{z_j^-}^{-1}(z_j))
  \end{align*}
  Note that both components have the norm $2\norm{\epsilon_j}$.

  The assumption in Equation~\eqref{eq:52} implies
  \begin{align*}
    \no{\alpha_jH} \leq 2\no{H}/r.
  \end{align*}
  Let $F$ be the identification of nearby tangent vectors induced by some chart at $z_j$ from Lemma~\ref{lem:Loclin:2}. Since $F_{z,z}=\id$ we see by Lemma~\ref{Hamflow} that there is a $\delta$ such that if $2\no{H}/r<\delta$ then we have the bound 
  \begin{align*}
    \norm{(D_{z_j}(\varphi_{\alpha_j}))^\dagger} \leq \sfd
  \end{align*}
  on the operator norm. With $K>2 \delta^{-1}$ we have $2\no{H}/r < \delta$ and it thus follows from Equation~\eqref{gradcalc} that
  \begin{align*}
    \norm{\nabla E}^2 \leq & \sum_j \norm{2\exp^{-1}_{z_j}(z_j^-) - 2(D_{z_j}(\varphi_{\alpha_j}))^\dagger(\epsilon_{j+1})}^2 \leq  \\ 
    \leq & 8 \sum_j \pare{\norm{\epsilon_j}^2 + \norm{(D_{z_j}(\varphi_{\alpha_j}))^\dagger(\epsilon_{j+1})}^2}\leq 8(E + (\sfd)^2 E) \leq 20E.
  \end{align*}
\end{proof}

We will also need to extend $S_r$ by $\chi$ to relate its gradient to $E$. However, before doing this we will get rid of the annoying fact that $\SLa_t$ varies with the points in $M$. That is, we will reduce the problem to local charts with $\SLa_t$ constant in those charts.

Since, we will make choices for each point in $\Lrb M$ at which we consider the gradient of $S_r$ we now fix such a point $\arw\in \Lrb M$, and consider in the following only points $\arz$ close to this $\arw$, and the goal is to prove Proposition~\ref{energy} at the point $\arw$.

To get rid of the varying behavior of $\SLa$ we now define a new choice of Lagrangians locally at $\arw$. Indeed, pick some charts as in Lemma~\ref{lem:Loclin:2} around each $w_j\in M$, say $h_j \co D_{\epsilon}(0) \to M'$. These induce canonical identifications:
\begin{align} \label{eq:45}
  \La(h_j) \co D_\epsilon(0) \times \La(d) = \La(TD_{\epsilon}(0)) \to \La(TM')_{\mid \im(h_j)}.
\end{align}
This defines a possibly different Riemannian structure (locally), but because of the bounds of the derivatives of $h_j$ these will be equivalent. I.e. there is a constant $C_M$, such that bounding the distance between two points in one structure bounds the distance by this same amount times this constant in the other structure. Now we define an alternate function $G$ to $S_r$ by
\begin{align} \label{eq:28}
  G(\arz) = \sum_{j\in \Z/r} \pare*{ \int_{\gamma_j}(\lambda-H^s dt) + \int_{\gamma^\llcorner(z_j^-,z_j,\La_j)} \lambda },
\end{align}
where $\La_j=\SLa_{j/r}(w_j)$ is now chosen to be constant in the chart $\La(h_j)$. So, this is the exact same function as $S_r$ except that we have replaced the dependence of $\SLa$ with the constants $\La_j$. By definition we have
\begin{align*}
  G(\arw)=S_r(\arw)
\end{align*}
However, much more importantly we have the following bound on their gradient difference at $\arw$.

\begin{Lemma}
  \label{lem:reducetolocal}
  For $r$ as in Equation~\eqref{eq:53} we have
  \begin{align*}
    \pare*{\norm{\nabla(S_r-G)}(\arw)}^2 \leq \tfrac{1}{100} E(\arw)
  \end{align*}
\end{Lemma}

\begin{proof}
  Let $\arz=\arz(u)$ be a path parameterized by unit arc length on $[-\epsilon,\epsilon]$ through $\arw=\arz(0)$. This implies that $\arz\,'(0)$ is a unit vector. We may prove the lemma by proving that $\absv{\pdd{u}(G-S_r)(\arz)}_{u=0}^2 \leq \tfrac{1}{100}E(\arw)$ for all such. Indeed, if $\arz\,'(0)$ is parallel to the gradient we get
  \begin{align*}
    \pare*{\norm{\nabla(S_r-G)}(\arw)}^2 = \absv*{D_\arw(G-S_r)(\arz'(0))}^2 =  \absv*{\pdd{u}(G-S_r)(\arz)}^2_{u=0}.
  \end{align*}
  Notice that $\norm{z_j'(0)} \leq \sqrt{\sum_{i\in \Z/r} \norm{z_i'(0)}^2} = 1$ for all $j\in \Z/r$, and so $\dist(z_j,w_j)\leq u$.

  Using the bounds on the derivative of $\SLa_t$ from Equation~\eqref{eq:30} we get (in $\La(TM')$ distances) that
  \begin{align*}
    \dist(\La_j,\SLa_{j/r}(z_j)) \leq C_\SLa u,
  \end{align*}
  which implies the bound
  \begin{align} \label{eq:43}
    \dist(\La_j,\SLa_{j/r}(z_j)) \leq C_M C_\SLa u
  \end{align}
  as Lagrangians in $\La(d)$ with the standard structure (here the factor $C_M$ is there to convert length bounds in the pull back structure to length bounds in the standard structure).

  The difference between $G$ and $S_r$ at $\arz$ close to $\arw$ is given by:
  \begin{align*}
    S_r(\arz)-G(\arz)= 
    \sum_{j\in \Z/r} \pare*{ \int_{\gamma^\llcorner(z_j^-,z_j,\SLa_{j/r}(z_j))} \lambda - \int_{\gamma^\llcorner(z_j^-,z_j,\La_j)} \lambda }.
  \end{align*}
  We now consider this in the local coordinates for each $j$, and at $u=0$ (where the Lagrangians are equal). We get by the chain rule that the terms in the differential (with respect to $u$) coming from the first two coordinates ($z_j$ and $z_j^-$ - i.e. moving the endpoints of the $\gamma^\llcorner$ paths) cancel and the only part left is the contribution from changing the Lagrangians in the last coordinate, and only the first integral actually depends on $z_j$ in this coordinate. That is; we have
  \begin{align} \label{eq:44}
    \pdd{u}\pare[\Big]{(S_r-G)(\arz)}_{u=0} = 
    \pdd{u}\pare*{\sum_{j\in \Z/r} \int_{\gamma^\llcorner(w_j^-,w_j,\SLa_{j/r}(z_j))} \lambda }_{u=0},
  \end{align}
  where $\SLa_{j/r}(z_j))$ lies in $\La(d)$ depending on $u$, who's derivative is bounded in Equation~\eqref{eq:43}. The derivative of the integral is bounded in Corollary~\ref{cor:Loclin:4}, and combined we get
  \begin{align} \label{eq:20}
    \absv*{\pdd{u}(G-S_r)(\arz)_{u=0}} \leq \sum_{j\in \Z/r} C'C_M C_\SLa \dist(w_j,w_j^-)^2 \leq C'C_MC_\SLa E(\arw).
  \end{align}
  Now squaring this and using Lemma~\ref{lem:Ebound} we get
  \begin{align*}
    \absv*{\pdd{u}(G-S_r)(\arz)_{u=0}}^2 \leq (C'C_MC_\SLa)^2 E(\arw)^2 \leq (C'C_MC_\SLa)^2\frac{2e(\arw)+8\nor{H}^2}{r} E(\arw) \leq \\
    \leq (C'C_MC_\SLa)^2\frac{ 2\beta+8\nor{H}^2}{r} E(\arw),
  \end{align*}
  which proves the lemma by adjusting $K$ (Similarly done as in the end of the proof of Lemma~\ref{erew}).
\end{proof}

To relate the gradient of $S_r$ at $\arw$ to $E$ we will also extend $G$ to $W$ (as we did for $E$ above). So we extend it by the formula
\begin{align} \label{eq:27}
  G^e(\arz,\arzm) = \sum_j\pare*{\int_{\gamma_j} (\lambda-H dt) + \int_{\gamma^\llcorner(z_j^-,z_j,\La_j)} \lambda }.
\end{align}
Here $\gamma_j$ is some smooth choice of paths depending on the end points $z_j$ and $z_{j+1}^-$ extending the Hamiltonian flow curves of $\alpha_jH$. As calculated in Equation~\eqref{piecedif} the choices of $\gamma_j$ does not matter for the gradient of $G^e$ on the image of $\chi$.

\begin{proof}[Proof of Proposition~\ref{energy}]
  As above we consider a fixed point $\arw \in \Lrb M$ at which we want to prove the proposition, and again we pick for each $j$ a chart $h_j$ as in Lemma~\ref{lem:Loclin:2} with $h_j(0)=w_j$ (pulling back $\SLa_{j/r}(z_j)$ to $i\R^d$), and we define $G$ as above. However, we now also pick charts $h_j^-$ using the same lemma around each $w_j^-$ (pulling back $\SLa_{j/r}(z_j^-)$ to $i\R^d$).

  Comparing with Lemma~\ref{lem:Enablabound} and Lemma~\ref{lem:reducetolocal} we see that it is sufficient to prove
  \begin{align*}
    3E \leq 4\norm{\nabla G(\arw)}^2 \leq 5E.
  \end{align*}
  Indeed, assuming this we get:
  \begin{align*}
    \norm{\nabla S_r (\arw)}^2 \leq (\norm{\nabla G(\arw)}+\norm{\nabla (S_r-G)(\arw)})^2 \leq (\sqrt{\tfrac{5}{4}}+\tfrac{1}{100})E \leq 2E
  \end{align*}
  and
  \begin{align*}
    \sqrt{E} \leq & \sqrt{\tfrac{4}{3}} \norm{\nabla G(\arw)} \leq \sqrt{\tfrac{4}{3}}(\norm{\nabla S_r(\arw)}+\norm{\nabla (G-S_r)(\arw)}) \leq \\
    \leq & \sqrt{\tfrac{4}{3}}(\norm{\nabla S_r(\arw)}+\sqrt{\tfrac{1}{100}E}) \quad \Rightarrow \quad E \leq 2\norm{\nabla S_r (\arw)}^2.
  \end{align*}
  So in the following we only consider $G$ and its extension $G^e$ close to the points $\arw$ and $\chi(\arw)=(\arw,\arwm)$ respectively.

  The gradient of $G^e$ with respect to $z_j$ only depends on the position of the points $z_j^-, z_j$ and $z_{j+1}$ and the gradient with respect to $z_{j+1}^-$ only depends on $z_j,z_{j+1}^-$ and $z_{j+1}$.

  We first consider the gradient of $G^e$ with respect to $z_j$ at $\arw$. For this we fix the remaining coordinates $z_i=w_i, i\neq j$ and $\arz^-=\arw^-$. We will calculate the gradient using the chart $h_j$, and we partly suppress this chart from the notation, and consider $g$ as the Riemannian metric induced by $h_j$ close to $0$ and $g^0$ as the standard structure. The symplectic structures agree, and we consider $z_j$ as points close to $0=w_j$ in the domain of $h_j$. We will successively replace $G^e$ by approximating functions $G_1^e$ and then $G_2^e$ defined for such $z_j$ close to $0$, and relate their gradients.

  Firstly, define:
  \begin{align*}
    G_1^e(z_j) =  {
      \int_{\gamma_{z_j,w_{j+1}^-}}(\lambda_0-H dt) } + 
      \int_{\gamma^\llcorner(w_j^-,z_j,i\R^d)} \lambda_0,
  \end{align*}
  These are the terms in $G^e$ which actually depends on $z_j$ and we have replaced the integration over $\lambda=h_j^*\lambda$ with integration over $\lambda_0$. Since the two paths concatenate to a path from $w_j^-$ to $w_{j+1}^-$ independent of $z_j$ the gradient of $G_1^e$ with respect to $z_j$ equals that of $G^e$. That is
  \begin{align}\label{eq:46} 
    \nabla_{z_j} G^e = \nabla_{z_j} G_1^e
  \end{align}

  Then define:
  \begin{align*}
    G_2^e(z_j) =  \pare*{
      \int_{\gamma_{z_j,w_{j+1}^-}}(\lambda_0-H dt) } + 
      \int_{\gamma_0^\llcorner(w_j^-,z_j,i\R^d)} \lambda_0,
  \end{align*}
  Here $\gamma^\llcorner_0$ means we use the Riemannian structure $g^0$ instead of $g$ to define the L-curve. This means that in $\C^d$ this L-curve consists of two straight lines the first parallel to $\R^d$ the other to $i\R^d$. Now since we want the gradient at $z_j=w_j=0$ in the chart Corollary~\ref{cor:Loclin:1} implies that
  \begin{align}  \label{eq:26}
    \norm{\nabla_{z_j} (G_1^e-G_2^e)}_{\mid z_j=w_j} \leq C_M\dist(w_j,w_j^-)^2,
  \end{align}
  with $C_M$ only depending on $M$. Indeed, the difference between these two functions are precisely that we use the two different metrics in the chart to define $\gamma^\llcorner$, and this was precisely the difference between the functions $F$ and $F^g$ in that subsection. Notice in particular that the gradient of the functions in question at $z_j=w_j$ do not depend on the metric since the metrics agree at this point.

  Inspecting the definition of $G_2^e$ we see that the gradient with respect to $z_j$ of the first integration term is 0. Indeed it is a flow curve ending on the zero section of $\C^d = T^*\R^d$, and we saw that this has gradient zero in Equation~\eqref{piecedif} with respect to varying this endpoint. The gradient of the second term was computed in Corollary~\ref{cor:Genfin:1} (since  $z_j=w_j=0$ in the chart), and we get that the gradient of $G_2^e$ is given by
  \begin{align*}
    \nabla_{z_j} G_2^e = J_0 (\epsilon_{y_j}^{h_j}) \qquad \brac*{= -J_0 (0,y_j^-) = (y_j^-,0) \textrm{ in coordinates $z^-=(x^-,y^-)$}},
  \end{align*}
  where $J_0$ is the standard complex structure on $\C^d$, and $\epsilon^{h_j}_{y_j}$ is simply the imaginary part of $z_j^-$.  Combining this with Equation~\eqref{eq:46}, Equation~\eqref{eq:26}, and the fact that $J_0=J$ at $0$ in the chart we get
  \begin{align*}
    \norm{\nabla_{z_j} G^e - J\epsilon_{y_j}^{h_j}}_{z_j=w_j} \leq C_M\dist(w_j,w_j^-)^2.
  \end{align*}

  Similarly, in coordinates $h_{z_j^-}$ we can define functions $G_1^e$ and $G_2^e$ depending on $z_j^-$ and get that
  \begin{align*}
    \norm{\nabla_{z_j^-} G^e - J\epsilon_{x_j}^{h_j^-}}_{z_j^-=w_j^-} \leq C_M\dist(w_j,w_j^-)^2.
  \end{align*}
  Notice in particular that even though there is a slight asymmetry in the definition of $G$ with respect to $z_j$ and $z_j^-$ - we don't see this when using Corollary~\ref{cor:Loclin:1}. Indeed, in that corollary the Lagrangian was also situated at $z_j$, but the resulting bounds were symmetric in $z_j$ and $z_j^-$. For $G$ Equation~\eqref{gradcalc} is
  \begin{align*}
    \nabla_{z_j} G = \nabla_{z_j} G^e + (D_{z_j}(\varphi_{\alpha_j}))^\dagger(\nabla_{z_{j+1}^-} G^e),
  \end{align*}
  so combining the above with this we get
  \begin{align} \label{eq:55}
    \norm{\nabla_{z_j} G - J\epsilon_{y_j}^{h_j} - &D_{z_j}(\varphi_{\alpha_j})^\dagger(J\epsilon_{x_{j+1}}^{h_{j+1}^-})} \leq \notag \\
    \leq & C_M\dist(w_j,w_j^-)^2+\norm{D_{z_j}(\varphi_{\alpha_j})^\dagger}C_M\dist(w_{j+1},w_{j+1}^-)^2 \leq \\
    \leq & C_M(\norm{\epsilon_j}^2+2\norm{\epsilon_{j+1}}^2). \notag
  \end{align}
  The latter is because we can adjust $K$ to make $\norm{D_{z_j}(\varphi_{\alpha_j})^\dagger}$ as close to 1 as we would like (Lemma~\ref{Hamflow} as usual). As indicated we now evaluate everything at $\arz=\arw$, so $\norm{\epsilon_j}=\dist(w_j,w_j^-)$.

  For brevity denote $v=J\epsilon_{y_j}^{h_j}$ and $w=D_{z_j}(\varphi_{\alpha_j})^\dagger (J\epsilon_{x_{j+1}}^{h_{j+1}^-})$. So the above formula bounds $\nabla_{z_j}G-v-w$. The important point now is that $v$ and $w$ are close to being orthogonal. Indeed, one is an $\epsilon_y$ and the other an $\epsilon_x$ - we will make this precise. Abbreviate $w'=J\epsilon_{x_{j+1}}^{h_{j+1}^-}$ and $\Phi:=D_{z_j}(\varphi_{\alpha_j})^\dagger$ such that $w=\Phi(w')$. In the local chart $h_j$ we can assume by adjusting $K$ and using Lemma~\ref{Hamflow} that $\Phi$ is $\epsilon'>0$ close to the identity. For as small an $\epsilon'>0$ as we would want. Similarly, let $V=\SLa_{j/r}(z_j)=i\R^d$ and $W'=\SLa_{(j+1)/r}(z_{j+1}^-)^{\perp}$ (as linear subspace in $\R^{2d}$ using the chart $h_j$) we can by Equation~\eqref{eq:49} see that
  \begin{align*}
    \dist(V,W'^{\perp}) \leq C_\SLa 2\nor{H}/r+\sqrt{C_\SLa}\sqrt{1/r}.
  \end{align*}
  So again by adjusting $K$ we can assume that this is as small as we would like. We end up with the abstract situation:
  \begin{itemize}
  \item We have two vectors $v,w'\in \R^{2k}$, which are in two linear subspaces $v \in V$ and $w\in W'$, which are almost orthogonal.
  \item We then apply a linear map $\Phi$ very close to the identity, which maps $w'$ to $w$ and $W'$ to some $W$.
  \item The point is that we still have $w\in W$ and $W$ is still almost orthogonal to $V$. We conclude that
    \begin{align*}
      \absv{\norm{v+w}^2 - \norm{w}^2 -\norm{v}^2} \leq \epsilon \norm{v}\norm{w},
    \end{align*}
    where $\epsilon>0$ is twice cosine to the angle between $V$ and $W$.
  \end{itemize}
  Using (in order at each step) first that $J$ preserves norm; then the $\epsilon'$ bound above on $I-\Phi$ together with
  \begin{align} \label{eq:56}
    \absv*{\norm{a}^2-\norm{b}^2}=\absv*{\inner{a-b,a+b}} \leq \norm*{a-b}(\norm*{a-b}+2\norm*{b});
  \end{align}
  then the above $\epsilon$ bound on the vectors $v$ and $w$ combined with $\Phi(w')=w$ and $\epsilon'\leq 1$; then Corollary~\ref{cor:Loclin:3}; and finally Equation~\eqref{eq:55} combined with Equation~\eqref{eq:56}, Corollary~\ref{cor:Loclin:3}, and the bound $\norm{\Phi} \leq 2$; we get {\small
    \begin{align*}
      & \absv*{\norm{\nabla_{z_j} G}^2 - \norm{\epsilon_{y_j}^{h_j}}^2- \norm{\epsilon_{x_{j+1}}^{h_{j+1}^-}}^2} = \absv*{\norm{\nabla_{z_j} G}^2 - \norm{v}^2- \norm{w'}^2} \leq \\
      \leq & \absv*{\norm{\nabla_{z_j} G}^2 - \norm{v}^2- \norm{\Phi(w')}^2} + \epsilon'\norm{w'}(\epsilon'\norm{w'}+2\norm{w'}) \leq \\
      \leq & \absv*{\norm{\nabla_{z_j} G}^2 - \norm{v+w}^2} + \epsilon\norm{v}\norm{w'}  + 3\epsilon'\norm{w'}^2 \leq  \\
      \leq & \absv*{\norm{\nabla_{z_j} G}^2 - \norm{v+w}^2} + 4\epsilon\norm{\epsilon_j}\norm{\epsilon_{j+1}} + 12\epsilon'\norm{\epsilon_{j+1}}^2 \leq \\
      \leq & C_M(\norm{\epsilon_j}^2+2\norm{\epsilon_{j+1}}^2)\pare[\Big]{C_M(\norm{\epsilon_j}^2+\norm{\epsilon_{j+1}}^2)+2\norm{\epsilon_j}+2\cdot 2\norm{\epsilon_{j+1}}} + \\
      & + 4\epsilon\norm{\epsilon_j}\norm{\epsilon_{j+1}}+12\epsilon'\norm{\epsilon_{j+1}}^2 \leq \tfrac{1}{100}(\norm{\epsilon_j}^2 + \norm{\epsilon_{j+1}}^2).
    \end{align*}}
  The very last for appropriately small $\epsilon$ and $\epsilon'$, which as argued above can be assumed for appropriate $K$
    
  Using this and Lemma~\ref{erew} we get 
  \begin{align*}
    \absv*{\norm{\nabla G}^2 - E} \leq & \absv*{\norm{\nabla G}^2 - E'} + \absv*{E-E'} \leq \\
    \leq & \tfrac{1}{100}\sum_j(\norm{\epsilon_j}^2 + \norm{\epsilon_{j+1}}^2) + \tfrac{1}{100} E \leq \tfrac{3}{100}E,
  \end{align*}
  which implies that
  \begin{align}\label{eq:2}
    3E \leq 4\norm{\nabla G}^2 \leq 5E
  \end{align}
  at $\arw$. 
\end{proof}


\section{Localization}\label{loclinsec}

In this section we prove a localization result, which will come in handy in the following sections. The localization result can be heuristically formulated as follows. For a family of Hamiltonians $H^\param,\param>0$; with certain bounds on derivatives (depending on $\param$) and a fixed Hamiltonian flow behavior at the boundary of $M$; the generalized finite dimensional approximations from Section~\ref{sec:gener-finite-dimens} has good index pairs in $\Lrb M$ for large $r>0$, small $\param$ and small intervals of action (this is Proposition~\ref{mainstabloc} below). This is especially helpful when $M=D_{1/2}T^*L \subset T^*N$ where we want to relate the spectra $W$ defined in $T^*N$ associated with $D_{1/2}T^*L$ to those defined inside $T^*L$. We start by describing the general setup for the family $H^\param$.

In this section we consider the same setup as the previous section except we have a family of Hamiltonians $H^\param$ for $\param\in ]0,1]$ satisfying the following assumptions.
\begin{itemize}
\item[{\textbf H1)}] {There exist a neighborhood $U \subset M'$ of $\partial M$ such that the Hamiltonian flow $(\varphi_t)^\param$ of $H^\param$ is
    \begin{itemize}
    \item[H1a)] independent of $\param$ on $U$,
    \item[H1b)] preserves the compact closure $\ovl{U}$, and
    \item[H1c)] has no periodic orbits (time 1) on $\ovl{U}$.
    \end{itemize}}
\item[{\textbf H2)}] {There is a constant $C_H>1$ such that for all $\param$ we have:
    \begin{itemize}
    \item[H2a)] {$\nor{H^\param}  \leq C_H$,}
    \item[H2b)] {$\no{H^\param} \leq \param^{-1}C_H$ and}
    \item[H2c)] {We restrict to the action interval $[a_\param,b_\param]$ (both smooth in $\param$ and both for regular for the action) such that $b_\param-a_\param \leq \param C_H$.}
    \end{itemize}}
\end{itemize}
H2c) is the narrowing we have alluded to in previous sections, and to accommodate the possibility (in interesting cases) of this we really need that H2b) does not simply bound $\no{H^\param}$ by $C_H$ (see Example~\ref{exm:4} below).

In this setup we have that the action and its approximation depends on $\param$. We also have that the energy type function $E$ considered in Section~\ref{sec:gener-finite-dimens} depends on $\param$. However, we will suppress some of these dependencies from the notation. The goal in this section is to prove the following proposition.
\begin{Proposition} \label{mainstabloc}
  With $K>1$ as in Proposition~\ref{energy} and $\beta$ large enough there is an $\param_0>0$ small enough such that: for any $0<\param<\param_0$ and 
  \begin{align*}
    r\in[2K C_H \param^{-1} , 3KC_H \param^{-1}]    
  \end{align*}
  a good index pair for the total index of $S_r=S_r^\param\colon\Lrb M \to \R$ (with its gradient) exists.
\end{Proposition}
Note, however, that \emph{even though} the critical points of the finite dimensional approximations are all the same - the dependence on $\SLa$ is so profound that changing it can change the Morse indices of non-degenerate critical points (a change in Maslov index changes the Conley-Zehnder index). It is the topic of several of the following sections to precisely describe how changing $\SLa$ in the definition of $S_r$ changes the \emph{stable} homotopy type of the Conley indices defined by this proposition.

The set of assumptions on $H^\param$ may seem somewhat technical and restrictive but the primary example to keep in mind is the following. However, we will see others (yet similar) in the following sections.
\begin{ex} \label{exm:4}
  As in Example~\ref{exa:Genfin:1} let $M=DT^*N$ and $M'=T^*N$ then let $h\co \R \to \R$ be smooth, convex, with $h(t)=\mu t +c$ for $t\geq 1-\epsilon$ and such that
  \begin{align*}
    H^\param(q,p) = s(h(\param^{-1}\norm{p})-c) + c
  \end{align*}
  is smooth for $\param >0$. Then this is linear at infinity with slope not depending on $\param$, but more importantly all the critical values of the associated action lies in the narrowing interval $[s(h(0)-c)+c,c]$, and it satisfies all of the above assumptions. As $\param$ tends to $0$ this narrows the bend close to the zero section, which makes the function look more and more like the non-smooth function $\mu\norm{p}+c$.
\end{ex}

\begin{Remark}
  \label{rem:Loclin:1}
  An important abstract idea used in the construction of the index pairs in this section is as follows. Assume that $(f,X)$ is a function and pseudo-gradient, and we would like to bound some function $F$ on a possible Conley index. Then in some cases we will be able to construct a cut-off function that satisfy the property of Lemma~\ref{cutoff} but which also bounds $F$. This can be done if we have a bound of the type:
  \begin{align*}
    \absv{X \cdot \nabla F} \leq c(X \cdot \nabla f).
  \end{align*}
  for a some $c>0$ - with equality only at critical points. Indeed, this will make the function
  \begin{align*}
    F - c f + c'
  \end{align*}
  satisfy the conditions in the lemma, and bound $F$ on the associated Conley index constructed in the lemma by
  \begin{align*}
    F \leq cf-c' \leq cb-c'
  \end{align*}
  where $[a,b]$ is the interval we wish to find Conley indices on.
\end{Remark}

The idea is to use this on the function $E$ (defined in Section~\ref{sec:gener-finite-dimens}) to bound the energy of the loops in the Conley index. However, we will need more cut-off functions designed to keep the $z_j$'s away from the boundary of $M$. So define
\begin{align*}
  Q(z) = - \dist(z,\partial M), \qquad z\in M.
\end{align*}
This is not smooth on all of $M$. However in Lemma~\ref{cutoff} we only need smoothness near the boundary. Because $M$ is compact we can find $\tau>0$ such that $Q$ is smooth on the collar
\begin{align} \label{eq:48}
  K_\tau=\{z\in M\mid Q(z)\geq-\tau\}.
\end{align}
By possibly making $\tau$ smaller we may also assume that $K_{\tau} \subset U$, where $U$ is the set on which $H^\param$ has no 1-periodic orbits and is independent of $\param$. We then define the functions $Q_j(\arz)=Q(z_j)$. The next lemma leads to the corollary, which is important to be able to bound $Q_j < -\tau/2$ on the Conley index that we will construct.

\begin{Lemma}
  \label{lem:Loclin:3}
  There exists a $c>0$ such that for any $r$ and $\param$ (as long as $E$ is defined) we have
  \begin{align*}
    E(\arz) > \frac{c}{r}
  \end{align*}
  when $z_j\in K_{\tau}$ for some $j\in \Z/r$.
\end{Lemma}

\begin{proof}
  This follows if we show that $T=\sum_j \norm{\epsilon_j} > \sqrt{c}$. So, define the function $\delta_\param \colon \Lambda M \to \R$ by
  \begin{align*}
    \delta_\param(\gamma) = \int_\gamma \norm{\gamma'-X_{H^\param}}dt.
  \end{align*}
  This is zero if and only if $\gamma$ is a 1-periodic orbit of the Hamiltonian flow. We can approximate $T$ by using $\delta_\param$ in the following way: let $\gamma$ be the closed curve which is the flow curve $(\varphi_{t-j/r})^\param(z_j)$ when $t\in[j/r,(j+1-1/k)/r]$ and when $t\in[(j+1-1/k)/r,(j+1)/r]$ it is the geodesic connecting $(\varphi_{(1-1/k)/r})^\param(z_j)$ and $z_{j+1}$. As $k$ tends to infinity $\delta_\param(\gamma)$ tends to $\Sigma$. So all we need is a lower bound on $\delta_\param$ for curves $\gamma$ which has some point in $K_\tau$. By cyclic symmetry we can assume that this point is $\gamma(0)$.

  Define a path by $\gamma_2(t)=(\varphi_{-t})^\param(\gamma(t))$. Now the flow of $X_{H^\param}$ is independent on $\param$ on $U$ and preserves it (condition H1 on the Hamiltonians). So, independently on $\param$ we get the bound
  \begin{align*}
    \delta_\param(\gamma) \geq & C \int_0^1 \norm{D_{\gamma(t)}(\varphi_{-t})^\param(\gamma'(t)-(X_{H^\param})_{\gamma(t)})}1_U(\gamma_2(t)) dt = \\
    = & C \int_0^1 \norm{\gamma'_2(t)}1_U(\gamma_2(t)) dt 
  \end{align*}
  where $C^{-1} > \norm{(D_z (\varphi_t)^\param)^{-1}}$ for all $z\in \ovl{U}$; and $1_U$ is the indicator function for $U$.

  Now we divide into two cases: Case 1: $\gamma$ lies entirely in $U$: Then we use
  \begin{align*}
    C \int_0^1 \norm{\gamma'_2(t)}1_U(\gamma_2(t)) dt = C \int_0^1 \norm{\gamma'_2(t)} dt \geq C\dist(\gamma_2(0),\gamma_2(1)) > c.
  \end{align*}
  Indeed, the flow has no periodic orbit on $\ovl{U}$ hence there is a lower bound on $\dist((\varphi_1)^\param(z),z)$ for $z$ in the compact set $\ovl{U}$.

  Case 2: $\gamma$ leaves $U$ at some point. We notice that $\gamma_2(0)=\gamma(0)\in K_\tau$, and use that
  \begin{align*}
    C \int_0^1 \norm{\gamma'_2(t)}1_U(\gamma_2(t)) dt \geq \dist(K_\tau,M-U) > 0
  \end{align*}
  because $\gamma_2$ has to move from $K_\tau$ to the complement of $U$ (both sets are compact in $M$).
\end{proof}

\begin{Corollary} \label{pjbound}
  There is a $k>0$ such that with $r$ as in Proposition~\ref{energy} we have
  \begin{align*}
    \absv*{\nabla S_r \cdot \nabla Q_j} < \sqrt{r}k \norm{\nabla S_r}^2
  \end{align*}
  when $z_j \in K_\tau$.
\end{Corollary}

\begin{proof}
  We see that $\norm{\nabla Q_j}=\norm{\nabla_{z_j} Q_j}=1$ when $z_j \in K_\tau$. So by the lemma above and Proposition~\ref{energy} we have on the set given by $z_j \in K_\tau$ that
  \begin{align*}
    \absv*{\nabla S_r \cdot \nabla Q_j} \leq \norm{\nabla S_r} < 2E \leq \sqrt{r}2c^{-1} \norm{\nabla S_r}^2.
  \end{align*}  
\end{proof}

\begin{proof}[Proof of Proposition~\ref{mainstabloc}]
  As explained in Remark~\ref{rem:Loclin:1} we will use the construction in Lemma~\ref{cutoff} to create good index pairs using $E$ and the functions $Q_j$ to define cut-off functions.

  The open manifold $\Lrb M$ has a natural compactification given by extending the definition to allow all $z_j$'s to lie in all of $M$ and allowing $e(\arz)=\beta$. We will refer to the new points in this extension as the boundary of $\Lrb M$.

  Fix $\beta \geq 66C_H^2K+8C_H^2+1$ with $K$ as in Proposition~\ref{energy}. The condition on $r$ in that proposition for this family of $H^\param$ depending on $\param\in ]0,1]$ is:
  \begin{align*}
    r > K(\param^{-1}C_H+C_\SLa^2(\beta+C_H^2)),
  \end{align*}
  which for fixed $\beta$ and $K$ (and $C_\Gamma$) we can assume is true for $r\geq 2 K C_H \param^{-1}$ for small $\param$. So fix $\param_0$ small enough for this to be true, and also fix $0<\param<\param_0$. We thus have from Proposition~\ref{energy} that for the $r$ stated in this proposition we have
  \begin{align} \label{eq:3}
    \absv*{\nabla S_r \cdot \nabla E} \leq \norm{\nabla S_r}\norm{\nabla E} \leq \sqrt{40}\norm{\nabla S_r}^2
  \end{align}
  with equality only at critical points. We will use the function
  \begin{align*}
    g = E + 10(S_r - b_\param)-\param
  \end{align*}
  as a cut-off function. Indeed, the set
  \begin{align*}
    A^E = \{g\leq 0\} \cap (S_r)^{-1}([a_\param,b_\param])
  \end{align*}
  will by construction contain all points in $(S_r)^{-1}([a_\param,b_\param])$ where $E=0$ in its interior - hence by Proposition~\ref{energy} it contains all critical points of $S_r$ in that interval. Using Equation~\eqref{eq:3} we see that on the boundary of $A^E$ where $E>0$ we have $\nabla S_r \cdot \nabla g > 0$. Also we have $E \leq 10(b_\param-a_\param)+s < 11\param C_H$ on $A^E$ (recall $b_\param-a_\param<\param C_H$ and $C_H>1$), which by Lemma~\ref{lem:Ebound} implies that
  \begin{align*}
    e(\arz) \leq 2rE(\arz) + 8C_H^2 < 22r\param C_H + 8C_H^2 \leq 66KC_H^2 + 8C_H^2 < \beta,
  \end{align*}
  since we also assumed in the proposition that $r\param \leq 3KC_H$. Hence $A^E$ is disjoint from the part of the boundary of $\Lrb M$ where $e(\arz)=\beta$.

  To keep the index pair away from the boundary of $\Lrb M$ defined by $z_j\in\partial M$ we use the functions $Q_j$ from the previous lemma to define additional cut-off functions
  \begin{align*}
    g_j =  Q_j + \frac{\tau}{3\param C_H} (S_r-b_\param) + \frac{\tau}{3}.
  \end{align*}
  The set
  \begin{align*}
    A_j = \{g_j \leq 0\} \cap (S_r)^{-1}([a,b])
  \end{align*}
  will have the critical points of $S_r$ in its interiors. Indeed, for $Q_j<-\tau$ we see that $g_j<0$ and thus $A_j$ contains the open set $(M-U)^r \cap \Lrb M$, which contains all 1-periodic orbits. Also when $g_j=0$ we have $Q_j \in [-2\tau/3, -\tau/3]$ (since $b_u-a_u<uC_H$ and $a_u\leq S_r \leq b_u$), which implies $z_j\in K_\tau$ and we may use Corollary~\ref{pjbound} and $2KC_H\param^{-1} \leq r\leq 3KC_H\param^{-1}$ to get the bound
  \begin{align*}
    \nabla S_r \cdot \nabla g_j > (-\sqrt{r}k+\frac{\tau}{3 \param C_H}) \norm{\nabla S_r}^2 \geq (-\sqrt{3KC_H}k\sqrt{\param^{-1}}+\frac{\tau}{3 C_H}\param^{-1} ) \norm{\nabla S_r}^2,
  \end{align*}
  which is greater than zero for appropriately small values of $\param$. So, we make $\param_0$ small enough to have this positive for all $\param<\param_0$.

  We have in fact proved that the critical points of $S_r$ lie in the
  set
  \begin{align*}
    A = A^E \cap \bigcap_{j} A_j,
  \end{align*}
  and since it avoids all parts of the boundary of $\Lrb M$ it is
  compact in the interior and the functions satisfies the requirements
  in Lemmas~\ref{cutoff}.
\end{proof}


\section{The Viterbo Isomorphism for Spectra.} \label{indexcalc}

In this section we calculate the homotopy type of the Conley index of the finite dimensional approximations defined in Section~\ref{Florlike} in some specific cases. This is how Viterbo originally related symplectic homology to homology of the loop space. The end result is that we calculate the stable homotopy type of the source spectrum in Equation~\eqref{eq:32} as $(\Lambda N)^{-TN}$ (even when $N$ is not oriented).

As in the construction of a single level of the spectrum $Z_a^b(H)$ we will start by assuming that $\alpha$ is the standard sub-division $\alpha_j=1/r, j\in \Z/r$ in the definition of $S_r \colon T^*\Lre N \to \R$ from Equation~\eqref{eq:7}.

Assume that $H\co T^*N \to \R$ is any Hamiltonian with
\begin{align*}
  H(q,p) = \mu \norm{p} + c
\end{align*}
for $(q,p)\notin DT^*N$. We will fix a specific Hamiltonian $H_{\mu}$ with slope $\mu$ at infinity. We can then take the convex combination homotopy $tH+(1-t)H_{\mu}$, and Corollary~\ref{cor:hominv} tells us that the spectra $Z(H)$ and $Z(H_\mu)$ associated with the total indices are homotopy equivalence using a contractible choice.

With this in mind we explicitly define $H_{\mu}(q,p)=h(\norm{p})$ where $h(t)=\frac{\mu+\epsilon}{2}t^2$ when $t<\frac{\mu-\epsilon}{\mu+\epsilon}$ for some small $\epsilon>0$ such that $[\mu-\epsilon,\mu]$ does not contain any geodesic length.

We still want $h(t)=\mu t+c$ outside $DT^*N$, but we also want $h$ to be convex so that all the $1$-periodic orbits will lie in the set where $h$ is quadratic. We achieve this by choosing $h''$ to be a smooth function with values in $[0,\mu+\epsilon]$ and constantly equal to $\mu+\epsilon$ when $t<\frac{\mu-\epsilon}{\mu+\epsilon}$ and zero when $t>1$ such that it integrates to $\mu$ over the interval $[0,1]$.

\begin{Lemma} \label{calcind}
  There exists a constant $D>0$ such that for any of the Hamiltonians $H_\mu$, with $\mu$ a geodesic length, we have that $r>D\mu$ implies existence of index pairs and
  \begin{align*}
    I(S_r,X_r) \simeq \Th(T\Lambda_r^\mu N)=(\Lambda_r^\mu N)^{T \Lambda_r^\mu N},
  \end{align*}
  where $\Lambda_r^\mu N$ is the manifold of piecewise geodesic loops in $N$ with each piece having length less than $\mu/r$
\end{Lemma}

Note that $\Th(\cdot)$ is often used for the Thom space, and is notationally convenient here. The proof of this lemma is a detailed version of Viterbos argument in \cite{MR1617648}, and we note that the only subtlety here is that we want this $D$ to be independent of $\mu$ and $\epsilon$.

\begin{proof}
  By construction of $H_\mu$ we may find a constant $C'>0$ such that $\no{H_\mu} < C'\mu$ not depending on $\epsilon$. This means that if $r>C'\mu 3\delta_0^{-1}$ then Equation~\eqref{eq:10} is satisfied and $S_r$ is defined, which by Lemma~\ref{infiniteflow} means we have good index pairs for $(S_r,X_r)$.
  
  Define $\ovl{\Lambda}_r^a N$ to be piecewise geodesics, each piece having length less than \emph{or equal} to $a/r$. We then define a discrete version of the Legendre transform, i.e. we define an embedding
  \begin{align*}
    i \colon \ovl{\Lambda}_r^{(\mu-\epsilon)} N \to T^*\Lre N,
  \end{align*}
  where $i$ is given by
  \begin{align*}
    (i(\arq))_j = \bigl(q_j,(\mu+\epsilon)^{-1} r\exp_{q_j}^{-1}(q_{j+1})\bigr).
  \end{align*}
  For $r>2\delta_0^{-1}\mu>\delta_0^{-1}(\mu+\epsilon)$ we have that $\ovl{\Lambda}_r^{\mu-\epsilon} N \subset \Lre N$, making this a section in the bundle $T^*\Lambda_r N \to \Lambda_rN$ restricted to $\ovl{\Lambda}_r^{(\mu-\epsilon)}N$. Furthermore, because 
  \begin{align*}
    \norm{(\mu+\epsilon)^{-1} r\exp_{q_j}^{-1}(q_{j+1})} \leq \tfrac{\mu-\epsilon}{\mu+\epsilon},    
  \end{align*}
  the point $(i(\arq))_j$ will lie in the set where $h$ is quadratic. In fact using the description of the flow lines for such Hamiltonians in section \ref{action} we see that we have chosen $p_j$ as the unique point in $T^*_{q_j}N$ such that $q_{j+1}^-=q_{j+1}$ (this is why we call this a discrete Legendre transform). So on the  image of $i$ all $\epsilon_{q_j}$ are 0. This implies that the image of $i$ contains all the critical points of $S_r$, because it contains all the curves with $\norm{p_j}\leq\mu-\epsilon$ and $\epsilon_{q_j}=0$ for all $j$.
  
  We will use the fiber directions ($\arp$ directions) as a normal bundle. In fact because $\epsilon_{q_j}=0$ for all $j$, Lemma~\ref{gradient} tells us that $\nabla_{p_j} S_r=0$, and because this is the only point in the fiber such that $\epsilon_{j+1}=0$, it tells us that this is the only critical point when restricting $S_r$ to the fiber. We will need that this is a global maximum in each fiber in the following very strong sense.

  Claim: for appropriate $D$ if $\arq\in\ovl{\Lambda}_r^{\mu-\epsilon}N$ is fixed then the function $S_r(\arq,\arp)$ goes to $-\infty$ as
  $\norm{\arp}$ goes to $\infty$ independently of $\mu$ and $\epsilon$.

  Proof of claim: the condition $\norm{\arp}\to \infty$ is equivalent to $\norm{p_j}\to \infty$ for some $j$. So we look at the terms in the
  definition of our finite dimensional approximation that involves $p_j$:
  \begin{align*}
    f(p_j) = \int_{\gamma_j} (\lambda - Hdt)  + p_{j+1}^-\epsilon_{q_{j+}}.
  \end{align*}
  Assume that $\norm{p_j}>1$, which means that the Hamiltonian flow of $(q_j,p_j)$ projects to a geodesic of length $\mu$. The integration part is easy to calculate and is as described in section \ref{action} $(\norm{p_j}h'(\norm{p_j})-h(\norm{p_j}))/r$, which is constant on the set $\norm{p_j}>1$. Because $\dist(q_j,q_{j+1})\leq(\mu-\epsilon)/r$ and $\dist(q_j,q_{j+1}^-)=\mu/r$, we are in the situation depicted in
  figure~\ref{negpar}.
  \begin{figure}[ht]
    \centering
    \includegraphics{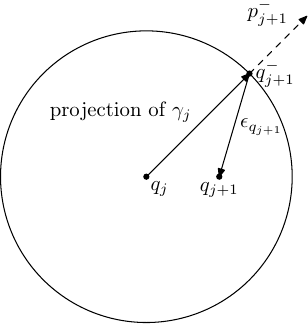}
    \caption{Position of points in $N$ when the norm of $p_j$ is
      larger than 1. The circle has radius $\mu/r$.}
    \label{negpar}
  \end{figure}
  Take the Riemannian structure we have on $N$ and multiply with $(r/\mu)^2$ such that lengths get multiplied with $r/\mu$, and take a normal chart around $q_j$ in this new metric. Then the circle in the picture is mapped to the unit circle in $\R^d$. Since the term $p_j^-\epsilon_{q_{j+1}}$ scales with the norm of $p_j$ (when $\norm{p_j}\geq 1$) it is enough to see that if $r/\mu$ is greater than some $D$ and $\norm{p_j}=1$ then this term is negative. This $D$ should be independent of $\epsilon$, because $\epsilon$ depends on $\mu$. So we have to argue that; if the Riemannian structure is flat enough then the pairing is negative for all possible $q_{j+1}$ in the open unit disc, and this is a little tricky since the pairing is of course 0 if $q_{j+1}$ is the boundary point $q_{j+1}=q_{j+1}^-$. This, of course, does not happen in our case, but we may be arbitrarily close for different $\mu$'s. So we now consider for fixed $(q_j,p_j)$ (and thus fixed $q_{j+1}^-$) with $\norm{p_j}=1$ varying $q_{j+1}$ in a small neighborhood of $q_{j+1}^-$. In fact, we will consider $q_{j+1}$'s outside the unit circle as well. In the flat case the pre-image of 0 of the term $p_j^-\epsilon_{q_{j+1}}$ is the tangent plane to the unit sphere, but for  a small perturbation it is some other sub-manifold. This manifold will by construction always contain the point $q_{j+1}=q_{j+1}^-$, which lies on the unit sphere; but moreover, it will also be parallel to the sphere. Indeed, parallel transport preserves the inner product, so the orthogonal complement of $p_{j+1}^-$ is always the tangent space to the unit sphere (even in the non-flat case). Using this tangency we see that for the Riemannian structure close enough to the flat one this manifold never enters the interior of the unit sphere. Now the pairing is negative on all of the interior, and a compactness argument gives us a choice of $D$ such that this works for all possible $(q_j,p_j)$ with $\norm{p_j}= 1$. So for such $D$, $S_r$ goes to $-\infty$ if $\norm{p_j}$ goes to $\infty$.

  Next we look at $S_r$ on the image of the embedding. Here the last term vanishes, and
  \begin{align*}
    S_r(i(\arq)) = \sum_j\int_{\gamma_j} (pdq -Hdt) =
    \frac{1}{2(\mu+\epsilon)}\sum_j r \norm{\exp_{q_j}^{-1}(q_{j+1})}^2.
  \end{align*}
  This is $(\mu+\epsilon)^{-1}$ times the energy functional
  \begin{align} \label{energyeq}
    e(\gamma) = \frac{1}{2}\int_0^1 \norm{\gamma'(t)}^2dt
  \end{align}
  evaluated on the piecewise geodesic $\arq$. This is positive and we conclude that if we look at the set defined by $S_r \geq -1$ intersected with one of the fibers, we get a bounded set diffeomorphic to a closed disc. This is true over every point in the compact set $\ovl{\Lambda}_r^{(\mu-\epsilon)} N$, so the set
  \begin{align*}
    A=\{(\arq,\arp)\mid \arq\in \ovl{\Lambda}_r^{(\mu-\epsilon)}N,
    S_r(\arq,\arp)\geq -1\}
  \end{align*}
  is compact and has points in each fiber. We also define
  \begin{align*}
    B=\{(\arq,\arp)\mid \arq\in \ovl{\Lambda}_r^{(\mu-\epsilon)}N,
    S_r(\arq,\arp) = -1\},
  \end{align*}
  which is thus the boundary sphere in each fiber. We wish to construct a new pseudo-gradient $X'$ on $A$ differing from $X_r$ only in a compact set such that $(A,B)$ is an index pair for $(S_r,X')$ and we may thus use Lemma~\ref{lem:234} to conclude that $A/B \simeq \Th(T\Lambda^{\mu}_rN)$ is the Conley index $I(S_r,X')$. We construct $X'$ only on $A$ since $A$ has no critical points on its boundary it is easy to extend $X'$ to a slightly larger open set and interpolate with $X$.
  
  Look at the gradient of $S_r$ restricted to the section we defined above, which were a constant times the energy. Here minus the gradient of the energy always flows in a direction where the longest geodesic becomes smaller or stays the same length. So in fact it flows the section strictly into the bundle over $\ovl{\Lambda}_r^{(\mu-\epsilon)} N$. It is not difficult to use that the fiber-wise gradient is non-zero away from this section to interpolate this to a pseudo-gradient that makes $B$ the exit set and $(A,B)$ and index pair. Note that it is only important what this pseudo-gradient is on the boundary of $A$ where there are no critical points.
\end{proof}

By considering the proof of the proposition on Conley indices with respect to intervals $[-1,b]$ we get the following corollary.

\begin{Corollary}
  \label{cor:2}
  The inclusion $I_{-1}^b(S_r,X) \to I_{-1}^{b'}(S_r,X)$ is homotopy equivalent to the Thom-space construction (the same as in the lemma above) on the inclusion of loops spaces.
\end{Corollary}

\begin{Remark}
  \label{rem:3}
  For the Hamiltonian $H_\mu$ the index $I_{-1}^b$ is given by
  \begin{align*}
    I_{-1}^b(S_r,X) \simeq \Th(T\Lambda_r^{\min(\sqrt{2(\mu+\epsilon) b},\mu)} N)
  \end{align*}
  Here the $x\mapsto\sqrt{2(\mu+\epsilon) x}$ is the conversion from $(\mu+\epsilon)^{-1}$ times energy to length. This is needed because the critical value corresponding to a geodesic was calculated in the proof to be $\mu^{-1}$ times the energy, and our notation for the loop spaces uses length.
\end{Remark}

A corollary of this construction which is important in \cite{immersions} is the following.

\begin{Corollary}
  \label{cor:Calculation:3}
  Let $(A_r,B_r)$ be an index pair for $(S_r,X_r)$ with $H$ as above. The inclusion $A_r \subset T^* \Lre N$ induces a map
  \begin{align*}
    A_r/B_r \to (T^* \Lre N)_+ \wedge A_r/B_r,
  \end{align*}
  which is canonically (contractible choice) homotopic to the map
  \begin{align*}
    \Th(T\Lre^{\mu-\epsilon} N) \to (T^*\Lre^{\mu-\epsilon} N)_+ \wedge \Th(T\Lre^{\mu-\epsilon} N) \subset (T^*\Lre N)_+\wedge \Th(T\Lre^{\mu-\epsilon} N).
  \end{align*}
  induced by the inclusion $DT\Lre^{\mu-\epsilon} \subset T^*\Lre N$ (and as usual identifying tangent vectors and cotangent vectors).
\end{Corollary}

Here $(-)_+$ means adding a disjoint base-point.

\begin{proof}
  The map is defined by taking the quotient of the diagonal map $A_r \to (A_r)_+ \wedge A_r/B_r = (A_r \times A_r)/(A_r \times B_r)$ and composing with the inclusion $(A_r)_+ \subset (T^*\Lre N)_+$ (still smashed with $A_r/B_r$). This is defined for any index pair, and the uniqueness proof using the negative gradient flow (in this case of $-X_r$) in the proof of Lemma~\ref{lem:indwell} extends to define a commuting diagram
  \begin{align*}
    \xymatrix{
      A_r/B_r \ar[r] \ar[d] & (T^* \Lre N)_+ \wedge A_r/B_r \ar[d] \\
      A'_r/B'_r \ar[r] & (T^* \Lre N)_+ \wedge A'_r/B'_r 
    }
  \end{align*}
  Here $(A_r',B_r')$ is an alternate index pair for $(S_r,X_r$). Furthermore, the proof of homotopy invariance in Lemma~\ref{hominv} similarly extends to a completely similar diagram, but where $(A_r,B_r)$ and $(A_r',B_r')$ are an index pair for each end of the homotopy.

  It follows that the above contractible choice identification by first changing $H$ to $H^\mu$ and then changing the pseudo-gradient to the one in the proof - proves the corollary. Indeed, the map of this type defined by using the index pair in the proof above is (up to contractible choice homotopy) the concrete map described in the corollary.
\end{proof}

When defining the generating function spectrum in Section~\ref{sec:gener-funct-spectr} we added some normal bundles to get rid of all the copies of $TN$ floating around. We see even more why this is important when comparing the proposition to the following corollary.

\begin{Corollary}
  \label{cor:Calculation:1}
  The Generating function spectrum $Z(H_\mu)$ is canonically (contractible choice) homotopy equivalent to $(\Lambda^{\mu} N)^{-TN}$. Furthermore, the inclusion $Z_{-1}^{b}(H_\mu) \to Z_{-1}^{b'}(H_\mu)$ for $0<b<b'$ is the obvious Thom-spectrum construction on the inclusion of loop spaces.
\end{Corollary}

\newcommand{\Lme}{\ovl{\Lambda}^{\mu-\epsilon}_r N}

\begin{proof}
  In Proposition~\ref{prop:Viterbo:1} we saw how adding copies of the normal bundle $\nu$ made us get an effective virtual bundle $-TN$ over the space of constants loops. Here we need to be a little more precise about the isomorphism since the Conley index is not defined until $r$ greater than or equal to some $r_0 \in \N_0$. So we now write a contractible choice formula for how to identify the stable bundle as $-TN$ and see that it will be compatible with the suspension maps (structure maps) in the spectrum.

  In the above proof we created (for a perturbed pseudo-gradient) a canonical (contractible choice) index pair canonically homeomorphic to
  \begin{align*}
    (A_r,B_r)=(DT(\Lme),ST(\Lme)).    
  \end{align*}
  Notice that we have a canonical isomorphism:
  \begin{align*}
    T(\Lme) \cong \bigoplus_{j\in \Z/r} T_{q_j}N.
  \end{align*}
  Recall the addition of $r+1$ copies of the normal bundles in Equation~\eqref{eq:1spect}. We now pick a specific way of identifying the relative Thom pair $(A_r,B_r)^{\nu^{r+1}-}$ from that equation with the pair $(\Lme,\varnothing)^{\nu \oplus \zeta^{rk}-}$. We will need the homeomorphism from Equation~\eqref{eq:1:pairident} and the isomorphism from Equation~\eqref{eq:15} to construct a bundle isomorphism:
  \begin{align} \label{eq:1}
    \pare*{\bigoplus_{j\in \Z/r} T_{q_j}N} \oplus \nu^{r+1} \cong \nu \oplus \zeta^{rk} \cong \nu \oplus (\zeta^k)^r.
  \end{align}
  Indeed,
  \begin{itemize}
  \item the first copy of $\nu$ on the left hand side is identified with the copy of $\nu$ on the right hand side (using identity).
  \item For each $j=0,\dots,r-1$ we take a parallel transport (contractible choice) along the piecewise geodesic defined by the sequence $q_{j},q_{j-1},\dots,q_0$ from $T_{q_j}N \cong T_{q_0}N$, this paired with the $(r+1-j)$th copy of $\nu$ and the isomorphism in Equation~\eqref{eq:15} produces an isomorphism
    \begin{align*}
      T_{q_j}N \oplus \nu \cong T_{q_0} N \oplus \nu \cong \zeta^k.
    \end{align*}
    Here we view $\zeta^k$ as the $j$th copy on the right hand side.
  \end{itemize}
  
  This identification is compatible with the suspensions isomorphisms in Equation~\eqref{eq:17}. Indeed, these where constructed by copying $q_0$ and putting in the new copy of $T_{q_0}N$, but also adding a copy of the normal bundle (we add this copy of the normal bundle as a new last factor). Note that the reason we choose the parallel transport above to go \emph{backwards} to $q_0$ along the string is because when adding a $q_0$ we push the remaining points forward in index ($q_j$ becomes $q_{j+1}$) - and this is precisely why this is compatible with the suspensions when increasing $r$.

  The last part of the corollary follows form Corollary~\ref{cor:2} and the fact that the construction in this proof respects the inclusion and quotients of Conley index pairs.
\end{proof}

Also, for the alternative in Remark~\ref{rem:Spectrum:2} defining alternate spectra, which we decorated with primes ${Z'}_a^b(H)$ we have the corresponding corollary.

\begin{Corollary}
  \label{cor:Calculation:4}
  The alternate Generating function spectrum $Z'(H_{\mu})$ is canonically (contractible choice) homotopy equivalent to $\Sigma^{\infty}(\Lambda^{\mu} N)_+$. Furthermore, the inclusion ${Z'}_{-1}^{b}(H_\mu) \to {Z'}_{-1}^{b'}(H_\mu)$ for $0<b<b'$ is the infinite suspension functor on the inclusion of loop spaces.  
\end{Corollary}

\begin{proof}
  Same as above, but the total bundle is now identified as trivial - hence we get a standard suspension spectrum. 
\end{proof}

\begin{Proposition}
  \label{prop:Calculation:1}
  There is a canonical (contractible choice) homotopy equivalence:
  \begin{align*}
    Z \simeq (\Lambda N)^{-TN},
  \end{align*}
  with $Z$ as in Equation~\eqref{eq:32}
\end{Proposition}

\begin{proof}
  Consider the function $f$ depicted in Figure~\ref{fig:functionf}, which we used to define $H^s$. Let $x\in I$ be the point at which the tangent with slope $\mu_N$ is tangent to this. Consider the Hamiltonian $H^{s_l}$, which defines the spectrum $Z(l)$. Let $y\in \R$ be the maximal critical level for the action below the regular level $-\ffd s_l$.
  \begin{figure}[ht]
    \centering
    \begin{tikzpicture}[scale=3] 
      \draw[->] (-0.2,-0.1) -- (1.2,-0.1) node [right] {$\norm{p_N}$};
      \draw[->] (0,0.1) -- (0,1.2);
      \draw[dotted] (0,-0.2) -- (0,0.1);    
      \draw[dotted] (0,0.275) -- (1.1,1.155);
      \draw (-0.02,1) node [left] {$s_l$} -- (0.02,1);
      \draw (-0.02,0.275) node [left] {$\ffd s_l$} -- (0.02,0.275);
      \draw (-0.02,0.375) node [left] {$-y$} -- (0.02,0.375);
      \draw[thick] (0,0.35) to[out=0,in=218.66] (0.2,0.435) to[out=38.66,in=218.66] (0.8,0.915) to[out=38.66,in=180] (1.2,1);
      \draw (1,-0.12) node [below] {$1$} -- (1,-0.08);
      \draw (0.8,-0.12) node [below] {$x$} -- (0.8,-0.08);
      \draw (0.5,-0.12) node [below] {$\eh$} -- (0.5,-0.08);
      \draw [->] (1.0,1.2) node [left] {slope $s_l\mu_N$}  -- (1.05,1.15);
      \draw (0.5,0.35) node [right] {$H_l'$}; 
    \end{tikzpicture}
    \caption{The Hamiltonian $H'_l$.}
    \label{fig:Hlprime}
  \end{figure}
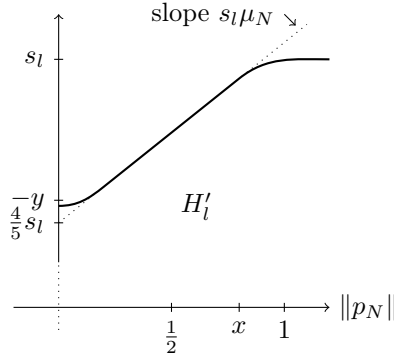
  Pick a smooth homotopy from this to a Hamiltonian $H'_l$ (illustrated in Figure~\ref{fig:Hlprime}), which satisfies:
  \begin{itemize}
  \item $H'_l$ depends only on$\norm{p_N}$,
  \item the homotopy (and thus also $H'_l$) is constantly equal to $H^{s_l}$ outside the set $D_{x}T^*N$,
  \item during the homotopy all critical values (for the associated actions) in the interval $]-\infty,y]$ are from periodic orbits outside $D_xT^*N$ (hence constant during the homotopy).
  \item $H_l'$ is convex in $\norm{p_N}$ from 0 until some small value close to zero, then it is linear with slope $s_l\mu_N$ until $x$ (and then concave by the above - this can be done smoothly due to Remark~\ref{rem:Viterbo:1}),
  \item it is quadratic close to the zero section (as $H_\mu$ above, but plus some constant),
  \end{itemize}
  Notice that the value $H_l'(q,0)$ on the zero section has to be less than $y$. If not it would violate the third bullet point.

  Such a homotopy and choice of $H'_l$ is canonical (up to a contractible choice) if we follow the canon:
  \begin{itemize}
  \item first remove the part of the definition of $H^{s_l}$ which depends on $\norm{p_L}$ by simply scaling and translating it such that the Hamiltonian becomes constantly equal to $\tfd$ inside $D_{1/2}T^*N$ (see Figure~\ref{segmentpic} as to why this does not violate the third point above),
  \item now the Hamiltonian is a function of $\norm{p_N}$ starting below $y$, then it has a convex part, and then a concave part. So we may choose the rest of the homotopy such that this is preserved.
  \end{itemize}
  The contractibility of these choices involves making the small bend at $0$ smaller and smaller.
  
  Let $\epsilon>0$ be such that at no point during the homotopy do we have a critical value in $]y,y+\epsilon]$. Now  because (definition of $y$) the interval $]y,-\ffd s_l=a_l^N]$ is regular for the action associated to $H^{s_l}$ we have that in the definition of $Z(l)$ we can replace $a_l^N$ with $y+\epsilon$. Now since there are no critical values above $-\ffd s_l$ for the action associated to $H_l'$ the above homotopy relates canonically (contractible choice)
  \begin{align*}
    Z(l) \simeq Z_{y+\epsilon}^{-\ffd s_l+\epsilon}(H'_l)
  \end{align*}
  for small $\epsilon>0$. In Proposition~\ref{mainstabloc} and Example~\ref{exm:4} we saw that if we make the bend close to zero depend on a small parameter $u>0$ we can (if we make $u$ small enough and consider appropriate $r$) get index pairs inside the sub-manifold $\Lrb D_{1/2}T^*N \subset \Lre T^*N$.

  This implies (using Lemma~\ref{lem:Homoindex:2}) that this index pair cannot see that the Hamiltonian is not linear outside of $D_{1/2}T^*N$ (continuing with the slope it already has at $\norm{p}=1/2$), and hence combining the above with Corollary~\ref{cor:Calculation:1} gives a canonical homotopy equivalence:
  \begin{align*}
    Z(l) \simeq (\Lambda^{s_l\mu_N} N)^{-TN}.
  \end{align*}
  Notice that the index pairs in Section~\ref{Florlike} are defined when $r>C\no{H'_l}$ for some $C$ (from Section~\ref{flowline}) and by definition of $C_{H_l'}$ we have $\no{H_l'}<\max(\param^{-1} C_{H_l'},C_{H_l'})$. Furthermore, the index pairs from Proposition~\ref{mainstabloc} are defined for
  \begin{align*}
    r\in [2KC_{H_l'} \param^{-1},3KC_{H_l'} \param^{-1}].
  \end{align*}
  However, in Remark~\ref{rem:Genfin:2} we fixed it such that the latter interval is contained in the solutions to the first equation (for small $\param$). So we can always find $\param$ small enough (and then $r$) so that both types of index pairs are defined at the same time. If this were not the case it would be difficult to compare them using Lemma~\ref{lem:Homoindex:2}.

  To finish the proof of the proposition we need to make this identification and the inclusion of loop spaces canonically compatible with the maps $\kappa_l \co Z(l) \to Z(l+1)$ defining $Z$.
  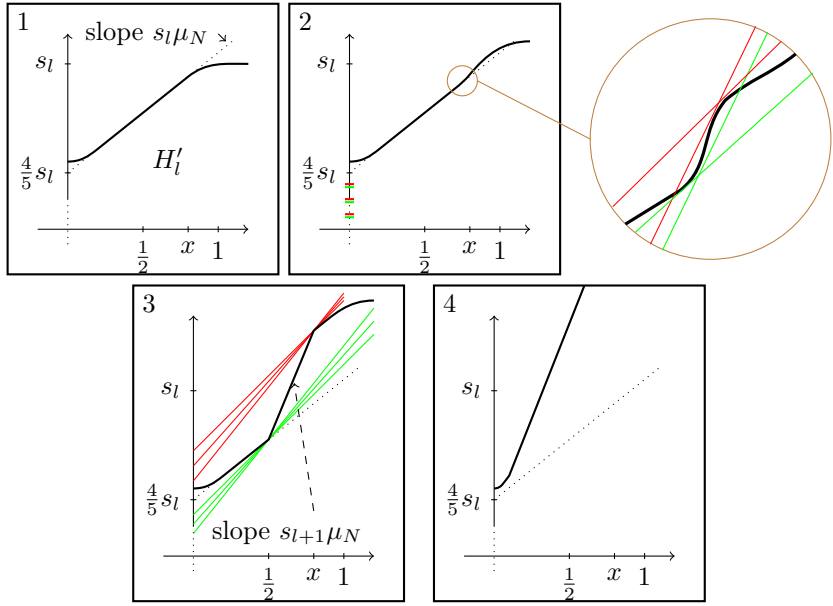
\begin{figure}[ht]
    \centering
    \begin{tikzpicture}[scale=2] 
      \draw (-0.4,1.4) node [below right] {$1$};
      \draw[thick] (-0.4,-0.4) -- (1.4,-0.4) -- (1.4,1.4) -- (-0.4,1.4) -- cycle;
      \draw[->] (-0.2,-0.1) -- (1.2,-0.1);
      \draw[->] (0,0.1) -- (0,1.2);
      \draw[dotted] (0,-0.2) -- (0,0.1);    
      \draw[dotted] (0,0.275) -- (1.1,1.155);
      \draw (-0.02,1) node [left] {$s_l$} -- (0.02,1);
      \draw (-0.02,0.275) node [left] {$\ffd s_l$} -- (0.02,0.275);
      \draw[thick] (0,0.35) to[out=0,in=218.66] (0.2,0.435) to[out=38.66,in=218.66] (0.8,0.915) to[out=38.66,in=180] (1.2,1);
      \draw (1,-0.12) node [below] {$1$} -- (1,-0.08);
      \draw (0.8,-0.12) node [below] {$x$} -- (0.8,-0.08);
      \draw (0.5,-0.12) node [below] {$\eh$} -- (0.5,-0.08);
      \draw [->] (1.0,1.2) node [left] {slope $s_l\mu_N$}  -- (1.05,1.15);
      \draw (0.5,0.35) node [right] {$H_l'$}; 
    \end{tikzpicture}
    \begin{tikzpicture}[scale=2] 
      \draw (-0.4,1.4) node [below right] {$2$};
      \draw[thick] (-0.4,-0.4) -- (1.4,-0.4) -- (1.4,1.4) -- (-0.4,1.4) -- cycle;
      \draw[->] (-0.2,-0.1) -- (1.2,-0.1);
      \draw[->] (0,0.1) -- (0,1.2);
      \draw[dotted] (0,-0.2) -- (0,0.1);    
      \draw[dotted] (0,0.275) -- (1.1,1.155);
      \draw (-0.02,1) node [left] {$s_l$} -- (0.02,1);
      \draw (-0.02,0.275) node [left] {$\ffd s_l$} -- (0.02,0.275);
      \draw[thick] (0,0.35) to[out=0,in=218.66] (0.2,0.435) to[out=38.66,in=218.66] (0.7,0.835) to[out=38.66,in=230] (0.8,0.930) to[out=50,in=180] (1.2,1.15);
      \draw (1,-0.12) node [below] {$1$} -- (1,-0.08);
      \draw (0.8,-0.12) node [below] {$x$} -- (0.8,-0.08);
      \draw (0.5,-0.12) node [below] {$\eh$} -- (0.5,-0.08);
      \draw[brown] (0.75,0.8875) circle (0.1);
      \draw[brown] (0.85,0.8875) -- (1.6,0.5);
      \draw[brown] (2.4,0.5) circle (0.8);
      \draw[very thick] (1.834,-0.066) -- (2.2,0.16) to[out=38.66,in=230] (2.5,0.76) to[out=40,in=218.66] (2.966,1.066);
      \draw[red] (2.0,-0.2) -- (2.7,1.25);
      \draw[green] (2.08,-0.24) -- (2.78,1.21);
      \draw[red] (1.75,0.05) -- (2.87,1.15);
      \draw[green] (1.90,-0.12) -- (3.08,0.94);
      \draw[red,thick] (-0.03,0.2) -- (0.03,0.2);
      \draw[green,thick] (-0.03,0.18) -- (0.03,0.18);
      \draw[red,thick] (-0.03,0.1) -- (0.03,0.1);
      \draw[green,thick] (-0.03,0.08) -- (0.03,0.08);
      \draw[red,thick] (-0.03,0.0) -- (0.03,0.0);
      \draw[green,thick] (-0.03,-0.02) -- (0.03,-0.02);
    \end{tikzpicture}

    \begin{tikzpicture}[scale=2] 
      \draw (-0.4,1.7) node [below right] {$3$};
      \draw[thick] (-0.4,-0.4) -- (1.4,-0.4) -- (1.4,1.7) -- (-0.4,1.7) -- cycle;
      \draw[->] (-0.2,-0.1) -- (1.2,-0.1);
      \draw[->] (0,0.1) -- (0,1.5);
      \draw[dotted] (0,-0.2) -- (0,0.1);    
      \draw[dotted] (0,0.275) -- (1.1,1.155);
      \draw[red] (0,0.6) -- (1.0,1.6);
      \draw[red] (0,0.5) -- (1.0,1.625);
      \draw[red] (0,0.4) -- (1.0,1.65);
      \draw[green] (0,0.175) -- (1.2,1.375);
      \draw[green] (0,0.1125) -- (1.2,1.4625);
      \draw[green] (0,0.05) -- (1.2,1.55);
      \draw (-0.02,1) node [left] {$s_l$} -- (0.02,1);
      \draw (-0.02,0.275) node [left] {$\ffd s_l$} -- (0.02,0.275);
      \draw[thick] (0,0.35) to[out=0,in=218.66] (0.2,0.435) to[out=38.66,in=218.66] (0.5,0.675)  -- (0.8,1.4) to[out=38.66,in=180] (1.2,1.6);
      \draw (1,-0.12) node [below] {$1$} -- (1,-0.08);
      \draw (0.8,-0.12) node [below] {$x$} -- (0.8,-0.08);
      \draw (0.5,-0.12) node [below] {$\eh$} -- (0.5,-0.08);
      \draw [->,dashed] (0.8,0.2) node [below] {slope $s_{l+1}\mu_N \qquad$}  -- (0.67,1.05);
    \end{tikzpicture}
    \begin{tikzpicture}[scale=2] 
      \draw (-0.4,1.7) node [below right] {$4$};
      \draw[thick] (-0.4,-0.4) -- (1.4,-0.4) -- (1.4,1.7) -- (-0.4,1.7) -- cycle;
      \draw[->] (-0.2,-0.1) -- (1.2,-0.1);
      \draw[->] (0,0.1) -- (0,1.5);
      \draw[dotted] (0,-0.2) -- (0,0.1);
      \draw[dotted] (0,0.275) -- (1.1,1.155);
      \draw (-0.02,1) node [left] {$s_l$} -- (0.02,1);
      \draw (-0.02,0.275) node [left] {$\ffd s_l$} -- (0.02,0.275);
      \draw[thick] (0,0.35) to[out=0,in=230] (0.1,0.435) -- (0.6,1.7);
      \draw (1,-0.12) node [below] {$1$} -- (1,-0.08);
      \draw (0.8,-0.12) node [below] {$x$} -- (0.8,-0.08);
      \draw (0.5,-0.12) node [below] {$\eh$} -- (0.5,-0.08);
    \end{tikzpicture}
    \caption{Homotopy from $H_l'$ to $H_{l+1}'+c$.} \label{fig:Hhomotopy}
  \end{figure}
  To argue this we consider the homotopy of Hamiltonians (illustrated in Figure~\ref{fig:Hhomotopy}) defined by the following steps.
  \begin{itemize}
  \item Firstly we create a small convex bend followed by a small concave bend right at the point $x$. This can create a lot of new critical points, but they are divided into canceling pairs (if the critical level was a single non-degenerate point it would be a canceling pair of cells in a CW structure on the Conley index). In the figure we have sketches associated tangents coming from the top concave part in red, and their canceling partner with the same slope from the bottom convex part of the bend in green.
  \item To begin with the associated green and red critical values are all greater than any critical value before and appear in pairs as indicated on frame 2 (where we have not drawn the associated tangents - as we have in frame 3 - but only their intersection with the 2. axis is indicated).
  \item We then slide the convex bend down creating a linear part with slope $S_{l+1}\mu_N$ and pushing up the outer bend. In doing so the intersection between the 2. axis and the red lines will pass through the value $\ffd s_l$ and the value of the Hamiltonian at the zero section.
  \end{itemize}
  At the end of this homotopy we have a Hamiltonian which is a translation of $H_{l+1}'$. Now, let $b>0$ be greater than any of the critical values during this homotopy, and let $a$ be regular and slightly smaller than the critical value corresponding to the top most red intersection in the last frame (i.e. corresponding to a point on the 2. axis right above this intersection point). Notice that through this homotopy all the critical value associated with tangents outside of $x$ intersect the 2. axis above the red lines. It follows that the map $\kappa_l \co Z(l) \to Z(l+1)$ can now be identified with the homotopy and quotient on Conley indices (plus the untwisting normal bundles $\nu^{r+1}$) given by the spectra maps:
  \begin{align} \label{eq:66}
    \kappa_l' \co Z(l) \simeq Z_{-y+\epsilon}^{\ffd s_l} (H_l') \simeq Z_{-y+\epsilon}^{b} (H_l') \simeq Z_{a}^b(H_{l+1}') \to Z_{-y+\epsilon}^{b}(H_{l+1}') \simeq Z(l+1).
  \end{align}
  Here the first homotopy equivalence is the one above. The second is simply extending the interval to be much larger, which for $H_l'$ does not include any new critical points. The second uses Corollary~\ref{cor:hominv} and the Hamiltonian homotopy we described above, which for large $b$ has the creation of red and green points within the interval, and we need to smoothly change the bottom value of the interval from $-y+\epsilon$ to $a$ so that it is always a regular value (this is possible by the concavity of the top part, and can be done precisely as in Section~\ref{sec:gener-funct-spectr} where we made this value equal to the intersection of the unique tangent with constant slope). The third map is the map collapsing away the red critical points. The last map translated the Hamiltonian and the regular values and then uses the identification from the first part of the lemma - this uses a homotopy from $H^{s_{l+1}\mu_L}$ to $H_{l+1}'$ backwards. This identifies $\kappa_l'$ as canonically (contractible choice) homotopy equivalent to $\kappa_l$. Indeed, the concatenation of all three homotopies can be undone while keeping the values we used regular.
 
  Now we show that $\kappa'_l$ is homotopy equivalent to the inclusion of loops (with Thom-constructions on top). Indeed, instead of including the new pairs directly and making the interval of action large (second step in $\kappa'_l$) we have an alternative: we can keep the interval very small around the bend at the zero section, and wait to include anything until it gets very close. In fact, we can choose not to include any of the green critical points until the very end (even when they get close), and we only include the red ones in the interval when we have to - that is when they actually enter the small interval. However, very soon after they have entered they leave the other end of the interval and gets collapsed away. Again Proposition~\ref{mainstabloc} tells us that by narrowing the bend and action interval at $0$ (and consequently prolonging the linear part from the $0$-bend to the ``green'' convex bend) the index pair does not see the difference before and after the red point passes through. So by passing all the red points through we get a sequence of homotopy equivalences, and then including the green points at the end simply is the inclusion we already identified in Corollary~\ref{cor:Calculation:1}.

  This alternate description is the same map as in Equation~\eqref{eq:66}. Indeed, including critical values at the top of an interval and collapsing away critical points at the bottom of an interval commutes. Furthermore, a very similar argument shows why the map coming from the concatenations of two such homotopies is canonically (up to contractible choice) identified with a single one going from slope $s_l\mu_N$ to $s_{l+2}\mu_N$.
\end{proof}

Similarly we have for the alternate spectra $Z'$ (Remark~\ref{rem:Spectrum:2} and Corollary~\ref{cor:Viterbo:1}) the following corollary.

\begin{Corollary}
  \label{cor:Calculation:2}
  We have a canonical homotopy equivalence
  \begin{align*}
    Z' \simeq \Sigma^{\infty}(\Lambda N)_+.
  \end{align*}
\end{Corollary}


\section{Stabilization of Generalized Finite Dimensional
  Approximations}\label{stabloc} 

In this section we describe natural stabilizations of the finite dimensional approximations defined in Section~\ref{loclinsec}. Indeed, we will cross $M$ with a standard symplectic disc $D^{2k}$, and by using a Hamiltonian with a single 1 periodic orbit on $D^{2k}$ the Hamiltonian Floer homology will be unchanged. However we will need to be able to manipulate the gradient a bit to be able to describe precisely what adding this extra factor does for our finite dimensional approximations and the Conley indices they define. Indeed, we will prove that: under a certain product assumption the situation is very similar to putting a trivial vector bundle on $M$ with a non-degenerate quadratic form. We will then prove that the Conley index will change by the relative Thom space construction using the negative eigenbundle of said quadratic form, which in general need \emph{not} be a trivial bundle.

Let $(M,\partial M)\subset M'$ and $H_M \colon M \to \R$ be as in Section~\ref{sec:gener-finite-dimens} (the Hamiltonian is decorated with $M$ to distinguish it in the following). Define $P=M\times (D^{2k},\omega_0)$ and $P'=M'\times \R^{2k}$ for some $k$ and $\lambda_0=ydx$ the standard Liouville 1-form on $\R^{2k}$. Also define $H \colon P' \to \R$ by $H(z_1,z_2) = H_M(z_1) + H_D(z_2)$, where 
\begin{align*}
  H_D(z_2) = \norm{z_2}^2.
\end{align*}
The Hamiltonian flow for $H_D$ is circular around 0 with revolution time $2\pi$, but we only flow for a time period of 1, so the only 1-periodic orbit is $0$, and this orbit has action 0. So the 1-periodic orbits for $H$ are the 1-periodic orbits for $H_M$ on the first factor and constantly equal to $0$ on the second factor, and the critical value of these orbits are the same as on $M$.

As before we need a compatible Riemannian structure $g$ on $P$, and in fact we define this as the product of such a structure on $M$ and the standard one on $\R^{2k}$. The corresponding finite version of the loop space $\Lrb (M\times D^{2k})$ will consist of curves denoted by $\arz=(\arzet,\arzto)$. So that $\arzet$ consists of $r$ points in the interior of $M$, and $\arzto$ consists of $r$ points in the interior of $D^{2k}$. We denote the energy on loops on $M$ by $e_M$ and that on $D$ by $e_D$. We see that we have:
\begin{align*}
  e(\arz) = e(\arzet,\arzto) = e_M(\arzet) + e_D(\arzto).
\end{align*}
Similarly we can define the relative energy from Equation~\eqref{eq:2EE} factor wise and we have
\begin{align*}
  E(\arz) = E_M(\arzet)+E_D(\arzto).
\end{align*}
Indeed, the flow and everything is defined factor wise. We now assume the $K>1$ from Proposition~\ref{prop:Calculation:1} works for all three domains $M\subset M'$, $P\subset P'$ and $D^{2k} \subset \R^{2k}$ simultaneously (the maximum of the three associated $K$'s). Notice, that this $K$ depends on the symplectic and Riemannian structures only.

For us to define finite dimensional approximations as in the Section~\ref{sec:gener-finite-dimens} we still need a section
\begin{align*}
  \SLa \colon P \to \La (T P),
\end{align*}
which we this time do not assume to be time-dependent (and it will be clear in this section why we did so before). In fact, we will assume even more regularity than this.
\begin{Definition} \label{product}
  The section $\SLa$ is said to be of \emph{product type} if it factors through the projection to $M$ and the inclusion 
  $\La(TM) \times \La(k) \subset \La(TP)$, where $\La(k)$ is the Grassmannian of Lagrangian subspaces in $\R^{2k}$.
\end{Definition}
Factoring through the projection to $M$ is equivalent to the section not depending on the second (contractible) factor $D^{2k}$. Factoring through the inclusion is equivalent to all the Lagrangians splitting as direct sums of two Lagrangians, one in each factor. For the rest of this section, $\SLa$ will be of product type and $S_r$ will be the finite dimensional approximation defined as in the previous section on $\Lrb P$. The assumptions we now have on $\SLa$ imply that we can write
\begin{align*}
  \SLa_t(z_1,z_2) = \SLa^1(z_1)\oplus\SLa^2(z_1).
\end{align*}
That is, time independent so we remove the $t$, not depending on $z_2$, and it is a direct sum of a section
\begin{align*}
  \SLa^1 \co M \to \La(TM),
\end{align*}
and what could be heuristically called a ``twisting'' map:
\begin{align} \label{eq:47}
  \SLa^2 \co M \to \La(k). 
\end{align}
Notice that the bound $C_{\SLa}$ from Equation~\eqref{eq:30} is a bound on the derivative of this ``twisting'' map.

In this case $S_r$ splits into two factors
\begin{align*}
  S_r(\arzet,\arzto) = S_r^M(\arzet) +  S_r^D(\arzet,\arzto)
\end{align*}
Here $S_r^M$ is the function defined in the previous section on $\Lrb M$ by only using the first factor $\SLa^1$ of $\SLa$. The function $S_r^D(\arzet,-)$ is the finite dimensional approximation on $D^{2k}$ defined by the Hamiltonian $H_D$, but using the Lagrangians given by the second factors $\SLa^2((z_1)_j)$, which depends on $j$ (time dependence from the point of view of the second factor). This is where the bounds we assumed on any time-dependent $\SLa$ in the previous section comes in. 

\begin{Lemma}
  \label{lem:Stabloc:1}
  With $K$ as in Proposition~\ref{energy} (for all three domains) we have for 
  \begin{align}\label{eq:57}
    r>K\pare*{\no{H} + (\beta C_\SLa^2)^2(\beta+\nor{H}^2)} 
  \end{align}
  that
  \begin{align} \label{eq:58}
    \norm{\nabla E}^2 \leq 20E \leq 40\norm{\nabla S_r}^2 \leq 80E \notag \\
    \norm{\nabla E_M}^2 \leq 20E_M \leq 40\norm{\nabla S^M_r}^2 \leq 80E_M \\
    \norm{\nabla E_D}^2 \leq 20E_D \leq 40\norm{\nabla_{\arzto} S^D_r}^2 \leq 80E_D \notag
  \end{align}
\end{Lemma}
Note that the only difference in the formula for $r$ (compared to Proposition~\ref{energy}) is that $C_\SLa$ is replaced by $\beta C_\Sla^2$.

\begin{proof}
  The hard part here is the last of the three inequalities in Equation~\eqref{eq:58}. Indeed, the first two are simply the old proposition for these two approximations (indeed, we always assumed $\beta>1$ and $C_\SLa>1$ and so the $r$ satisfying Equation~\eqref{eq:57} also satisfy Equation~\eqref{eq:53}). The third and last is a little more subtle: a priori we have that how large we need $r$ may depend on $\arzet$ in the first factor. However, since the map in Equation~\eqref{eq:47} has derivative bounded by $C_\SLa$ it follows that for any $\arzet \in \Lrb M$ we can assume that the piece-wise geodesic defined in $\La(k)$ by the points $\SLa^2((z_1)_j)$ has energy bounded by $(C_\SLa)^2\beta$ (since the energy of $\arzet$ is bounded by $\beta$). Now we simply consider the non-compact family of all piece-wise (number of pieces not fixed) geodesics in $\La(k)$ with this bound on the energy. Then this family satisfy the bounds we assumed in Equation~\eqref{eq:31}, and Equation~\eqref{eq:30} is satisfied (in fact by the $0$ bound) since the section does not depend on $z_2$. This means that independent of what $\arzet\in \Lrb M$ is we have the bounds in the old proposition on the last factor alone (now viewed as a time dependent section), but the old bound $C_\SLa$ had to be replaced with the new bound $\beta C_\SLa^2$.
\end{proof}

This was the most important reason for not allowing $K$ to depend on $\beta$. Indeed, that would have introduced some circular reasoning here.

A big reason for the subtleties involving the second factor $S_r^D$ is that we have a mixing of the gradients:
\begin{align} \label{gradar}
  \begin{array}{rl}
  \nabla_{\arzet} S_r &= \nabla S_r^M + \nabla_{\arzet} S_r^D \\
  \nabla_{\arzto} S_r &= \nabla_{\arzto} S_r^D
  \end{array}
\end{align}
The second term in the first line is a little troublesome, and the next part is to get rid of this ``mixed'' part of the gradient, and then use this to prove that we basically get a relative Thom construction on the Conley indices.

To be able to actually have index pairs we now consider the ``narrowing'' case in Section~\ref{loclinsec}. That is we replace $H_M$ with a family $H_M^\param$ and a narrowing interval $[a_\param,b_\param]$ (satisfying H1 and H2 from Section~\ref{loclinsec}). However, here in the product case we define
\begin{align*}
  H^\param(z_1,z_2) = H_M^\param(z_1) + H_D(z_2).
\end{align*}
So, the Hamiltonian on the second factor will not be narrowed. Indeed, we don't have to do this since the ``interval'' of critical action values is as narrow as intervals gets. Indeed, any periodic orbit for $H^\param$ is a periodic orbit on $M$, but $0$ on the second factor. So, without narrowing this second factor we in fact have that $H^\param$ with values $[a_\param<b_\param]$ does satisfy H1 and H2. However, since we assumed that $\partial M$ was smooth, and this is not exactly the case for $\partial P$ (it has corners), we will see a slight elaboration to compensate for this in the argument below.

Again we will need functions that we can use to create cut-off functions keeping index pairs away from the boundary of $P$. Precisely as explained in Remark~\ref{rem:Loclin:1}, and since the boundary has corners it is convenient to do each part separately and define:
\begin{align*}
  Q_j^M(\arz) = - \dist((z_1)_j,\partial M) \qquad \textrm{and} \qquad Q_j^D(\arz) = - \dist((z_2)_j,S^{2k-1}).
\end{align*}
We similarly define $K_\tau^M$ and $K_\tau^D$ as in Equation~\eqref{eq:48}.
\begin{Lemma}
  \label{lem:Stabloc:2}
  Similarly to Corollary~\ref{pjbound} we have a constant $k>0$ such that if $r$ is as in Lemma~\ref{lem:Stabloc:1} then we have
  \begin{align*}
    \nabla S_r^M \cdot \nabla Q_r^M \leq k\sqrt{r} \norm{\nabla S_r^M}^2
  \end{align*}
  when $(z_1)_j \in K_\tau^M$ and
  \begin{align*}
    \nabla_{\arzto} S_r^D \cdot \nabla Q_r^D \leq k\sqrt{r} \norm{\nabla S_r^D}^2
  \end{align*}
  when $(z_2)_j \in K_\tau^D$.  
\end{Lemma}

\begin{proof}
  Since $r$ satisfies Equation~\eqref{eq:57} it satisfies Equation~\eqref{eq:53} and thus the first is simply Corollary~\ref{pjbound}. For the second we use that Lemma~\ref{lem:Loclin:3} provides a lower bound on $E_D>c/r$ which combined with the inequality in Equation~\eqref{eq:58} gives
  \begin{align*}
    \nabla_{\arzto} S_r^D \cdot \nabla Q_r^D \leq \norm{\nabla S_r^D} < 2\sqrt{r}c^{-1}\norm{\nabla S_r^D}^2.
  \end{align*}
\end{proof}

We now have all the functions needed to create good index pairs for $S_r$ on $\Lrb P$ using its gradient, but we will need to deform the gradient through pseudo-gradients to obtain a homotopy that essentially removes the unwanted mixed term in Equation~\eqref{gradar}. We also want to scale the term in the first factor so that we can argue that we essentially get a Conley index that fibers over the first factor. Hence we will prove a lemma similar to Proposition~\ref{mainstabloc} but with a family of pseudo-gradient suited for this.

\begin{Lemma} \label{7}
  With $K$ as in Lemma~\ref{lem:Stabloc:1} there exists a $\beta>0$ large enough and an $\param_0>0$ small enough so that the following holds.
  
  For any $0<\param<\param_0$, $r\in[2KC_H \param^{-1} , 3KC_H \param^{-1}]$, $0< t_1 \leq 1$ and $0 \leq t_2 \leq 1$ there exist a good index pair for the index $I_{a_\param}^{b_\param}(S_r,X)$ of $S_r\colon\Lrb  P \to \R$, where
  \begin{align*}
    X = (t_1 \nabla S_r^M + t_2 \nabla_{\arzet}S_r^D) \oplus \nabla_{\arzto} S_r^D.
  \end{align*}
\end{Lemma}

The proof is very similar to the proof of Proposition~\ref{mainstabloc}, but with a few extra complications.

\begin{proof}
  Again we write an explicit formula for $\beta$. however, the formula is a little different to accommodate the new proof, but the idea is essentially the same:
  \begin{align*}
    \beta =  104 K C_H^2 + 8C_H^2 + 1.
  \end{align*}
  The splitting of $E$ into $E_M+E_D$ makes the gradient split into:
  \begin{align*}
    \nabla E = \nabla E_M \oplus \nabla E_D.
  \end{align*}
  Now bounding the mixed term in Equation~\eqref{gradar} can be done by using Corollary~\ref{cor:Loclin:4}. Indeed, this term comes only from the fact that when moving points on $M$ the Lagrangians on the other factor changes. We have $C_\SLa$ as a bound (Equation~\eqref{eq:30}) on how fast the Lagrangians can change depending on $\arzet$. We also have a bound $C'$ in Corollary~\ref{cor:Loclin:4} on how much changes in the Lagrangians change $S_r$. Combined we get
  \begin{align} \label{eq:59}
    \norm{\nabla_{\arzet} S_r^D} \leq C_{\SLa}C'\sum_j \dist((z_2)^-_j,(z_2)_j)^2 = C_\SLa C'E_D
  \end{align}
  Note that this is really a simpler version of Equation~\eqref{eq:20}. Indeed, since the factor $D^{2k}$ has the standard structure we do not need the constant $C_M$ to translate bounds between the structures.

  Again we see that for $r>2KC_H\param ^{-1}$ and small $\param$ we have
  \begin{align*}
    r > 2KC_H \param^{-1} > K(C_H\param^{-1}+(\beta C_\SLa^2)^2(\beta + C_H^2)),
  \end{align*}
  So Equation~\eqref{eq:59} together with Lemma~\ref{lem:Stabloc:1} now proves:
  \begin{align*}
    \absv{X \cdot \nabla E} & \leq \bigr(t_1\norm{\nabla_{\arzet} S_r^M} + t_2\norm{\nabla_{\arzet}S_r^D}\bigl) \norm{\nabla E_M} + \norm{\nabla_{\arzto} S_r^D}\norm{\nabla E_D} \leq \\
    & \leq \bigr(t_1\sqrt{2E_M} + t_2C_{\SLa}C'E_D\bigl)\sqrt{20E_M} + \sqrt{2E_D}\sqrt{20E_D} \leq \\ 
    & \leq 7(t_1 E_M+t_2C_{\SLa}C'E_D\sqrt{E_M}+E_D) \leq 8(t_1E_M + E_D),
  \end{align*}
  The last inequality follows for large $r$ (small $\param$) where the middle term is much smaller than the last terms ($E_M$ is very small by Lemma~\ref{lem:Ebound}). Similarly we have
  \begin{align} \label{eq:60}
      X \cdot \nabla S_r & \geq (t_1 \nabla_{\arzet} S_r^M + t_2 \nabla_{\arzet} S_r^D) \cdot (\nabla_{\arzet} S_r^M + \nabla_{\arzet} S_r^D) + \norm{\nabla_{\arzto} S_r^D}^2 \geq \notag \\
      & \geq t_1 E_M/2 - (t_1+t_2)\sqrt{E_M/2}C_{\SLa}C'E_D - t_2(C_{\SLa}C'E_D)^2 + E_D/2 \geq \\
      & \geq t_1/2 E_M + E_D/3 \geq \tfrac13(t_1E_M+E_D)\notag
  \end{align}
  Again the two middle terms are swallowed by the last term (both $E_M$ and $E_D$ are small for large $r$). So by making $\param_0$ small we make $r$ larger and this makes $E=E_M+E_D$ small. So for small $\param_0$ we can assume:
  \begin{align*}
    X \cdot \nabla E \leq 8(t_1 E_M + E_D) \leq 24(X \cdot \nabla S_r).
  \end{align*}
  This proves that $X$ is a pseudo-gradient, because at non-critical points we have $E_M+E_D>0$. It also proves that we can use $E$ as a cut-off function in the same way we did in the proof of Proposition~\ref{mainstabloc}. Indeed, this time we can use
  \begin{align*}
    g = E + 25(S_r - b_\param)-\param
  \end{align*}
  as a cutt-off function. Indeed, again since critical points has $E=0$ these are inside the set $g<0$ and $X(g) = X \cdot \nabla g > 0$ by the above. So using this function to cut-off as in Lemma~\ref{cutoff} we have $g\leq 0$ implies $E \leq 25(b_\param-S_r)+\param \leq 26 \param C_H$, and using this with Lemma~\ref{lem:Ebound}
  \begin{align*}
    e(\arz) \leq 2rE(\arz) + 8C_H^2 < r 52 \param C_H + 8C_H^2 \leq 104 K C_H^2 + 8C_H^2 < \beta,
  \end{align*}

  For the other part of the boundary we use Lemma~\ref{lem:Stabloc:2} and $\norm{\nabla Q_j^M}\leq 1$ and get on the set $Q_j^M(\arz)>-\tau$ that 
  \begin{align*}
    X \cdot \nabla Q_j^M & \leq t_1 \nabla_{\arzet} S_r^M \cdot \nabla Q_j^M + t_2 C_{\SLa}C'E_D\norm{\nabla Q_j^M} \leq \\
    & \leq t_1\sqrt{r}k\norm{\nabla_{\arzet} S_r^M}^2+ t_2C_{\SLa}C'E_D \leq \\
    & \leq \sqrt{r}k\pare*{t_12E_M + \frac{t_2C_\SLa C'}{\sqrt{r}k}E_D} \leq \\
    & \leq 6\sqrt{r}k(X\cdot \nabla S_r).
  \end{align*}
  The last inequality follows from Equation~\eqref{eq:60} and large $r$ (small $s$), this is still independent of $t_1$ and $t_2$ as the lemma stipulates.
  
  Similarly (yet easier) we get from Lemma~\ref{lem:Stabloc:2} (when $Q_j^D>-\tau$) that
  \begin{align*}
    X \cdot \nabla Q_j^D = \nabla S_r^D \cdot \nabla Q_j^D < \sqrt{r}k \norm{\nabla S_r^D}^2 \leq \sqrt{r}k(X\cdot \nabla S_r^D).
  \end{align*}
  Now we have enough cut-off functions to get a compact pair in the interior of $\Lrb (M\times D^{2k})$, just as in the proof of Proposition~\ref{mainstabloc}. Indeed, the construction of cut-off functions using $Q_j^D$ and $Q_j^M$ is completely analogous to the construction in that proof.
\end{proof}

\begin{Remark} \label{quadrem}
  Because $\nabla_{v z_2} H_D = v \nabla_{z_2} H_D$ for $v\in \R_+$, we see that flow curves for the Hamiltonian flow of $H_D$ is preserved under scaling. So if $\gamma$ is a flow curve then $v \gamma$ is a flow curve. Since the Lagrangian $\SLa^2$ does not depend on $\arzto$ the L-curves scale as well. This means that the curve over which we integrate $\lambda_0$ (in the formula for $S_r^D$) scales proportionally with $\arzto$, so the integral of the canonical 1-form scales quadratically. Furthermore, $H_D$ is quadratic. So we conclude
  \begin{align} \label{quad}
    S_r^D(\arzet,v\arzto) = v^2 S_r^D(\arzet,\arzto).
  \end{align}
  Because this is a smooth function it must be equal to its Hessian at $0$ (for fixed $\arzet$). So $S_r^D(\arzet,-)$ is in fact a quadratic form in $\arzto$, and we can thus extend it uniquely to all of $(\R^{2k})^{r}$ - as that quadratic form - although we will not need this until the next section.

  If the critical point $0$ were degenerate for this quadratic form it would not be an isolated critical point. So in fact this is  a non-degenerate quadratic form.
\end{Remark}

\begin{Definition}\label{zetadef}
  Let $W^- \to \Lrb M$ be the vector bundle with fiber at $\arzet$ the negative eigenspace for the quadratic form $S_r^D(\arzet,-)$.
\end{Definition}

\begin{Lemma} \label{thomsusp}
  Assume that $K,\beta,\param_0,\param$ and $r$ satisfy the conditions in the previous lemma. Then the Conley index $I_{a_\param}^{b_\param}(S_r)$ is canonically (contractible choice) the relative Thom space of $W^-$ on a Conley index pair for $I_{a_\param}^{b_\param}(S_r^M)$.
\end{Lemma}

Note that here we claim this for the gradient, but to be able to prove it we use the above family of pseudo-gradients.

\begin{proof}
  First we use homotopy invariance from Lemma~\ref{hominv} to realize that if we can prove this for $X$ as in the above lemma with $t_2=0$ and $t_1$ very small the lemma will follow. So, choose a good index pair $(A,B)$ for $(S_r^M,X)$, with respect to the narrow interval $[a_\param,b_\param]$. We will extend this to an index pair for $S_r$ with $t_1$ small enough. Let $W_{\arzet}^{\pm}$ be the negative/positive eigenbundle of $S_r^D(\arzet,-)$. It is easy to construct index pairs very close to zero for a non-degenerate quadratic form on $D^{2kr}$, so we do this fiber-wise
  \begin{align*}
    A_{\arzet}&=D_\epsilon E_{\arzet}^-\times D_\epsilon E_{\arzet}^+ \\
    B_{\arzet}&=S_\epsilon E_{\arzet}^-\times D_\epsilon E_{\arzet}^+.
  \end{align*}
  Since $e(A) \in [0,\beta[$ and $A$ is compact we have $e(A)\in[0,\beta-c]$, so we can find an $\epsilon>0$ such that $A_{\arzet}$ is contained in $\Lrb (M\times D^{2k})$ for all $\arzet\in A$. Define
  \begin{align*}
    A' &= \bigcup_{\arzet \in A} A_{\arzet} \\
    B' &= (\bigcup_{\arzet \in B} A_{\arzet})\cup(\bigcup_{\arzet\in A}
    B_{\arzet})
  \end{align*}
  for such an $\epsilon$.

  Claim: For sufficiently small $t_1$ (and $t_2=0$) this is an index pair for $(S_r,X)$. I1 and I2 from the definition of index pair have been taken care of. I3 is because critical points of $S_r$ are of the form $(\arzet,0)$, where $\arzet$ is a critical point for $S_r^M$. To get I4 we need to carefully pick $t_1$. Indeed, $t_1$ controls the speed of the flow on the first factor (the base). For $t_1$ equal to zero (where $X$ is not a pseudo-gradient) any point in $B'$ will on the second factor flow entirely out of $\Lrb (M\times D^{2k})$ (except if $\arz \in B$ so that $B'$ is all of $A_{\arz}$ in this ``fiber''). Because $A'$ is compact we can choose $t_1$ very small such that the escaping above for points in $B'$ is not changed, and such that other points in the boundary of $A'$ still flows directly into $A'$. However for small $t_1$ we see by the fact that the flow projected to the base $\Lrb M$ does not depend on $\arzto$ that $B'$ is precisely the exit set of $X$

  The quotient $A'/B'$ is the wanted Conley index.
\end{proof}


\section{Quadratic Forms Associated With the Action in $\R^{2k}$.} \label{sec:quadr-forms-assoc}

In the previous section we saw how the negative eigenspace of the quadratic form (described in Remark~\ref{quadrem}) given by the finite dimensional approximations (with the Hamiltonian $H_D(z) = \norm{z}^2$)
\begin{align*}
  S_r \co (\R^{2k})^r \to \R
\end{align*}
are important for understanding the Conley index of the generalized approximations from Section~\ref{loclinsec} on products as in Section~\ref{stabloc}. This approximation depended on a single time dependent Lagrangian (i.e. a loop in $\La(k)$), and we had a family (in that section parameterized by a subspace in $\Lrb M$) of such - with bounded energy. So in this section we assume that we have a map $B \to \Lambda \La(k)$, which we for convenience write as:
\begin{align*}
  l_b \in \Lambda \La(k) \qquad \textrm{for each }b\in B,
\end{align*}
and we will be assuming that the energy of each $l_b$ is bounded by some $C_l$. We let $S_r^b$ denote the quadratic form defined by finite dimensional approximation using $l_b$. Let (as in the previous section) $W^- \to B$ denote the negative eigenbundle of this non-degenerate (for large $r$) quadratic form over $B$. Recall that the trivial vector bundle of dimension $n$ was denoted $\zeta^n$ over any base. In this section we prove the following proposition.
\begin{Proposition}
  \label{prop:Quad:1}
  For $r$ odd and large enough (only depending on the bound $C_l$) the virtual vector bundle class of $W^- - \zeta^{k(r+1)}$ is classified by the map
  \begin{align} \label{eq:51}
    B \xrightarrow{l} \Lambda \La(k) \to \Lambda \La \simeq \Lambda \LUO{} \xrightarrow{\pi_\Omega} \Omega \LUO{} \simeq Z\times BO.
  \end{align}
\end{Proposition}
We will assume throughout this section that $r$ is odd. We will \emph{not} assume that $B$ is compact, but since the energy is bounded by some constant $C_l$ we can \emph{almost} assume that $B$ is compact. Indeed, the space of loops with energy less than $C_l$ in $\La(k)$ deformation retracts onto a finite CW complex, we will refer to this as: ``the fact that $l$ has compact homotopy type''. Firstly, we start by explaining the maps in this composition. 

The infinite Lagrangian Grassmannian is defined as $\La = \lim_{k\to \infty} \La(k)$, here the maps $\La(k)\to \La(k+1)$ are given by adding the trivial Lagrangian $\R\subset \C$ in a new factor. We will refer to this as a \emph{standard stabilization}. The free loop of the inclusion $\La(k) \subset \La$ is the first map after $l$ in Equation~\eqref{eq:51}. The next map, the homotopy equivalence $\La \simeq \LUO{}$, is classical (in fact $\LUO{(k)}\simeq \La(k)$) and can be found in e.g. \cite{MR1698616}. By Bott periodicity we have $\La \simeq \LUO{} \simeq \Omega^6 \Or$ (see e.g. \cite{MR0163331}). Since this is a loop space we have a homotopy equivalence
\begin{align} \label{homequi}
  \ev_0 \times \pi_\Omega \colon \Lambda \La \to \La \times \Omega \La,
\end{align}
where $\ev_0$ is evaluation at the base point, and $\pi_\Omega$ is homotopic to point-wise multiplication with the homotopy loop-inverse of
$\ev_0$. This is $\pi_\Omega$ in Equation~\eqref{eq:51}. The same Bott periodicity as above shows that
\begin{align*}
  \Omega \LUO{} \simeq \Omega^7 O \simeq \Z \times BO,
\end{align*}
but we will discuss this homotopy equivalence more explicitly below.

First we reduce the computation to an easier to understand family of quadratic forms. So, define
\begin{align*}
  Q_r^b \co (\R^{2k})^r \to \R
\end{align*}
as the quadratic form defined similarly to $S_r^b$, but using the zero Hamiltonian $H_D=0$. Note that the argument in Remark~\ref{quadrem} shows why this is a quadratic form. Indeed, for any $H_D$ which is a quadratic form on $\R^{2k}$ these finite dimensional approximations are quadratic forms.

\begin{Lemma}
  \label{lem:Quad:4}
  Let $K$ be as in Proposition~\ref{energy} (for $K=K(k)$ associated to $M=D^{2k}\subset \R^{2k}$). For $r>2KC_l^2$ the quadratic forms $Q_r^b$ has kernel given by the constant loops $\arz \in (\R^{2k})^r$. I.e. those $\arz$ where $z_{j+1}=z_j$ for all $j\in \Z/r$.
\end{Lemma}

We could prove this lemma using simple linear algebra, and even get a more explicit and better bound. However, this is \emph{rather} cumbersome, and the work we have already done is extremely helpful.

\begin{proof}
  Use Proposition~\ref{energy} for $Q_r^b$ with $M=D^{2k}$ (and $K=K(k)$), $\C_H=1$, $C_\SLa=C_l$, and $\no{H}=0$ we get for $r>K(0+C_l^2(1+1))$ and $\arz \in \Lrb D^{2k}$ that
  \begin{align*}
    \norm{\nabla Q_r^b}^2 \geq E(\arz).
  \end{align*} 
  This means that for $\arz$ close to $0$ (energy less than $1$) only the periodic orbits (constant loops) are critical points. However, since this is true close to zero it means that it is true everywhere (since it is a quadratic form).
\end{proof}

This allows us to define $V^- \to B$ as the vector bundle with fiber the negative eigenbundle of $Q^b_r$. Indeed, since the kernel is of constant (in $b$) dimension this makes perfect sense.

\begin{Lemma} \label{quadcon}
  With $K$ as above and $r>K(2+3C_l^2)$ there is a canonical (contractible choice) isomorphism
  \begin{align*}
    W^- \cong V^- \oplus \zeta^{2k}
  \end{align*}
  of real vector bundles over $B$.
\end{Lemma}

\begin{proof}
  For each $b\in B$ define a continuous family $A_c^b$ of quadratic forms for $c \in I$ by finite dimensional approximation using the Hamiltonians $H^c = c\norm{z_2}^2$. Then $Q_r^b=A_0^b$ and $S_r^b=A_1^b$. The argument in Remark~\ref{quadrem} shows why all of these are still quadratic forms.

  Applying again Proposition~\ref{energy} (as above) to $A_c^b$ for each $c$ with $C_H=2$, $C_\SLa=C_l$ and $\no{H^c}=2c$ we get for $r>K(2+C_l^2(1+2))\geq K(2c+C_l^2(1+2))$ that for all $c$ the kernel of the quadratic form consists precisely of the periodic orbits. Since the Hamiltonian flow for $H^c$ has the only the trivial 1-periodic orbit $0$ for $0<c<2\pi$, these are all non-degenerate. So by smoothness in $b$ and $c$ (and a parallel transport argument) they have isomorphic negative eigenbundles, but we need to see what happens at $c=0$. The critical points of $Q_r^b=A_0^b$ are precisely the constant curves, so the Hessian is degenerate, and the kernel as a bundle over $B$ is the trivial bundle of dimension $2k$. We prove the lemma by proving that for a small perturbation of $c=0$ in positive direction, this kernel becomes part of the negative eigenspace.

  We do this point wise in $b$. So, fix $b\in B$. Denote by $V_-, V_0$ and $V_+$ the negative, zero and positive eigenspace of $A_0^b$. It is enough to prove that the first order change in $c$ at $c=0$ of $A_c^b$ is negative definite on the kernel $V_0$. Indeed, if so we can restrict $A_0^b$ to the sphere of $V_0\oplus V_-$ and what we see is a non-positive function on a closed manifold, which is then perturbed to the first order to be negative on the set where it is zero. This will imply that the function is in fact going to be negative on the entire sphere for very small $c$, and thus $A_c^b$ is negative definite on $V_0\oplus V_-$ for small $c>0$.

  To prove this negativity to the first order in $c$ on $V_0$, we look at $A_c^b$ on $V_0$ for $c$ close to zero. The kernel $V_0$ is the set of constant curves, so we assume that $z_j=z_{j+1}$ for all $j\in \Z/r$. We need to take a look at the precise definition of
  \begin{align*}
    A_c^b = \sum_j(\int_{\gamma_j} \lambda_0 - H^cdt + \int_{\gamma^{\llcorner}_j} \lambda_0)
  \end{align*}
  For $c=0$ all of this is zero (on $V_0$) because $\gamma_j$ and $\gamma^{\llcorner}_j$ are constant, and $H^c$ is zero. We want to prove that the dominating term when perturbing to positive $c$ is $-H^c$, which is negative.

  Because any time independent Hamiltonian is constant on its flow curves, we can rewrite this as
  \begin{align*}
    A_c^b & = \int_{\sum_j(\gamma_j+\gamma^{\llcorner}_j)} \lambda_0 + \frac{1}{r}\sum_j H^c(z_j) = \\
    & = \int_{\sum_j(\gamma_j+\gamma^{\llcorner}_j)} \lambda_0 - H^c(z_0).
  \end{align*}
  The curves $\gamma_j$ are the $1/r$ time flow curves of $H^c$, so they have lengths of order $\norm{\nabla H^c}/r$ which is of order $c\norm{z_0}/r$, and since $z_j=z_{j+1}$, and $\gamma^{\llcorner}_j$ connects the endpoint of $\gamma_j$ with $z_{j+1}$, the same is true for $\gamma^{\llcorner}_j$. This means that the integral, which is minus\footnote{we are integrating $\lambda_0=ydx$ hence we get \emph{minus} the symplectic area} the symplectic area enclosed by the closed curve obtained by concatenating $\gamma_j$ and $\gamma^{\llcorner}_j$, is of order $(c\norm{z_0}/r)^2$. We have $r$ of these terms summed, but this is still of order $(c\norm{z_0})^2/r$. The term $H^c(z_0)$ is equal to $c\norm{z_0}^2$, so this is the dominating term (for small $c$) and the lemma follows.
\end{proof}

\subsection{Special cases}

We will in the following compute natural representatives for the negative eigenbundles in some special cases of $Q_r^b$. So, in this subsection we will assume that $B=\{b_0\}$ and that $k=1$. So we are only considering a single loop of Lagrangians $l=l_{b_0} \in \Lambda \La(1) \cong \Lambda S^1$, which for each $t\in I$ is defined by an argument $l(t) \in \R/\pi$ ($\R P^1$). This means that $e^{2il(t)} \in \C$. We also denote the associated quadratic form from the lemma above simply by $Q_r=Q_r^{b_0}$.

\begin{Lemma}
  \label{lem:Quad:33}
  The quadratic form $Q_r$ is given by
  \begin{align*}
    Q_r(\arz) = \frac12 \sum_j  (y_{j+1}+y_j)(x_{j+1}-x_j) - \frac14 \sum_j \im(e^{-2il(j/r)}(z_{j+1}-z_j)^2)),
  \end{align*}
  and the first sum is minus the symplectic area inside the loop given by connecting the $z_j$ in order by straight lines. The second sum is the sum of the differences of symplectic area of the straight line connection and the L-curve connecting the two points $z_j$ and $z_{j+1}$.
\end{Lemma}

In the following we will denote the two sums by 
\begin{align*}
  A=\sum_j (y_{j+1}+y_j)(x_{j+1}-x_j) \qquad \textrm{and} \qquad T=\sum_j \im(e^{-2il(j/r)}(z_{j+1}-z_j)^2).
\end{align*}
$A$ for area and $T$ for triangle area.

\begin{proof}
  Since there is no Hamiltonian term in the definition of $Q_r$ it is equal to minus the symplectic area bounded by the concatenation of the L-curves from $z_j$ to $z_{j+1}$ defining a zig-zaggy loop in $\R^{2k}$.

  The first sum is easy since it is the sum of the integration of the standard Liouville 1-form $\lambda_0=ydx$ over the straight line connections.

  The second part: pick either of the two numbers representing $l(j/r)$ in $[0,2\pi[$. Then the area of the $j$th triangle can be computed as:
  \begin{align*}
    -\tfrac12 \re(e^{-i\pi l(j/r)}(z_{j+1}-z_j))\im(e^{-\pi l(j/r)}(z_{j+1}-z_j))
  \end{align*}
  Indeed, multiplying $e^{-i\pi l(j/r)}$ onto $z_{j+1}-z_j$ simply rotates the vector into a position where it looks like the Lagrangian $l(j/r)$ equals the real axis, and in this case this product of real part and imaginary part computes the symplectic area of the triangle. This formula is the same as minus 1 fourth of the imaginary part of the square ($\re(a)\im(a)=\tfrac12 \im(a^2)$), which even removes the ambiguity of the choice of representative for $l(j/r)$.     
\end{proof}

We will need to consider a finite dimensional versions of Fourier coefficients. Indeed, let $\rho = e^{i2\pi/r}$ be the standard $r$'th root of unity. Use this to define the complex vector spaces $E_m$ by
\begin{align*}
  E_m = \{ (b\rho^{mj})_{j\in \Z/r} \mid b\in \C \} \subset \C^{r} = (\R^{2})^r,
\end{align*}
for any $m\in \Z/r$. With this we have
\begin{align}
  \label{eq:63}
  \C^r = \bigoplus_{m\in \Z/r} E_m
\end{align}
For $\arz=(z_j)_{j\in\Z/r} \in E_{m}$ and $\arw=(w_j)_{j\in\Z/r} \in E_{m'}$ one may readily check that
\begin{align}
  (z_{j+1})_{j\in\Z/r} &\in E_m \label{shift}\\
  (\ovl{z_j})_{j\in\Z/r} &\in E_{-m} \notag\\
  \re(z_j)_{j\in\Z/r} &\in E_m\oplus E_{-m} \notag \\
  \im(z_j)_{j\in\Z/r} &\in E_m\oplus E_{-m} \notag \\ 
  (z_j\cdot w_j)_{j\in\Z/r} &\in E_{m+m'}. \notag
\end{align}
and if $m\neq 0$ 
\begin{align*}
  \sum_j z_j &= 0.
\end{align*}
Notice that the second to last fact makes sense only because we have $k=1$.

\begin{Lemma}
  \label{lem:Quad:5}
  The sum $A$ splits orthogonally on the decomposition from Equation~\eqref{eq:63}. Furthermore, for $\arz = (b\rho^{mj})_{j\in \Z/r} \in E_m$ we have that
  \begin{align*}
    \tfrac12 A(\arz) = -r\norm{b}^2 \sin(2\pi m /r).
  \end{align*}
\end{Lemma}

\begin{proof}
  First we assume $\arz\in E_m$ and $\arw\in E_{m'}$ with $m\neq \pm m'$ then
  \begin{align*}
    A(\arz+\arw) = &\frac12 \sum_j (y^w_{j+1}+y^z_{j+1}+y^w_j+y^z_j)(x^w_{j+1}+x^z_{j+1}-x^w_j-x^z_j) = \\
    = & \frac12 \sum_j (y^w_{j+1}+y^w_j)(x^w_{j+1}-x^w_j) + \frac12 \sum_j (y^z_{j+1}+y^z_j)(x^z_{j+1}-x^z_j)
  \end{align*}
  since the rules above implies that summing mixed terms (in $w$ and $z$) like e.g. the term $\sum_j(x^w_{j+1}y_j^z)$ is $0$. Indeed, the products of the real part of something in $E_m$ with the imaginary part of something in $E_{m'}$ has components in $E_{\pm(m\pm m')}$, but no other $E_n$ - especially not $E_0$. Hence the sum is zero.

  The case $m=-m'$ is a little more tricky. However, for $\arz =(b\rho^{mj})_{j\in \Z/r} \in E_m$ and $\arw = (b'\rho^{-mj})_{j\in \Z/r} \in E_{-m}$ unit vectors we have that the points $(z_j+w_j)$ are all contained in the real 1-dimensional vector space spanned by $(b+b')$ in $\C$ - hence the enclosed area is zero. Now, if we establish the second part of the lemma it will follow that for these vectors:
  \begin{align*}
    A(\arz+\arw) = 0 = A(\arz)+A(\arw)
  \end{align*}
  since the formula for $A$ proves that $A(\arw)=-A(\arz)$. Having this for all unit vectors in the two subspaces proves orthogonality.

  For the formula we use that $\tfrac12 A$ is minus the symplectic area bounded by connecting the points $z_j$ to $z_{j+1}$ by straight lines. For $\arz\in E_m$ this is the formula given since each piece forms a triangle with $0$ with signed symplectic area $\norm{b}^2\sin(2\pi m/r)$. So, in a sense this is the formula you get if you integrate the other standard primitive for $\omega_0$ on $\R^{2}$ given by $\tfrac12 ydx-\tfrac12 xdy$ over this closed curve.
\end{proof}

Now we will use these linear subspaces to identify natural choices of negative eigenbundles in 3 very important cases. Define
\begin{align*}
  E_- = \bigoplus_{j=1}^{(r-1)/2} E_j \qquad \textrm{and} \qquad E_+=\bigoplus_{j=(r+1)/2}^{r-1} E_j.
\end{align*}
This gives the splitting
\begin{align} \label{eq:61}
  \C^r = E_- \oplus E_0 \oplus E_+.
\end{align}
Lemma~\ref{lem:Quad:5} tells us that this is in fact the splitting into negative, zero, and positive eigenspaces of $A$.

\begin{Lemma}
  \label{lem:Quad:1}
  Let $r>2$ be odd, $k=1$, $B=\{b_0\}$, and $l$ be a constant path of Lagrangians. Then
  \begin{align*}
    E_- \oplus E_0 \oplus E_+
  \end{align*}
  is a splitting into negative, zero, and positive vector spaces of $Q_r$.
\end{Lemma}

Notice that we say \emph{a} splitting and \emph{vector} spaces instead of eigenspaces. Indeed, we are not claiming that these are sums of eigenspaces - only that $Q_r$ restricted to each is negative definite, zero, and positive definite respectively. This, however, implies that there are canonical isomorphisms to the eigenspaces by taking orthogonal projections.

\begin{proof}
  Since $L$ is constant we can rotate and assume $L=\R$. Indeed, $Q_r$ is preserved and the splitting in Equation~\eqref{eq:61} is preserved under rotations.

  In this case Lemma~\ref{lem:Quad:33} provides
  \begin{align*}
    Q_r(\arz) = \frac12 \sum_j (y_{j+1}+y_j)(x_{j+1}-x_j) - \frac14 \sum_j \im((z_{j+1}-z_j)^2))    
  \end{align*}
  Now if we assume that $\arz \in E_-$ then by the rules in Equation~\eqref{shift} above we have $(z_{j+1}-z_j)_{j\in \Z/r} \in E_-$ and then $(z_{j+1}-z_j)^2 \in E_- \oplus E_+$. So, we avoid $E_0$. this is by the fourth rule preserved by taking imaginary part, and hence by the last rule we have that the sum is actually $0$. Hence restricting $Q_r$ to $E_-$ we have
  \begin{align*}
    Q_r(\arz) = \frac12 \sum_j (y_{j+1}+y_j)(x_{j+1}-x_j) = \tfrac12 A(\arz).
  \end{align*}
  The formula in Lemma~\ref{lem:Quad:5} for this is negative for each $m=1,\dots,(r-1)/2$ - and hence $Q_r$ is negative on $E_-$.

  The same argument on $E_+$ shows that $Q_r$ is positive definite on $E_+$ and $E_0$ is part of the kernel since translating $\arz$ by a $c\in \C$ preserves $Q_r$ - hence $E_0$ consists of critical points for $Q_r$.
\end{proof}

Now we need to modify these spaces a little bit for the next case. Indeed, let $\sqrt[n]{\rho}=e^{-2 i \pi/(nr)}$ then inside $E_{(r+1)/2}$ we have the real 1 dimensional subspace
\begin{align*}
  R = \{(b\rho^{(r-1)/2j})_{j\in \Z/r} \mid b\sqrt[4]{\rho} \in (1-i)\R \}
\end{align*}
and its orthogonal complement inside $E_{(r+1)/2}$ is
\begin{align*}
  R^{\perp} = \{(b\rho^{(r-1)/2j})_{j\in \Z/r} \mid b\sqrt[4]{\rho} \in (1+i)\R \}.
\end{align*}
Both are real lines in $\C^r$. To ease notation we will write
\begin{align*}
  E_+ \ominus R = \pare*{\bigoplus_{m=(r+3)/2}^{r-1} E_m} \oplus R^{\perp} \subset E_+
\end{align*}
with codimension 1 in $E_+$. We thus have a new splitting of $\C^r$ as
\begin{align*}
  (E_-\oplus R) \oplus E_0 \oplus (E_+ \ominus R),
\end{align*}
which has ``moved'' the line $R$ from the $E_+$ part to the $E_-$ part.

\begin{Lemma}
  \label{lem:Quad:2}
  Let $r>2$ be odd, $k=1$, $B=\{b_0\}$, and $l$ be the Maslov index 1 loop of Lagrangians defined by $l(t)=e^{i\pi t}\R$. Then 
  \begin{align*}
    (E_-\oplus R) \oplus E_0 \oplus (E_+ \ominus R),
  \end{align*}
  is a splitting into negative, zero, and positive vector spaces of $Q_r$.
\end{Lemma}

\begin{proof}
  Lemma~\ref{lem:Quad:33} gives us an explicit formula for $Q_r$
  \begin{align*}
    Q_r(\arz) = \frac12 \sum_j (y_{j+1}+y_j)(x_{j+1}-x_j) - \frac14 \sum_j \im(\rho^{-j}(z_{j+1}-z_j)^2)) = \tfrac12 A -\tfrac14 B.
  \end{align*}
  We can no longer argue that the last part vanishes on $E_+$. However, since multiplying with $(\rho^{-j})_{j\in \Z/r} \in E_{-1}$ moves us from $E_m$ to $E_{m-1}$ we get something very close. Since the symmetry is broken we get different cases when dealing with $E_-$ and $E_+$.
  
  Claim: the quadratic form $T$ splits orthogonally on $E_{(r+1)/2} \oplus \cdots \oplus E_{r-1} = E_+$ and is zero on all factors except $E_{(r+1)/2}$. To see this, let $\arw \in E_m \subset E_+$ and $\arz \in E_{m'} \subset E_+$ then we have
  \begin{align*}
    T(\arw+\arz) = &\sum_j \im(\rho^{-j}(z_{j+1}+w_{j+1}-z_j-w_j)^2) = \\
    = & \sum_j \im(\rho^{-j}(z_{j+1}-z_j)^2) + \sum_j \im(\rho^{-j}(w_{j+1}-w_j)^2)
  \end{align*}
  unless $m=m'=(r+1)/2$. Indeed, for any other pair $(m,m')$ the mixed terms before taking imaginary part (in $z$ and $w$) like e.g. $\sum_j(\rho^{-j}z_{j+1}w_j)$ are by the rules of the spaces $E_m$ zero since $(z_{j+1}w_j)$ lies in $E_k$ for $k=2,\dots,r-2$ and hence $(\rho^{-j}z_{j+1}w_j)$ lies in $E_k$ for $k=1,\dots,r-2$. However for $m=m'=(r+1)/2$ we have $(z_{j+1}w_j)$ lying in $E_1$ and hence $(\rho^{-j}z_{j+1}w_j)\in E_0$ and the sum is no longer $0$. This shows the entire claim (for $m=m'\neq (r+1)/2$ this proves $B(2\arz)=2B(\arz)$ hence $B(\arz)=0$ since it is quadratic).

  Now we compute $T$ on $E_+$, where it is in fact non-zero. Indeed, for $\arz=\sum_{m=(r+1)/2}^{r-1}\alpha_m\rho^{mj} \in E_+$ we have since only the $E_{(r+1)/2}$ part contributes that
  \begin{align*}
    B(\arz) = & \sum_j\im(\alpha_{(r+1)/2}^2 \rho^{-j}\rho^{(r+1)j}(\rho^{(r+1)/2}-1)^2) = \\
    = &\sum_j\im(\alpha_{(r+1)/2}^2(\rho^{(r+1)/2}-1)^2) = \\
    = & r\im(\alpha_{(r+1)/2}^2(-(1+\re(\sqrt[4]{\rho}))\sqrt[4]{\rho})^2) = \\ 
    = & r(1+\re(\sqrt[4]{\rho}))^2\im((\alpha_{(r+1)/2}\sqrt[4]{\rho})^2) = \\
    = & r(1+\re(\sqrt[4]{\rho}))^2\im((\alpha_{(r+1)/2}\sqrt[4]{\rho})^2) = \\
    = & r(1+\re(\sqrt[4]{\rho}))^2(\pi_{R^\perp}(\arz)^2-\pi_R(\arz)^2).
  \end{align*}
  This was the reason for the definition of $R$. Indeed, $R$ is the negative eigenvector in $E_{(r+1)/2}$ of $T$. Combined with Lemma~\ref{lem:Quad:5} this imply
  \begin{align*}
    Q_r(\arz) = -r\smashoperator{\sum_{m=(r+1)/2}^{r-1}}\norm{\alpha_m}^2 \sin(2\pi m/r) + r(1+\re(\sqrt[4]{\rho}))^2(\pi_{R^\perp}(\arz)^2-\pi_R(\arz)^2),
  \end{align*}
  which for $\pi_R(\arz)=0$ consist only of positive terms - hence $Q_r$ is positive on $E_+\ominus R$.

  Claim: the quadratic form $T$ splits orthogonally on $E_1\oplus \cdots E_{(r+1)/2} = E_- \oplus E_{(r+1)/2}$ and is zero on all except $E_{(r+1)/2}$. This is similar to above if $(m,m')\neq((r+1)/2,(r+1)/2)$ the mixed terms cannot have components in $E_0$. Indeed, $m+m'-1\in \{1,\dots,r-1\}$. Again combined with Lemma~\ref{lem:Quad:5} and the calculation of $A$ on $E_{(r+1)/2}$ above we get for any $\arz \in E_-\oplus E_{(r+1)/2}$ the same formula as above
  \begin{align*}
    Q_r(\arz) =  -r\smashoperator{\sum_{m=(r+1)/2}^{r-1}}\norm{\alpha_m}^2 \sin(2\pi m/r) + r(1+\re(\sqrt[4]{\rho}))^2(\pi_{R^\perp}(\arz)^2-\pi_R(\arz)^2).
  \end{align*}
  However, now we see that this is negative if $\pi_{R^\perp}(\arz)=0$. Indeed, the \emph{one} term that is in fact positive is related to the component in $R$, but the sum of all the contributions from $R$ to $Q_r$ are:
  \begin{align*}
    &-r\norm{\alpha_{(r+1)/2}}^2\sin(\pi (r+1)/r) + r(1+\re(\sqrt[4]{\rho}))^2\norm{\alpha_{(r-1)/2}}^2 = \\
    = & -r\norm{\alpha_{(r+1)/2}}^2((1+\re(\sqrt[4]{\rho}))^2+\sin(\pi (r+1)/r))  < 0
  \end{align*}
\end{proof}

Notice that even though we get the same formula for $Q_r$ in the two cases in the proof above it is not true that this formula works generally for any $\arz\in\C^r$. Indeed, there are many possible interacting terms. However, all we need is to know that the restriction is either positive or negative.

Since $\ovl{R} \subset E_-$ we may define $E_-\ominus \ovl{R}$  analogous to above. The last case is the conjugate of the second case.

\begin{Lemma}
  \label{lem:Quad:3}
  Let $r>2$ be odd, $k=1$, $B=\{b_0\}$, and $l$ be the Maslov index -1 loop of Lagrangians defined by $l(t)=e^{-i\pi t}\R$. Then 
  \begin{align*}
    (E_-\ominus \ovl{R}) \oplus E_0 \oplus (E_+\oplus \ovl{R}),
  \end{align*}
  is a splitting into negative, zero, and positive vector spaces of $Q_r$.
\end{Lemma}

\begin{proof}
  Since conjugation of $\arz$ and $l$ changes the sign on everything and
  \begin{align*}
    (E_-\ominus \ovl{R}) \oplus E_0 \oplus (E_+\oplus \ovl{R}) = \ovl{(E_+ \ominus R)} \oplus \ovl{E_0} \oplus \ovl{(E_-\oplus R)}
  \end{align*}
  this is the same as the lemma above.
\end{proof}

\subsection{The general case}

We now go back to the general parameterized case where $l \co B \to \Lambda \La(k)$ describes a family of loops, and thus a family of quadratic forms $Q_r^b, b\in B$. To be able to use the concrete computations in the previous subsection we will start by arguing that $l$ is homotopic to another map on a standard form after stabilizing.

\begin{Lemma}
  \label{lem:Quad:7}
  The stabilization of $l \co B \to \La(k)$ to a map $l' \co B \to \La(k')$ using the standard inclusion $\La(k) \subset \La(k')$ for $k'>k$ does not change the virtual bundle $W^--\zeta^{k(r+1)}$ over $B$ considered in Proposition~\ref{prop:Quad:1}.
\end{Lemma}

\begin{proof}
  By Lemma~\ref{quadcon} we have $W^- - \zeta^{k(r+1)} = V^- - \zeta^{k(r-1)}$ (as virtual bundles). So if we can argue that this is unchanged by stabilization we are done. A single stabilization $\La(k) \to \La(k+1)$ is given by $l'_b(t)=l_b(t)\oplus \R \subset \C^k\times\C$, and in this case everything is completely defined coordinate wise. So, the associated quadratic forms satisfy:
  \begin{align*}
    (Q_r^b)'(\arzet,\arzto) = Q_r^b(\arzet) + Q_r(\arzto), 
  \end{align*}
  where $(Q_r^b)'$ is defined on $(\C^{k+1})^r$ using $l'$, $Q_r^b$ was the original quadratic form defined by $l$ on $(\C^k)^r$ and $Q_r$ is defined for $\arzto \in (\C^1)^r$ using the constant Lagrangian $\R \subset \C$. By Lemma~\ref{lem:Quad:1} the negative eigenspace of $Q_r$ is isomorphic to $\R^{(r-1)}$. Hence a stabilization adds a trivial bundle of dimension $\R^{(r-1)}$ to the negative eigenbundles, but since we are subtracting $\zeta^{k(r-1)}$ the increase in $k$ cancels this out.
\end{proof}

For any $D_+\subset \R^k$, $D_-\subset \R^k$ and $D_0 \subset \R^k$ pairwise orthogonal and $D_+ \oplus D_- \oplus D_0 = \R^k\subset \C^k$, we define the curve $\gamma_{(D_+,D_-,D_0)} \in \Omega \La(k)$ by 
\begin{align*}
  \gamma_{(D_+,D_-,D_0)}(t) = e^{i\pi t} D_+ \oplus e^{-i\pi t} D_- \oplus D_0 \in \La(k),
\end{align*}
for $t\in [0,1]$. The space of such curves will be denoted $\Omega^S \La(k)$ ($S$ for standard form), and is a sub-space of the based loops $\Omega \La(k)$. Over this space we have the three canonical vector bundles $D_+$, $D_-$ and $D_0$.

\begin{Lemma} \label{factor}
  Any map $l\colon B \to \Lambda \La(k)$, is after stabilization homotopic (using a homotopy with energy bounds only depending on $C_l$) to a map into the subspace
  \begin{align*}
    \La(k_1) \times \Omega^S \La(k_2) \subset \Lambda \La(k_1+k_2)
  \end{align*}
  for large enough $k_1$ and $k_2$. Furthermore, the map $\pi_\Omega$ restricted to this subspace is the projection to the second factor $\Omega^S\La(k_2) \subset \Omega \La \simeq Z\times BO$ and the virtual bundle this map classifies is $D_+ - D_-$.
\end{Lemma}

\begin{proof}
  Consider the homotopy equivalence in Equation~\eqref{homequi}. The inverse to this can be described as the limit of injective maps
  \begin{align*}
    \La(n) \times \Omega \La(n) \to \Lambda \La(2n).
  \end{align*}
  as $n$ tends to infinite. To make sense of this we intertwine the factors such that: if the copy of $\C^n$ which the Lagrangians in the first factor is a subspace in has standard basis $e_1,\dots,e_n$ and the second factor has standard basis $f_1,\dots,f_n$ then the standard basis for $\C^{2n}$ on the right hand side is $e_1,f_1,e_2,\dots,e_n,f_n$. This defines the injection and this way the maps are compatible with the standard stabilizations $\La(n) \times \Omega \La(n) \subset \La(n+1) \times \Omega \La(n+1)$ and $\Lambda \La(2n) \subset \Lambda\La(2n+2)$. So, that we can take the limit. This is, indeed, a homotopy inverse to the map in Equation~\eqref{homequi} since the product on $\La$ is induced by such direct sums, and one may rearranging factors by a homotopy since $\Gl_n(\C)$ is connected. It follows that since $l$ has compact homotopy type there is a $k_1>k$ large enough such that after stabilization into $2k_1$ we have that $l$ is homotopic to a map into
  \begin{align*}
    \La(k_1) \times \Omega \La(k_1).
  \end{align*}
  We now argue that by increasing $k_1$ to some $k_2$ we can assume that the map to the last factor lands in $\Omega^S\La(k_2)$.

  So, let $f\co B \to \Omega \La(n)$ be any given map (with bounded energy). The proof of this claim follows the standard Morse theoretic proof of Bott periodicity (see e.g. \cite{MR0163331}): multiplication with $e^{-i\pi t/2}$ on the based loops gives a homeomorphism of $\Omega \La(n)=\Omega(\La(n),\R^n,\R^n)$ to $\Omega (\La(n),\R^n, i\R^n)$. Here $\Omega(X,q,q')$ denotes paths in $X$ starting at $q$ ending at $q'$. So by abuse of notation we now consider $f$ a map into the later space. In \cite{MR0163331} part IV paragraph 24 the space of minimal geodesics for this space is computed to be (with some notational change to fit the current context)
  \begin{align*}
    \Omega^{\textrm{min}} (n)= \Omega^{\textrm{min}} (\La(n),\R^n,i\R^n) = \{ \gamma \mid \gamma(t) = e^{i\pi t/2} W \oplus e^{-i\pi t/2} W^{\perp} , W \subset \R^n \}.
  \end{align*}
  The embedding of this space into $\Omega (\La(n),\R^n,i\R^n)$ has high connectivity on the components where $\dim(W)$ and $\dim(W^{\perp})$ are both high.

  To be able to use this high connectivity statement we stabilize $f$ (in a non-trivial way) by
  \begin{align*}
    \gamma(t)=e^{i\pi t/2} \R \oplus e^{-i\pi t /2} \R \subset \C^2.
  \end{align*}
  That is, we compose with the map
  \begin{align*}
    \oplus \gamma \colon \Omega \La(n,\R^n,i\R^n) \to \Omega \La(n+2,\R^{n+2},i\R^{n+2})
  \end{align*}
  given by direct sum with $\gamma$. This increases both the dimension of $W$ and its complement by 1. So, if we do this $m>0$ times with $m$ large enough, we can assume that the map $(\oplus \gamma)^{\circ m} \circ f$ factors through $\Omega^{\textrm{min}}(n+2m)$.

  Going back with the homeomorphism (above) to our version of $\Omega\La(n)$ we see that the stabilization we did corresponds to having stabilized with
  \begin{align*}
    e^{i\pi t/2} (\gamma (t)) = e^{i\pi t} \R \oplus \R
  \end{align*}
  $k$ times. We have thus argued: after $m$ stabilizations of this type, the map is homotopic to a map which factors through the following subspace
  \begin{align*}
    e^{i\pi t/2} \Omega^{\textrm{min}}(n+2m) = \{ \gamma \mid \gamma(t)=e^{i\pi t} W
    \oplus W^{\perp} \} \subset \Omega^S \La(n+2m).
  \end{align*}
  I.e. the part of $\Omega^S\La(n+2m)$ with $D_-=0$. Now by further stabilizing with
  \begin{align*}
    \gamma_2(t) = (e^{-i\pi t} \R)^{\oplus m},
  \end{align*}
  one has in total stabilized $f$ with something homotopic to a standard stabilization. Indeed, it is easy to get the two ``twistings'' to cancel out. So we have now homotoped the map $f$, stabilized in the standard way, to a map into $\Omega^S(\La(n+3m))$.

  The fact that we can assume this homotopy to have bounded energy can be argued as follows. First compose $l$ with a deformation retraction of all curves with energy less than $C_l$ to a compact subspace. Then use the above on this compact subspace in $\Lambda \La(k)$. Now this homotopy has compact image hence bounded in energy.

  The last statement in the lemma follows from the fact that these highly connected inclusions of Grassmannians into $\Omega\La\simeq \Z \times B\Or$ used in the proof, are the standard way of identifying the stable bundle with the difference of two actual bundles.
\end{proof}
  
What we in fact proved was that the map is homotopic to a map factoring through $\Omega^S \La(n+3m)$, where $D_-$ is the trivial bundle of dimension $m$ (given by the last $m$ factors in $\R^{n+3m}$). This is well-known; indeed, any virtual bundle over a compact space can be written as the difference between an actual vector bundle and a trivial vector bundle.

\begin{proof}[Proof of Proposition~\ref{prop:Quad:1}]
  Any stabilization does by Lemma~\ref{lem:Quad:7} not change the problem and by abuse of notation we still denote such a map $l$. Stabilizing enough times we can assume by Lemma~\ref{factor} above that we have a bounded energy homotopy $l^t\co B \to \Lambda \La(k), t \in I$ from $l^0=l$ with $\im(l^1) \subset \La(k_1)\times \Omega^S\La(k_2)$. Assume the bound on energy is given by some constant $C_l'$. We can use Lemma~\ref{quadcon} (with this stabilized $k$) to conclude that for $r>K(2+3C_l')$:
  \begin{itemize}
  \item Both $W^-$ and $V^-$ are defined over $B\times \{0\}$ and Lemma~\ref{quadcon} relates them by $W^-\cong V^-\oplus \zeta^{2k}$ with $k=k_1+k_2$, and the proposition will follow if we prove the corresponding statement for $V^-$. This corresponding statement is: the virtual vector bundle $V^- -\zeta^{k(r-1)}$ defined over $B=B\times \{0\}$ is classified by the map in the proposition.
  \item The vector bundle $V^-$ are defined over $B\times I$ and since $B\times\{0\} \subset B\times I \supset B\times \{1\}$ are homotopy equivalences we need only prove the corresponding proposition for $V^-$ over $B\times \{1\}$. This by Lemma~\ref{factor} above reduces us to having to prove: as virtual vector bundle classes over $B\times \{1\}$ we have
    \begin{align*}
      V^- - \zeta^{k(r-1)} = D_+- D_-,
    \end{align*}
    where the $D_+$ and $D_-$ are the canonical bundles over $\Omega^S \La(k_2)$.
  \end{itemize}
  This statement we can prove by pasting together the 3 cases we considered in the previous subsection. Indeed, for a given $b\in B$ we pick a basis in $l^1_b(0)$ for $\C^{k_1}$ (which by the constancy in this first factor is in $l_b(t)$ for all $t\in I$). Then pick a basis in $\R^{k_2}$ for $\C^{k_2}$ such that the first $k_2^+$ vectors is a basis for $D_+=D_+(b)$ the next $k_2^0$ a basis for $D_0=D_0(b)$ and the last $k_2^-$ is basis for $D_-=D_-(b)$. So, $k_2^++k_2^0+k_2^-=k_2$. In this basis we have:
  \begin{align*}
    l^1_b(t) = \R^{k_1} \oplus e^{-i\pi t}\R^{k_2^+}\oplus \R^{k_2^0} \oplus e^{i\pi t}\R^{k_2^-}
  \end{align*}
  This means that $l^1_b$ is a product of $k_1+k_2$ copies of the cases considered in the previous subsection. The quadratic form $Q_r^b$ in this basis splits into a sum of each part (it defines an orthogonal decomposition). So, (with $E_-,E_0$ and $E_+$ in that subsection) we conclude by Lemma~\ref{lem:Quad:1}, Lemma~\ref{lem:Quad:2} and Lemma~\ref{lem:Quad:3} that
  \begin{align*}
    (\R^{k_1} \otimes E_-) \oplus (\R^{k_2^+}\otimes (E_-\oplus R)) \oplus (\R^{k_2^0}\otimes E_-) \oplus (\R^{k_2^-}\otimes (E_-\ominus \ovl{R})
  \end{align*}
  is a negative vector space for $Q_r^b$ of maximal dimension, and hence its orthogonal projection to the negative eigenspace of $Q_r^b$ is an isomorphism. Since acting by $O(k_1)\times O(k_2^+) \times O(k_2^0) \times O(k_2^-)$ in the obvious way does not change this vector space we see that it in fact is canonically defined without picking the basis.

  Had $D_+$ and $D_-$ both been $0$ then this vector space would simply be $\R^{k_1+k_2} \otimes E_- \cong \R^{(k_1+k_2)(r-1)}$. However, with the vector spaces non-trivial we in fact can write the above as:
  \begin{align*}
    (\R^{k_1+k_2} \otimes E_-) \oplus (\R^{k_2^+}\otimes R) \ominus (\R^{k_2^-}\otimes \ovl{R}).
  \end{align*}
  Here $V \ominus V'$ means that $V' \subset V$ and we take the orthogonal complement of $V'$ inside $V$. In standard coordinates this can be written as
  \begin{align*}
    (\R^{k_1+k_2} \otimes E_-) \oplus (D_+ \otimes R) \ominus (D_- \otimes \ovl{R}),
  \end{align*}
  which as a virtual vector bundle over $B$ is simply $D_+-D_- + \zeta^{k(r-1)}$. Indeed, $R$ and $\ovl{R}$ are constant lines hence trivial line bundles over $B$.
\end{proof}

 
\section{The Maslov Bundle and the homotopy type of the Target}  \label{maslov}

In this section we define the virtual Maslov bundle $\eta$ over $\Lambda L$ associated to the embedding $L\subset T^*N$, and prove the following proposition.
\begin{Proposition}
  \label{prop:Homtyp:1}
  The spectrum $W$ (from Section~\ref{indexcalc}) is homotopy equivalent to
  \begin{align*}
    (\Lambda L)^{-TL+\eta}
  \end{align*}
  where $\eta$ is the virtual bundle defined below.
\end{Proposition}

Let $j\colon L\to T^*N$ be any Lagrangian embedding (immersion is enough for the definition of $\eta$.). We will define the virtual Maslov bundle relative to this embedding. It is a generalization of the Maslov index related to curves of Lagrangian subspaces in $\R^{2n}$ (see e.g. \cite{MR1698616}). In fact the bundle is a canonical virtual vector bundle over $\Lambda L$, such that the dimension of this bundle on each component is precisely the Maslov index.

The projection $T^*N \to N$ will be denoted $\pi$. For any point $q\in L$ the tangent space $T_{q}L$ is mapped by $j_*$ to a Lagrangian subspace of $T_{j(q)}(T^*N)$, and by abuse of notation we use this to define a section
\begin{align*}
  j_* \colon L \to \La(T(T^*N))_{\mid L},
\end{align*}
where $\La(T(T^*N))$ is defined in definition~\ref{ladefn}. A stabilization of this map with a vector bundle $V \to N$ will be denoted by $j_*\oplus V$ and is defined by
\begin{align*}
  (j_* \oplus V) (q) = j_*(T_{q}L) \oplus (\pi^*V) \subset T_{j(q)}(T^*N) \oplus (\pi^*V) \oplus (\pi^*V)^*,
\end{align*}
which is also Lagrangian (in the obvious symplectic structure).

Let again $\nu$ be the normal bundle of $N$ for some embedding $N \to \R^k$, we get a canonical symplectic trivialization
\begin{align*}
  T(T^*N) \oplus (\pi^*V) \oplus (\pi^*V)^* \to T^*N \times (\R^{2k},\omega_0).
\end{align*}
This is defined by using the Riemannian metric induced from the embedding to split the tangent space of $T^*N$ at $z$ into $T_{\pi(z)}N \oplus T^*_{\pi(z)}N$, then mapping $V_z=T_{\pi(z)}N \oplus V_{\pi(z)}$ isomorphically to $\R^k$ by the obvious map, and mapping $V_z^*=T^*_{\pi(z)}N \oplus V_{\pi(z)}^*$, by the inverse of the dual to this map, to $i\R^k$. If we compose this trivialization with $j_*\oplus V$ we get a map from $L$ to $\La(k)$, and since all embeddings are isotopic for $k$ sufficiently large, we have a map unique up to homotopy
\begin{align*}
  F \colon L \to \La(k).
\end{align*}
Since the inclusions $\{*\}\times \La(k) \subset \La(k) \times \La(k) \subset \La(2k)$ and $\La(k) \times \{*\} \subset \La(k) \times  \La(k) \subset \La(2k)$ are homotopic, we see that the maps to $\La(2k)$, given by $\R^k \oplus F$ and $F\oplus \R^k$ are homotopic. This implies that after enough stabilization we can homotopy this map to be trivial/horizontal (equal to the horizontal $TN$) in the tangent space of $T^*N$. If we subsequently stabilize by a copy of the bundle $TN$ we have essentially stabilized by a trivial bundle.

Since we can globally homotopy horizontal directions in $T^*N$ to verticals by the homotopy which multiplies with $e^{-tJ\pi/2}$ in $T(T^*N)$ and since the same is true close to $L\subset T^*N$ we have argued the following lemma.

\begin{Lemma}
  \label{lem:Maslov:1}
  For $DT^*L \subset DT^*N$ the canonical sections $\SLa^L,\SLa^N_{\mid DT^*L} \co DT^*L \to \La(DT^*L)$ (from Example~\ref{SLdef}) satisfy
  \begin{align*}
    \SLa^L \oplus \R^k \simeq \SLa^N_{\mid DT^*L} \oplus F \co DT^*L \to \La(TM)\times \La(k)
  \end{align*}
  for large enough $k$.
\end{Lemma}

So, $F$ measures the stable difference of these Lagrangians sections. The above discussion really tells us that any two sections in $\La(TM) \to M$ (with $M$ compact) has such a stable difference map to $\La(k)$ for large enough $k$. It will actually be more convenient to use the inverse map to $F$ in the following. Indeed, this is given by taking the conjugate of the Lagrangian (this specific description will not be important, but motivates the notation) so we denote this map $\ovl{F}$. This satisfies:
\begin{align}
  \label{eq:65}
  \SLa^N_{\mid DT^*L} \oplus \R^k \simeq \SLa^L \oplus \ovl{F} \co DT^*L \to \La(TM)\times \La(k),
\end{align}
which is more natural to use in the following.

\begin{Definition} \label{maslovbundle}
  The Maslov bundle $\eta$ is the virtual vector bundle classified by the map
  \begin{align} \label{eq:50}
    \xymatrix{
      \Lambda L \ar[r]^{\Lambda \ovl{F}\quad} & \Lambda \La(k) \ar[r] &
      \Lambda \La \ar[r]^{\pi_\Omega} & \Omega \La
      \ar[r]^{\simeq\quad} & \Z\times \BO.
    }
  \end{align}
\end{Definition}
This is the same map as we saw in Proposition~\ref{prop:Quad:1}.

\begin{proof}[Proof of Proposition~\ref{prop:Homtyp:1}]
  Consider the Hamiltonians used in Section~\ref{vitconst}. Recall the narrowing process we used in Section~\ref{indexcalc} to compute the spectrum $Z$. Using this same idea we can narrow $H^{s_l}$, but in the neighborhood $M=D_{1/2}T^*L$, and argue as in Proposition~\ref{prop:Calculation:1} that by narrowing enough we can have the Conley index completely defined inside $\Lrb D_{1/2}T^*L \subset \Lre T^*N$. However, to compare it to something defined on $T^*L$ we need to adjust the Riemannian structure. So, let $g^v$ be the convex combination Riemannian structure from the one induced from $D_{1/2}T^*L \subset T^*N$ to the one induced from $D_{1/2}T^*L \subset T^*L$. Simultaneously, let $\SLa^v$ be a homotopy of sections guaranteed by Equation~\eqref{eq:65} above. Now consider the domain:
  \begin{align*}
    P = D_{1/4}T^*L \times D^{2k} \subset T^*N \times \R^{2k}
  \end{align*}
  and the Hamiltonians as in Section~\ref{stabloc} (narrowing $H^{s_l}$ depending on a parameter $\param>0$ on the first factor and not depending on $\param$ on the second factor). The added thing here is that now we have a family of underlying structures $g^v$ and $\SLa^v$ for $v\in I$. The $K$ in Proposition~\ref{energy} can by a compactness argument be chosen to work for all $v\in I$. Similarly, we can pick the upper bound on the parameter $\param_0$ and $\beta$ in Proposition~\ref{mainstabloc} to work for all $v\in I$.

  We may further assume that: since $\SLa^v$ for $v=0$ and $v=1$ are on product form this $\param_0$ and $\beta$ works for Lemma~\ref{7} (and thus Lemma~\ref{thomsusp}) in these two cases as well. Now in the following let $r$ and $\param$ be such that all these Propositions and Lemmas applies to get good index pairs. We will need to use Proposition~\ref{prop:Quad:1} on the second factor in the product cases ($v=0$ and $v=1$), and by compactness of $DT^*L$ and $I$ we have a bound on the energy of the loops on the last factor for our entire homotopy. As in Section~\ref{stabloc} this bound can be written as $\beta C_\SLa^2$ for $C_\Sla$ a bound on the derivatives of all the Lagrangians sections in the homotopy. Now we can make sure that $r$ is also large enough and odd for Proposition~\ref{prop:Quad:1} to apply for the quadratic forms on the second factors when $v=0$ and $v=1$.

  The argument that $W(l) \simeq (\Lambda^{s_l\mu_L} L)^{-TL+\eta}$ (here $\eta$ is restricted to this finite length part of the loop space) is now divided into steps
  \begin{itemize}
  \item Consider the Conley index $I_{a_{s_l}^L}^{b_{s_l}}(S_r,X_r)$ (associated to $H^{s_l}$, which is not yet narrowed) used to define $W(l)$ (adjusted by $r+1$ copies of $\nu$ such that increasing $r$ gives actual suspensions - this adjustment can simply be carried along during all the following steps).
  \item The narrowing process (down to the narrow parameter $\param$ fixed above) proves that this is homotopy equivalent to a similar Conley index, but with index pair inside $\Lrb D_{1/2}T^*L$ ($r$ is completely fixed during all these steps).
  \item Adding a new factor of $D^{2k}$ and doing the finite reduction with the section $\SLa^{v=0}=\SLa^N\oplus \R^k$ on $P$ instead changes the homotopy type of the Conley index by a Thom-space construction (Lemma~\ref{thomsusp}) of the negative eigenbundle associated to the second factor. This bundle is trivial since the finite dimensional approximation on the second factor is constant. For this particular $r$ this is a standard $(r+1)k$ fold suspension by Proposition~\ref{prop:Quad:1} (the bundle is the trivial $\zeta^{(r+1)k}$ - this follows from Proposition~\ref{prop:Quad:1} with the map $l\co B \to \Lambda \La(k)$ constant).
  \item The homotopy of structures for $v\in I$ now by Lemma~\ref{hominv} (as usual) provides a Homotopy equivalence of this with the similar index defined on $\Lrb (D_{1/2}T^*L \times D^{2k})$, but with the structure from the inclusion $D_{1/2}T^*L \subset T^*L$, and the Lagrangian section $\SLa^{v=1} = \SLa^L\oplus \ovl{F}$. Now Lemma~\ref{thomsusp} and Proposition~\ref{prop:Quad:1} states that we get a relative Thom-space construction using a representative (more specifically $W^-$) of the virtual bundle $\eta+\zeta^{(r+1)k}$ on the index pair.
  \item Now Proposition~\ref{prop:Calculation:1} applied to $L$ in instead of $N$ with this extra Thom-space construction from the second factor proves the homotopy equivalence.
  \end{itemize}

  The identification of the maps $W(l)\to W(l+1)$ can be done completely analogously to the identification in Section~\ref{indexcalc} of the map $Z(l) \to Z(l+1)$. Indeed, the construction there can be done close to $L$ as well, such that we after the narrowing step above we identify this map with the inclusion of a Conley index (defined close to $L$) into a slightly larger one, and the above remaining steps are all easily compatible with inclusion (and quotients) to Conley indices of smaller intervals of action.
\end{proof}

Note that the identification of this spectra depends on the choice of homotopy $\SLa^v$. So, the identification might not be canonical.

\begin{proof}[Proof of Theorem~\ref{thm:2}]
  This follows from combining Proposition~\ref{prop:Viterbo:1}, Proposition~\ref{prop:Calculation:1}, and the above proposition.
\end{proof}

Again we may also consider the case in Remark~\ref{rem:Spectrum:2}, and use the notation for this alternate map of spectra $Z'\to W'$, which at each point in the limits is given by maps $Z'(l) \to W'(l)$.

\begin{Corollary}
  \label{cor:Maslov:3}
  The spectrum $W'$ appearing in Corollary~\ref{cor:Viterbo:1} satisfy
  \begin{align*}
    W'\simeq  (\Lambda L)^{TN-TL+\eta},
  \end{align*}
  and the alternate transfer map is the same on the level of homology (up to a shift) for oriented $N$.
\end{Corollary}

\begin{proof}
  This is completely analogous to the above, except the actual bundles showing up in the target is changed to this because we add on less copy of $\nu$ (and grade a little differently). The statement about the map on homology follows from naturality of the Thom-isomorphism for the oriented normal bundle $\nu$ - recall that the difference is precisely adding an extra copy of $\nu$ or not.
\end{proof}

The following corollary to the above proof of the proposition is needed in \cite{immersions}. This is a generalization of Corollary~\ref{cor:Calculation:3}.

\begin{Corollary}
  \label{cor:Immersion}
  Let $(A_r,B_r)$ be an index pair for $S_r$ with a narrow $H_l'$ as above. The inclusion $A_r \subset \Lrb D_{1/2} T^*L$ induces a map
  \begin{align*}
    A_r/B_r \to (T^* \Lre L)_+ \wedge A_r/B_r,
  \end{align*}
  which induces a spectrum map
  \begin{align*}
    W'(l) \to (T^*\Lre L)_+ \wedge W'(l),
  \end{align*}
  which is canonically (contractible choice) stably homotopic to the map
  \begin{align*}
    W'(l) \to (\Lambda^{s_l\mu_L} L)_+ \wedge W'(l) \subset (\Lambda L)_+ \wedge W'(l)
  \end{align*}
  induced by the diagonal defined using the homotopy equivalence above.
\end{Corollary}

Note here that by the diagonal we mean that: for any Thom space $X^V=(X,\varnothing)^{V/}$ the diagonal map induces the map $X^V \to X_+ \wedge X^V$, and this induces a similar map of Thom-spectra (see Appendix~\ref{cha:app} for a definition of Thom-spectra where this is easily incorporated).

\begin{proof}
  If we disregard (project away from) the second factor $D^{2k}$ the proof above shows that the map $A_r/B_r \to (T^* \Lre L)_+ \wedge A_r/B_r$ is the map we already identified in Corollary~\ref{cor:Calculation:3}, but with $N$ replaced by $L$. Indeed, the index pair on the product is in the first factor (see proof of Lemma~\ref{thomsusp}) an index pair for the first factor. We also now that the $k(r+1)$th space in the spectrum is this pair adjusted by adding extra bundles, but that does not change the map to $(\Lambda L)_+\wedge (\dots)$. Hence the map defined
  \begin{align*}
    (\Lambda^{s_l\mu_L} L )^{TN-TL+\eta} \to (\Lambda^{s_l\mu_L} L)_+ \wedge (\Lambda^{s_l\mu_L} L )^{TN-TL+\eta}
  \end{align*}
  is canonically identified (contractible choice) with the diagonal (using the above identification).
\end{proof}


\appendix
\section{}\label{cha:app}

In Section~\ref{sec:gener-funct-spectr} we defined the notion of a spectrum. To give an idea of what these are and how to think of them in the present context we relate these to CW spectra defined by Adams in \cite{MR0402720}, which we briefly describe here. We then relate this to Morse homology. We also describe a construction of Thom-spectra, which is closely related to the construction we see in the paper; and finally we discus the mapping cylinders and homotopy colimits appearing in the constructions.

\subsection{Spectra and CW-Spectra}

Let $Z=(Z_n,\sigma_n)$ be a spectrum and inductively take a based CW approximation $c_n \co Z_n' \to Z_n$ (see \cite{MR1867354}) such that it extends the previous making the diagrams
\begin{align*}
  \xymatrix{
    \Sigma Z_n \ar[r]^{\sigma_n} & Z_{n+1} \\
    \Sigma Z'_n \ar[r]^{\sigma'_n} \ar[u]^{\Sigma c_n} & Z'_{n+1} \ar[u]^{c_{n+1}}
  }
\end{align*}
commute on the nose and making $\sigma'_n$ a CW inclusion. Note that the non-base-point cells in $\Sigma Z_n'$ are 1-1 with a shift of 1 in dimension to the ones in $Z_n'$, and the suspension isomorphism on $\tH_*^{\CW}$ is given by the corresponding degree 1 shift on the chain complex level ($\tC^{\CW}_*$). Now the definition of spectrum homology in Section~\ref{sec:gener-funct-spectr} is recovered by
\begin{align*}
  H_*(Z) \cong \colim_{n\to \infty} \tH^{\CW}_{*+n}(Z_n').
\end{align*}
The CW spectrum $Z'=(Z_n',\sigma_n')$ is by construction weakly homotopy equivalent to $Z$. In \cite{MR0402720} Adams defines a category of such CW spectra, and the reason to have $\sigma_n'$ be a CW inclusion is that then one may think of this as actually having cells. Indeed, the cells are increased in dimension each time the space is suspended, and then we add new cells. Keeping track of the degrees (a cell of dimension $k$ added at level $n$ has degree $k-n$) one can recovered the homology of $Z'$ using a single cellular chain complex with one generator per cell in $Z'$. This is precisely the same as the limit chain complex
\begin{align*}
  \tC^{\CW}_*(Z') = \colim_{n\to\infty} \tC^{\CW}_{*+n}(Z'_n) = \cup_{n\in\N}\tC^{\CW}_{*+n}(Z'_n).
\end{align*}
Here the colimit turns into a union precisely because the maps we are taking the colimit of are injective on the chain complexes. Note that since colimits commute with taking homology this, indeed, does recover $H_*(Z')$.

This colimit/union idea is what Adams uses to define maps between spectra. Indeed, he defines maps as you would between colimits. This implies for example that any CW spectrum level-wise inclusion $Y' \subset Z'$ (commuting with structure maps) is an \emph{isomorphism} if and only if all cells in $Z'$ eventually appear in $Y'$. If all the cells do not appear one may use a spectrum version of Hurewitz (which follows from the usual one adapted to this setting) and conclude that it is a homotopy equivalence of spectra if and only if it induces an isomorphism when passing to homology. However, we have omitted the general definition of maps in this paper since all the maps we construct are actually constructed at some level $n$. Except in the case of the constructed homotopy colimits in Equation~\eqref{eq:62} and Equation~\eqref{eq:64}, which we describe explicitly below.

Two of the most important operations on spectra is wedging $\vee$ and smashing $\wedge$. The wedge $\vee$ is easy - you simply take the wedge level-wise and use $\Sigma (X\vee Y) = (\Sigma X) \vee (\Sigma Y)$, and one can consider this as the spectrum version of disjoint union since the parts never touch \emph{except} at the base-point (we really only consider the non base-point cells as cells). The smash product is much more subtle and requires a lot of structure to define properly and we omit it here - the subtleties are related to the reordering of suspension factors mentioned in the proof of Proposition~\ref{prop:spectrum}.

\subsection{Relation to Morse homology}

The CW spectrum view-point is particularly good when relating to Morse theory. Indeed, the Conley indices $I_a^b(S_r,X_r)$ used to define the spectra in the paper can (in the Morse setting) be CW approximated by using a single cell per critical point. However, since the action on the infinite dimensional manifold of loops in $T^*N$ has infinite Morse indices it is only natural that the dimension of these cells goes up as we increase the ``fineness'' (number of points $r$) of the finite dimensional approximations. However, in this case the CW homology considered above gets a single generator per critical point, which is precisely what Morse homology has. Relating the differentials of these and those in Floer homology is more than a little subtle; indeed, the two approaches are counting the same gradient trajectories, but the signs may differ. It was thought to be the same signs when $N$ is oriented, but the homotopy constructions in this paper has revealed that this is only true if $N$ is also spin (see \cite{MySympfib} for more details on this).

\subsection{Thom-spectra}

The following specific construction is formulated in the way it is used in the paper, which is why it may look a little warped compared to standard definitions. However, the reader familiar with some other construction should easily be able to relate it to this.

Let $f\colon X \to \Z \times B\Or$ be any map from a space $X=\cup_{l\in \N} X_l$ where $X_l$ is of finite homotopy type (homotopy equivalent to a finite CW complex). This could of course be $\Lambda L$ and $\Lambda^{s_l \mu_L} L$, and the map $f$ could be given by the map described in Definition~\ref{maslovbundle}, which is precisely what came up in the paper.

First note that if $X$ is not connected and $f$ has range in different components in $\Z \times B\Or$ we may simply split it up into components and wedge the resulting components of the Thom-spectra together. So in the following we assume that $f$ is a map from $X$ to $\{d\} \times B\Or$.

Then we describe how to define it in the case that all of $X$ is in fact of finite homotopy type. Indeed, in this case there exists an $n\geq 0$ such that the map $f$ is homotopic to a map $f'\co X \to \{d\} \times B\Or(n)$. We then (dependent on this homotopy) define the Thom-spectra as the shifted suspension spectrum defined by
\begin{align*}
  (X^f)_{n+d} = D(f'^*\gamma_n)/S(f'^*(\gamma_n)).
\end{align*}
That is, the $(n+d)$th space is the Thom space of the canonical bundle $\gamma_n \to B\Or(n)$ pulled back to $X$, and for all $n'\geq n+d$ 
\begin{align*}
  (X^f)_{n'} = \Sigma^{n'-n-d}(X^f)_{n+d}
\end{align*}
and the structure maps are the identity. So if the bundle is oriented the homology will (by the Thom isomorphism) be isomorphic to the homology of $X$ but shifted by the (virtual) dimension $d\in \Z$.

In the general case let $n_l$ be strictly increasing and such that $f_{\mid X_l}$ factors through $\{d\} \times B\Or(n_l)$ up to homotopy, we may inductively assume that this is compatible with the chosen homotopy of the restriction to $X_{l-1}$. In fact, we may we assume that there is a homotopy to a map $f'\co X \to \{d\}\times B\Or$ such that the image of $f'_{\mid X_l}$ is in $\{d\} \times B\Or(n_l)\subset \{d\}\times B\Or$. Now we define for each $l$ the spectrum $X^{f,l}$ as above, which means it is the degree shifted suspension spectrum of:
\begin{align*}
  (X^{f,l})_{n_l+d} = X_l^{(f_l^*\gamma_{n_l})} = D(f^*_l\gamma_{n_l})/S(f_l^*\gamma_{n_l}).
\end{align*}
There are canonical maps of spectra from this to the next. Indeed, if we consider the restriction of the Thom space construction over $X_{l+1}$ to $X_l$, which is functorial, we see that we have a canonical inclusion
\begin{align*}
  \Sigma^{n_{l+1}-n_l}(X^{f,l})_{n_l+s} = & \Sigma^{n_{l+1}-n_l}D(f^*_l\gamma_{n_l})/S(f_l^*\gamma_{n_l}) = \\
  = &  D(f^*_l\gamma_{n_{l+1}})/S(f_l^*\gamma_{n_{l+1}}) \subset (X^{f,l}).
\end{align*}
Se Section~\ref{sec:gener-funct-spectr} for an explanation of the second equality. This means that we have canonical inclusions of spectra
\begin{align*}
  X^{f,l} \subset X^{f,l+1}
\end{align*}
and $X^f$ is ``simply'' the union of all these. However, in Equation~\eqref{eq:64} where a similar sequence turned up we used a homotopy colimit instead of a union. This was because we had similar maps up to homotopy, but since they were constructed as quotients of Conley index pairs they were not inclusions. However, by taking \emph{homotopy} colimits we replace them by inclusions. The next sections explains some aspects of this construction.

\subsection{Limits of spectra and mapping cylinders}

Assume we have a sequence of spectra $X^l \subset X^{l+1}$ where each inclusion is a level-wise cofibration. Then we can take their limit as simply the union $X=\cup_l X^l$. However, what appears in the paper is a sequence of maps $f^l \co Z^l \to Z^{l+1}$ (of spectra) which are not cofibrations, but where we have diagrams
\begin{align*}
  \xymatrix{
    Z^l \ar[r]^{f^l} \ar[d] & Z^{l+1} \ar[d] \\
    X^l \ar[r]^{g^l} & X^{l+1}
  } 
\end{align*}
where the vertical maps are homotopy equivalences (sometimes given by a contractible choice) and the diagram homotopy commutes (again sometimes given by a contractible choice). Then if we want to take a limit of the maps $Z^l \to Z^{l+1}$ which has a (contractible choice in the case where the above were such) homotopy equivalence to $X$ we take a homotopy colimit. This we define by letting $Z'^l$ be the mapping telescope of $Z^0 \xrightarrow{f^0} Z^1 \xrightarrow{f^1} \cdots \xrightarrow{f^{l-1}} Z^l$. Now we define
\begin{align*}
  \hocolim_{l\to \infty} Z^l = \cup_{l} Z'^1
\end{align*}
since $Z'^l \subset Z'^{l+1}$ is a cofibration. This is (contractible choice) homotopy equivalent to $X$ by the fact that $\hocolim_{l\to \infty} X^l$ deformation retracts onto $X$, and having the above diagrams means one can construct a map (contractible choice in that case) between mapping cylinders
\begin{align*}
  Z'^l \to X'^l,
\end{align*}
which are homotopy equivalences, and compatible with the inclusions $Z'^l \subset Z'^{l+1}$ and $X'^l \subset X'^{l+1}$.

There is a slight subtlety about these homotopy colimits for readers not very familiar with spectra. Indeed, above we noted that maps of spectra are defined only on an equivalent ``sub-spectrum'' this means that if we fix a level $k\in \N$ and look at the mapping cylinders $Z'^l$ $k$th level $(Z'^l)_k$ then for $l$ increasing the part of $Z^{l'}_k$ inside this for fixed $l'<l$ can get smaller and smaller (even be empty at some point). However, taking their union (and considering how this looks for CW spectra) we see that all cells appear at some point. This also means that strictly speaking the ``inclusion'' of $Z'^l$ into $Z'^{l+1}$ is not actual a level-wise inclusion. It is an inclusion of an isomorphic ``sub-spectrum'' of $Z'^l$ into $Z'^{l+1}$, but that is just as good as an inclusion in the category of spectra.

Combining this subtlety with the fact that our level-spaces in Section~\ref{sec:gener-funct-spectr} were already defined as mapping cylinders does not make this less confusing. This however, does not change anything in the above.


\bibliographystyle{../../hplain}
\bibliography{../../Mybib}

\end{document}